\documentclass[3p,10pt]{elsarticle}

\usepackage{amssymb,amsthm}
\usepackage{latexsym,amsmath,subfigure,tikz,siunitx,lineno}

\journal{}

\newcommand{\eps}{\varepsilon}
\newcommand{\set}[1]{\{#1\}}

\newcommand{\epsb}{\varepsilon_\mathrm{b}}
\newcommand{\mub}{\mu_\mathrm{b}}

\newcommand{\p}{\partial}
\newcommand{\mH}{\mathbf{H}}
\newcommand{\mG}{\mathbf{G}}
\newcommand{\mB}{\mathbb{B}}
\newcommand{\mU}{\mathbf{U}}
\newcommand{\mV}{\mathbf{V}}
\newcommand{\mW}{\mathbf{W}}
\newcommand{\me}{\mathbf{e}}
\newcommand{\mf}{\mathbf{f}}
\newcommand{\mg}{\mathbf{g}}
\newcommand{\mr}{\mathbf{r}}
\newcommand{\vt}{\boldsymbol{\theta}}
\newcommand{\vv}{\boldsymbol{\vartheta}}
\newcommand{\vx}{\boldsymbol{\xi}}
\newcommand{\vn}{\boldsymbol{\nu}}

\theoremstyle{plain}
\newtheorem{thm}{Theorem}[section]
\newtheorem{lem}[thm]{Lemma}

\theoremstyle{remark}
\newtheorem{rem}{Remark}[section]
\newtheorem{ex}{Example}[section]
\newtheorem{discuss}{Discussion}[section]

\DeclareMathOperator*{\music}{MUSIC}
\DeclareMathOperator*{\inc}{inc}
\DeclareMathOperator*{\obs}{obs}
\DeclareMathOperator*{\scat}{scat}

\DeclareMathOperator*{\noise}{noise}
\DeclareMathOperator*{\range}{Range}

\begin{document}

\begin{frontmatter}



\title{A novel study on the MUSIC-type imaging of small electromagnetic inhomogeneities in the limited-aperture inverse scattering problem}
\author{Won-Kwang Park}
\ead{parkwk@kookmin.ac.kr}
\address{Department of Information Security, Cryptography, and Mathematics, Kookmin University, Seoul, 02707, Republic of Korea.}

\begin{abstract}
We apply MUltiple SIgnal Classification (MUSIC) algorithm for the location reconstruction of a set of {two-dimensional circle-like} small inhomogeneities in the limited-aperture inverse scattering problem. Compared with the full- or limited-view inverse scattering problem, the collected multi-static response (MSR) matrix is no more symmetric (thus not Hermitian), and therefore, it is difficult to define the projection operator onto the noise subspace through the traditional approach. With the help of an asymptotic expansion formula in the presence of small inhomogeneities and the structure of the MSR-matrix singular vector associated with nonzero singular values, we define an alternative projection operator onto the noise subspace and the corresponding MUSIC imaging function. To demonstrate the feasibility of the designed MUSIC, we show that the imaging function can be expressed by an infinite series of integer-order Bessel functions of the first kind and the range of incident and observation directions. Furthermore, we identify that the main factors of the imaging function for the permittivity and permeability contrast cases are the Bessel function of order zero and one, respectively. This further implies that the imaging performance significantly depends on the range of incident and observation directions; peaks of large magnitudes appear at the location of inhomogeneities for permittivity contrast case, and for the permeability contrast case, peaks of large magnitudes appear at the location of inhomogeneities when the range of such directions are narrow, while two peaks of large magnitudes appear in the neighborhood of the location of inhomogeneities when the range is wide enough. The numerical simulation results via noise-corrupted synthetic data also show that the designed MUSIC algorithm can address both permittivity and permeability contrast cases.
\end{abstract}


\end{frontmatter}




\section{Introduction}
In this paper, we consider the fast identification of a set of two-dimensional small inhomogeneities, the dielectric permittivities (or magnetic permeabilities) of which differ from those of homogeneous space in the limited-aperture inverse scattering problem. This is an old but very interesting problem for scientists and engineers because it has a potentially wide application range in physics, geophysics, nondestructive evaluations, material engineering, and medical sciences, e.g., biomedical imaging \cite{A1,A3}, stroke detection \cite{PFTYMPKE,SVMM}, ground-penetrating radar \cite{A6,LSL}, synthetic aperture radar (SAR) imaging \cite{CB,JCPKJC}, breast cancer detection \cite{HSM1,RSAAP}, damage detection of civil structures \cite{CPL,KJFF}, and landmine detection \cite{CGSMMRW,DEKPS}. In essence, the main aim of the limited-aperture inverse scattering problem is the identification of unknown shape, location, or physical properties, which involves estimating the electric conductivity, dielectric permittivity, or magnetic permeability.

To solve the problem, various iterative-based schemes have been investigated, e.g., Newton-type methods \cite{K,K3}, Gauss-Newton methods \cite{ASKK,CDLR}, level-set techniques \cite{AKM,DL}, the Born iterative method \cite{IBA,L7}, conjugate gradient method \cite{BSSSNST,HWTT}, Levenberg-Marquardt algorithm \cite{BDK,CM2}, and optimization schemes \cite{ASSS,FMP}. As we have already seen \cite{KSY,PL4}, the iteration process must begin with a good initial guess that is close to the unknown target. If not, serious situations will arise, e.g., non-convergence issues, local minimizing problems, and high computational costs, which means that the success of iterative-based schemes is significantly dependent on the generation of a good initial guess. Accordingly, it is natural to develop both a mathematical theory and a numerical technique for generating a good initial guess and, correspondingly, various non-iterative techniques for identifying the location or imaging the shape of unknown inhomogeneities.

The multiple signal classification (MUSIC) algorithm is a well-known non-iterative technique in the inverse scattering problem. Traditionally, MUSIC is used in signal processing problems to estimate the individual frequencies of multiple time-harmonic signals \cite{C}, and in pioneering research \cite{D}, it has been used to identify the locations of a number of point-like scatterers. MUSIC has also been applied in various inverse scattering problems, e.g., the identification of three-dimensional small inhomogeneities \cite{AILP}, localization of small inhomogeneities hidden in three-dimensional half-space \cite{IGLP}, detection of internal corrosion \cite{AKKLV}, imaging of thin curve-like inhomogeneities \cite{P-MUSIC1}, perfectly conducting cracks \cite{PL1}, extended targets \cite{HSZ1}, and anisotropic inhomogeneities \cite{ZC}. Based on this utility, it has been applied to various real-world problems, such as damage diagnosis on complex aircraft structures \cite{BYG}, eddy-current nondestructive evaluation \cite{HLD}, impedance tomography \cite{H3}, remote sensing in safety/security applications \cite{KCC}, rebar detection \cite{SLA}, damage imaging of aircraft structures \cite{FZSY}, indoor localization problems \cite{MISH}, through-wall imaging \cite{TBAOU}, super-resolution fluorescence microscopy for single-molecule localization \cite{AM2}, multi-frequency imaging \cite{MNPT}, time-reversal MUSIC for the imaging of extended targets \cite{LH2}, imaging of anisotropic scatterers in a multi-layered medium \cite{SCC}, magnetoencephalography (MEG) from human cortical neural activities \cite{SNPM}, microwave imaging \cite{P-MUSIC6}, and breast cancer detection \cite{RSAAP}. We also refer to other studies \cite{A2,AIL1,CA,F3,K1,MGS,ML,OBP,P-MUSIC3,P-MUSIC5,PKLS,RSCGBA,S2} for more applications of MUSIC.

Let us emphasize that most studies have addressed full- and limited-view inverse problems. There also exists a considerable number of interesting limited-aperture inverse scattering problems; refer to \cite{AHP2,AJMP,AH,CAB,INS,KLAHP,MB,O,P-SUB13,Z} and the references therein. However, to the best of our knowledge, MUSIC has not been applied to the limited-aperture problem. Notice that MUSIC is based on the characterization of the range of the so-called multi-static response (MSR) matrix, which is symmetric but not Hermitian in full- and limited-view inverse problems. The difficulties that arise in the application of MUSIC in the limited-aperture problem come from the non-symmetric property of the MSR matrix. Accordingly, the imaging function of MUSIC is yet to be designed because an appropriate method to generate the projection operator onto the noise subspace does not exist.

In this study, we apply the MUSIC algorithm to the limited-aperture inverse scattering problem to identify or image small electromagnetic inhomogeneities, the dielectric permittivities (or magnetic permeabilities) of which differ compared with a homogeneous background. The first goal of this study is to design a MUSIC imaging function. This is based on the fact that the far-field pattern, which is the element of the MSR matrix, can be represented by an asymptotic expansion formula in the existence of small inhomogeneities, and the structures of left and right singular vectors of the MSR matrix are associated with nonzero singular values. The next goal is to analyze the mathematical structure of the imaging function to certify its feasibility and explore any fundamental limitations. For this, we prove that the imaging function is represented by an infinite series of first-kind integer-order Bessel functions and the configuration of incident and observation directions. On the basis of the analyzed structure, we identify that the main factors of the imaging function for the permittivity and permeability contrast cases are Bessel functions of orders zero and one, respectively. This further implies that the imaging performance is significantly dependent on the range of incident and observation directions; peaks of large magnitudes appear at the location of inhomogeneities for permittivity contrast case, and for the permeability contrast case, peaks of large magnitudes appear at the location of inhomogeneities when the range of such directions is narrow, while two peaks of large magnitudes appear in the neighborhood of the location of inhomogeneities when the range is wide enough. In light of this, a least condition for the proper imaging performance of MUSIC in the limited-aperture inverse scattering problem can be explored. The final goal is to exhibit numerical simulation results with synthetic data corrupted by random noise to demonstrate the feasibility and limitations of the designed imaging function and to support theoretical results.

This paper is organized as follows. In Section \ref{sec:2}, we briefly survey the two-dimensional direct scattering problem in the presence of well-separated small electromagnetic inhomogeneities and introduce the asymptotic expansion formula of the far-field pattern. In Section \ref{sec:3}, the imaging function of MUSIC in a limited-aperture inverse scattering problem is designed, the mathematical structure of an imaging function is established, and some properties (including pros and cons) of the imaging function are discussed. In Section \ref{sec:4}, corresponding numerical simulation results with noise-corrupted synthetic data generated by the Foldy-Lax framework \cite{HSZ4} are exhibited. In Section \ref{sec:5}, a short conclusion is provided, including an outline of future research.

\section{Two-dimensional direct scattering problem and traditional MUSIC algorithm}\label{sec:2}
First, we briefly survey two-dimensional direct scattering from a set of small electromagnetic inhomogeneities located in homogeneous space $\mathbb{R}^2$. {Throughout this paper, we denote such inhomogeneities as $\Sigma_s$, $s=1,2,\cdots,S$}, and assume that each $\Sigma_s$ can be expressed as 
\[\Sigma_s=\mr_s+\alpha_s\mB_s,\]
where $\mr_s$ denotes the location of $\Sigma_s$, and $\mB_s$ is a simple connected smooth domain containing the origin. We let $\Sigma$ denote the collection of $\Sigma_s$. {For the sake of simplicity, we assume all $\Sigma_s$ are balls (i.e., all $\mB_s$ are unit circles centered at the origin say, $\mathbb{S}^1$) with same radius $\alpha$ (i.e., $\alpha_s=\alpha$ for all $s=1,2,\cdots,S$)} and are characterized by the dielectric permittivity and magnetic permeability at the given angular frequency $\omega=2\pi f$, where $f$ denotes the ordinary frequency measured in \texttt{hertz}.

Let $\epsb$ and $\mub$ denote the value of dielectric permittivity and magnetic permeability of $\mathbb{R}^2$, respectively. {Analogously, we denote $\eps_s$ and $\mu_s$ be those of $\Sigma_s$}. Then, respectively, we can introduce the piecewise constants of dielectric permittivity and magnetic permeability as
\[\eps(\mr)=\left\{\begin{array}{ccl}
\eps_s&\mbox{~for~}&\mr\in\Sigma_s\\
\epsb&\mbox{for}&\mr\in\mathbb{R}^2\backslash\overline{\Sigma}
\end{array}\right.\quad\mbox{and}\quad\mu(\mr)=\left\{\begin{array}{ccl}
\mu_s&\mbox{~for~}&\mr\in\Sigma_s\\
\mub&\mbox{for}&\mr\in\mathbb{R}^2\backslash\overline{\Sigma}.
\end{array}\right.\]
With this, we denote $k$ be the background wavenumber that satisfies $k^2=\omega^2\epsb\mub$. Throughout this paper, we assume that all $\Sigma_s$ are well-separated from each other and, correspondingly, that $k$ satisfies 
\begin{equation}\label{Separation}
|\mr_s-\mr_{s'}|\gg\frac{3}{4k}\quad\text{for}\quad s\ne s'
\end{equation}

For a fixed wavenumber $k$, let $u(\mr,\vt)$ be the time-harmonic total field satisfying
\[\nabla\cdot\bigg(\frac{1}{\mu(\mr)}\nabla u(\mr,\vt)\bigg)+\omega^2\eps(\mr)u(\mr,\vt)=0,\]
with transmission conditions at the boundaries of $\Sigma_s$. Note that $u(\mr,\vt)$ can be split into the incident field $u_{\inc}(\mr,\vt)$ and scattered field $u_{\scat}(\mr,\vt)$ as \[u(\mr,\vt)=u_{\inc}(\mr,\vt)+u_{\scat}(\mr,\vt).\]
In this study, we consider the plane-wave illumination. Let $u_{\inc}(\mr,\vt)=e^{ik\vt\cdot\mr}$ be the given incident field with direction $\vt\in\mathbb{S}^1$ and $u_{\scat}(\mr,\vt)$ be the corresponding scattered field satisfying the Sommerfeld radiation condition, which can be expressed as
\[\lim_{|\mr|\to\infty}\sqrt{|\mr|}\left(\frac{\p u_{\scat}(\mr,\vt)}{\p|\mr|}-i ku_{\scat}(\mr,\vt)\right)=0\quad\text{uniformly in all directions}\quad\vv=\frac{\mr}{|\mr|}.\]

\subsection{Representation formula for the far-field pattern and the MUSIC algorithm: dielectric permittivity contrast case}
{First, let us consider the dielectric permittivity contrast case, i.e., $\eps(\mr)\ne\epsb$ and $\mu(\mr)=\mub$. Notice that based on \cite{CK,RC}, the induced source can be an electric line current, which generate 2D monopole radiation, also known as single-layer potential, the scattered field can be expressed by the single-layer potential with unknown density function $\varphi$ as
\begin{equation}\label{SingleLayer}
u_{\scat}(\mr,\vt)=\int_\Sigma\Gamma(\mr,\mr')\varphi(\mr',\vt)d\mr'=\sum_{s=1}^{S}\int_{\Sigma_s}\Gamma(\mr,\mr')\varphi(\mr',\vt)d\mr',
\end{equation}
where $\Gamma(\mr,\mr')$ denotes the fundamental solution to the Helmholtz equation, which can be expressed as
\[\Gamma(\mr,\mr')=-\frac{i}{4}H_0^{(1)}(k|\mr-\mr'|)=-\frac{i}{4}J_0(k|\mr-\mr'|)+\frac{1}{4}Y_0(k|\mr-\mr'|).\]
where $J_n$ and $Y_n$ denote the Bessel and Neumann functions of integer order $n$, respectively.}

{Let $u_\infty(\vv,\vt)$ be the far-field pattern of the scattered field $u_{\scat}(\mr,\vt)$ with observation direction $\vv\in\mathbb{S}^1$ that satisfies
\begin{equation}\label{FarFieldPattern}
u_{\scat}(\mr,\vt)=\frac{e^{ik|\mr|}}{\sqrt{|\mr|}}u_\infty(\vv,\vt)+O\left(\frac{1}{\sqrt{|\mr|}}\right)\quad\text{uniformly in all directions}\quad\vv=\frac{\mr}{|\mr|},\quad |\mr|\longrightarrow\infty.
\end{equation}
Based on \eqref{SingleLayer} and the asymptotic behavior of the Hankel function, $u_\infty(\vv,\vt)$ can be expressed as
\[u_\infty(\vv,\vt)=-\frac{1+i}{4\sqrt{k\pi}}\int_\Sigma e^{ik\vt\cdot\mr'}\varphi(\mr',\vt)d\mr'=-\frac{1+i}{4\sqrt{k\pi}}\sum_{s=1}^{S}\int_{\Sigma_s}e^{ik\vv\cdot\mr'}\varphi(\mr',\vt)d\mr'.\]}

Notice that, since $\varphi(\mr',\vt)$ is unknown, there are limitations with respect to designing the MUSIC algorithm, which means we need an alternative expression of $\varphi(\mr',\vt)$. Based on \cite{AK2}, the far-field pattern $u_\infty(\vv,\vt)$ can be represented as an asymptotic expansion formula, which plays a key role in designing the MUSIC algorithm.

\begin{lem}[Asymptotic formula: dielectric permittivity contrast case]
 For sufficiently large $\omega$, $u_\infty(\vv,\vt)$ can be represented as 
 \begin{equation}\label{AsymptoticEps}
 u_\infty(\vv,\vt)=\alpha^2|\mB_s|\frac{k^2(1+i)}{4\sqrt{k\pi}}\sum_{s=1}^{S} \left(\frac{\eps_s-\epsb}{\sqrt{\epsb\mub}}\right)e^{-ik(\vv-\vt)\cdot\mr_s}+o(\alpha^2).
 \end{equation}
where $|\mB_s|=\pi$ denotes the area of a unit circle.
\end{lem}

{Now, we introduce the traditional MUSIC algorithm to identify the locations of $\mr_s$ from a set of measured far-field patterns such that
\begin{equation}\label{SetFarField}
\Psi=\set{u_\infty(\vv_n,\vt_n):\vt_n=-\vv_n\in\mathbb{S}^1,~n=1,2,\cdots,N},
\end{equation}
where $\vt_n$ are given by
\[\vt_n=-\left[\cos\frac{2\pi n}{N},\sin\frac{2\pi n}{N}\right]^T.\]
From the collection of far-field pattern data $\Psi$, an MSR matrix can be generated as $\mathbb{K}\in\mathbb{C}^{N\times N}$:
\begin{equation}\label{MSR}
\mathbb{K}=\left[\begin{array}{cccc}
u_{\infty}(\vv_1,\vt_1) & u_{\infty}(\vv_1,\vt_2) & \cdots & u_{\infty}(\vv_1,\vt_N)\\
u_{\infty}(\vv_2,\vt_1) & u_{\infty}(\vv_2,\vt_2) & \cdots & u_{\infty}(\vv_2,\vt_N)\\
\vdots&\vdots&\ddots&\vdots\\
u_{\infty}(\vv_N,\vt_1) & u_{\infty}(\vv_N,\vt_2) & \cdots & u_{\infty}(\vv_N,\vt_N)\\
\end{array}\right].
\end{equation}
From \eqref{AsymptoticEps} and \eqref{SetFarField}, the elements of $\mathbb{K}$ can be approximated as
\begin{equation}\label{Approximation1}
u_\infty(\vv_m,\vt_n)=\alpha^2\pi\frac{k^2(1+i)}{4\sqrt{k\pi}}\sum_{s=1}^{S} \left(\frac{\eps_s-\epsb}{\sqrt{\epsb\mub}}\right)e^{ik(\vt_m+\vt_n)\cdot\mr_s}.
\end{equation}
Now, assume that $N>S$. Then the singular value decomposition (SVD) of $\mathbb{K}$ is given by
\[\mathbb{K}=\mathbb{USV}^*=\sum_{n=1}^{N}\sigma_n\mU_n\mV_n^*\approx\sum_{n=1}^{S}\sigma_n\mU_n\mV_n^*,\]
where the superscript $*$ represents the Hermitian operator, $\mU_s$ and $\mV_s$ are the left and right singular vectors of $\mathbb{K}$, respectively, and $\sigma_s$ denotes the singular value of $\mathbb{K}$ such that
\[\sigma_1\geq\sigma_2\geq\cdots\geq\sigma_{S}>0\quad\text{and}\quad\sigma_s\approx\rho\approx0\quad\text{for}\quad s>S.\]
Then, $\set{\mU_1,\mU_2,\cdots,\mU_{S}}$ spans the basis for the signal space of $\mathbb{K}$. Therefore, we can define the projection operator onto the null (or noise) subspace $\mathbf{P}_{\noise}:\mathbb{C}^{N\times1}\longrightarrow\mathbb{C}^{N\times1}$. This projection is given explicitly by
\[\mathbf{P}_{\noise}=\mathbb{I}_N-\sum_{s=1}^{S}\mU_s\mU_s^*.\]}

{Let us denote $\Omega$ as the region of interest and $\mr\in\Omega$. Then, based on \eqref{Approximation1}, we define a vector $\mW_\eps(\mr)\in\mathbb{C}^{N\times1}$ as
\[\mW_\eps(\mr):=\bigg[e^{ik\vt_1\cdot\mr},e^{ik\vt_2\cdot\mr},\cdots,e^{ik\vt_N\cdot\mr}\bigg].\]
Accordingly, there exists some $N_0^{(\eps)}\in\mathbb{N}$ such that, for any $N\geq N_0^{(\eps)}$ and $s=1,2,\cdots,S$
\[\mW_\eps(\mr)\in\mbox{Range}(\mathbb{K\overline{K}})=\mbox{Range}(\mathbb{KK^*})=\mbox{Range}(\mathbb{K^*K})\quad\text{if and only if}\quad\mr=\mr_s,\]
With this, the imaging function of the MUSIC algorithm is defined as
\begin{equation}\label{ImagingfunctionEps}
 f_{\music}(\mr)=\frac{1}{|\mathbf{P}_{\noise}(\mW_\eps(\mr))|},\quad\text{where}\quad|\mathbf{P}_{\noise}(\mW_\eps(\mr))|^2=\mathbf{P}_{\noise}(\mW_\eps(\mr))\cdot\overline{\mathbf{P}_{\noise}(\mW_\eps(\mr))}.
\end{equation}
Then, the map of $f_{\music}(\mr)$ will have peaks of large and small magnitudes at $\mr=\mr_s\in\Sigma_s$ and $\mr\in\mathbb{R}^2\backslash\Sigma$, respectively. A more detailed description can be found in \cite{C,AK2}.}

\subsection{Representation formula for the far-field pattern and the MUSIC algorithm: magnetic permeability contrast case}
{Next, we consider the magnetic permeability contrast case, i.e., $\eps(\mr)=\epsb$ and $\mu(\mr)\ne\mub$. In this case, based on the \cite{CK,RC}, the induced source can be magnetic dipoles along the transverse direction, which has a freedom of two, also known as double-layer potential. Hence, 
the scattered field can be expressed by the double-layer potential with unknown density function $\psi$ as
\[u_{\scat}(\mr,\vt)=\int_\Sigma\frac{\p\Gamma(\mr,\mr')}{\p\vn(\mr')}\psi(\mr',\vt)d\mr'=\sum_{s=1}^{S}\int_{\Sigma_s}\frac{\p\Gamma(\mr,\mr')}{\p\vn(\mr')}\psi(\mr',\vt)d\mr'\]
and the far-field pattern of \eqref{FarFieldPattern} can be defined analogously.}

{Same as the dielectric permittivity contrast case, the complete form of $\psi(\mr',\vt)$ is unknown so that its alternative expression is required. Based on \cite{AK2}, the far-field pattern $u_\infty(\vv,\vt)$ can be represented as an asymptotic expansion formula.}
\begin{lem}[Asymptotic formula: magnetic permeability contrast case]
 For sufficiently large $\omega$, $u_\infty(\vv,\vt)$ can be represented as 
 \begin{equation}\label{AsymptoticMu}
u_\infty(\vv,\vt)=\alpha^2|\mB_s|\frac{k^2(1+i)}{4\sqrt{k\pi}}\sum_{s=1}^{S} \bigg(\frac{\eps_s-\epsb}{\sqrt{\epsb\mub}}-\vv\cdot\mathbb{M}(\mr_s)\cdot\vt\bigg)e^{-ik(\vv-\vt)\cdot\mr_s}+o(\alpha^2).
\end{equation}
where $\mathbb{M}(\mr_s)$ is a $2\times2$ diagonal matrix with components $2\mub/(\mu_s+\mub)$.
\end{lem}

{Now, we introduce the traditional MUSIC algorithm to identify the locations of $\mr_s$ from a set of measured far-field patterns $\Psi$ of \eqref{SetFarField}. From the expression \eqref{AsymptoticMu}, the elements of $\mathbb{K}$ of \eqref{MSR} can be approximated as
\begin{equation}\label{Approximation2}
u_\infty(\vv_m,\vt_n)=\alpha^2\pi\frac{k^2(1+i)}{4\sqrt{k\pi}}\sum_{s=1}^{S} \left(\frac{2\mub}{\mu_s+\mub}(\vv_m\cdot\vt_n)\right)e^{ik(\vt_m+\vt_n)\cdot\mr_s}.
\end{equation}
On the basis of the above approximation \eqref{Approximation2} and with the assumption $N>2S$, the SVD of $\mathbb{K}$ can be written as
\[\mathbb{K}=\mathbb{USV}^*=\sum_{n=1}^{N}\sigma_n\mU_n\mV_n^*\approx\sum_{s=1}^{2S}\sigma_s\mU_s\mV_s^*\]
and correspondingly, in contrast to the permittivity contrast case, $\set{\mU_1,\mU_2,\cdots,\mU_{2S}}$ spans the basis for the signal space of $\mathbb{K}$. Therefore, the projection operator onto the null (or noise) subspace $\mathbf{P}_{\noise}:\mathbb{C}^{N\times1}\longrightarrow\mathbb{C}^{N\times1}$ is given explicitly by
\[\mathbf{P}_{\noise}=\mathbb{I}_N-\sum_{n=1}^{2S}\mU_n\mU_n^*.\]}

{Now, based on \eqref{Approximation2}, we define a vector $\mW_\mu(\mr)\in\mathbb{C}^{N\times1}$ such that for $\mathbf{c}_n\in\mathbb{R}^2\backslash\set{\mathbf{0}}$, $n=1,2,\cdots,N$,
\[\mW_\mu(\mr):=\bigg[(\mathbf{c}_1\cdot\vt_1)e^{ik\vt_1\cdot\mr},(\mathbf{c}_2\cdot\vt_1)e^{ik\vt_2\cdot\mr},\cdots,(\mathbf{c}_N\cdot\vt_N)e^{ik\vt_N\cdot\mr}\bigg].\]
Then, there exists some $N_0\in\mathbb{N}$ such that, for any $N\geq N_0$ and $s=1,2,\cdots,S$
\[\mW_\mu(\mr)\in\mbox{Range}(\mathbb{K\overline{K}})=\mbox{Range}(\mathbb{KK^*})=\mbox{Range}(\mathbb{K^*K})\quad\text{if and only if}\quad\mr=\mr_s,\]
With this, the imaging function of the MUSIC algorithm is defined as
\begin{equation}\label{ImagingfunctionMu}
 f_{\music}(\mr)=\frac{1}{|\mathbf{P}_{\noise}(\mW_\mu(\mr))|}.
\end{equation}
Then, the map of $f_{\music}(\mr)$ will have peaks of large and small magnitudes at $\mr=\mr_s\in\Sigma_s$ and $\mr\in\mathbb{R}^2\backslash\Sigma$, respectively.}

\begin{rem}
To apply MUSIC to determine unknown inhomogeneities, the condition of \eqref{SetFarField} must be satisfied, i.e., incident and observation directions must be opposite, and their total numbers must be the same. This means that the MSR matrix $\mathbb{K}$ must be complex symmetric (but not Hermitian in general). This condition holds in full- and limited-view inverse scattering problems but does not hold in the limited-aperture inverse scattering problem.
\end{rem}

\section{MUSIC algorithm in the limited-aperture problem: introduction, analysis, and various properties}\label{sec:3}
\subsection{MSR matrix in the limited-aperture configuration}
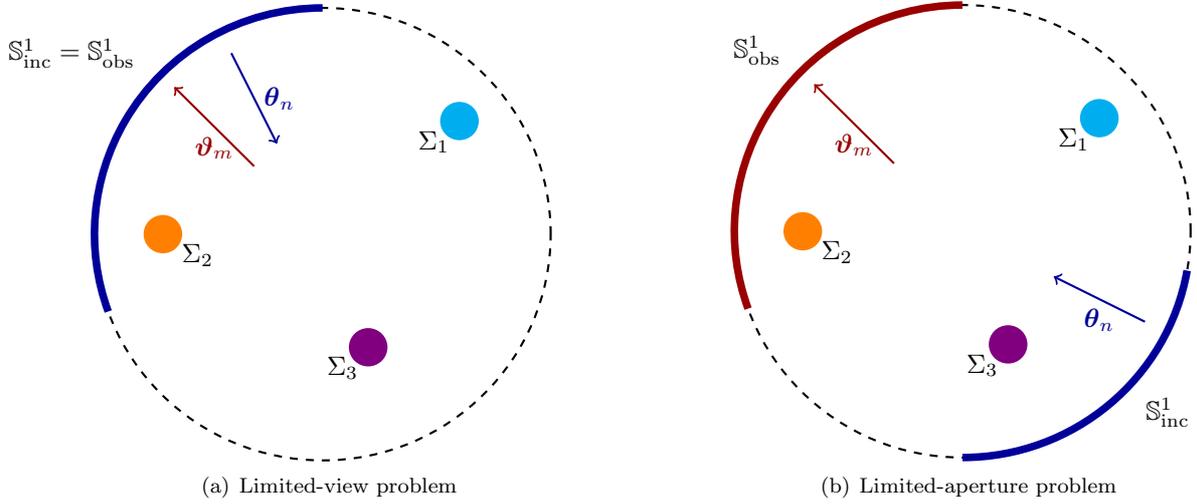
\begin{figure}[h]
\begin{center}
\subfigure[Limited-view problem]{\begin{tikzpicture}[scale=3]
\draw[black,thick,dashed] (0,0) circle (1cm);
\draw[blue!60!black,thick,solid,line width=1mm] (0,1) arc (90:200:1cm);
\node at (-1.1,0.8) {$\mathbb{S}_{\inc}^1=\mathbb{S}_{\obs}^1$};
\draw[cyan,thick,fill=cyan] (0.6,0.5) circle (0.08cm) node[below left,black] {$ \Sigma_1 $};
\draw[orange,thick,fill=orange] (-0.7,0.0) circle (0.08cm) node[below right,black] {$~\Sigma_2 $};
\draw[violet,thick,fill=violet] (0.2,-0.5) circle (0.08cm) node[below left,black] {$ \Sigma_3 $};
\draw[blue!60!black,thick,solid,->] (-0.4,0.8) -- node[right] {$\vt_n$} (-0.2,0.4);
\draw[red!60!black,thick,solid,->] (-0.3,0.3) -- node[below] {$ \vv_m $} (-0.65,0.65);
\end{tikzpicture}\qquad\qquad}\qquad
\subfigure[Limited-aperture problem]{
\begin{tikzpicture}[scale=3]
\draw[black,thick,dashed] (0,0) circle (1cm);
\draw[blue!60!black,thick,solid,line width=1mm] (0,-1) arc (-90:-10:1cm);
\draw[red!60!black,thick,solid,line width=1mm] (0,1) arc (90:200:1cm);
\node at (0.9,-0.8) {$\mathbb{S}_{\inc}^1$};
\node at (-0.9,0.8) {$\mathbb{S}_{\obs}^1$};
\draw[cyan,thick,fill=cyan] (0.6,0.5) circle (0.08cm) node[below left,black] {$ \Sigma_1 $};
\draw[orange,thick,fill=orange] (-0.7,0.0) circle (0.08cm) node[below right,black] {$~\Sigma_2 $};
\draw[violet,thick,fill=violet] (0.2,-0.5) circle (0.08cm) node[below left,black] {$ \Sigma_3 $};

\draw[blue!60!black,thick,solid,->] (0.8,-0.4) -- node[below] {$\vt_n$} (0.4,-0.2);
\draw[red!60!black,thick,solid,->] (-0.3,0.3) -- node[below] {$ \vv_m $} (-0.65,0.65);
\end{tikzpicture}\quad}
\caption{\label{Configuration}Illustration of limited-view and limited-aperture inverse scattering problems.}
\end{center}
\end{figure}

In this section, we consider the imaging function of the MUSIC algorithm to identify the locations of $\mr_s$ from a set of measured far-field patterns such that
\[\Psi=\set{u_\infty(\vv_m,\vt_n):\vv_m\in\mathbb{S}_{\obs}^1,~\vt_n\in\mathbb{S}_{\inc}^1,~m=1,2,\cdots,M,~n=1,2,\cdots,N},\]
where $\mathbb{S}_{\obs}^1$ and $\mathbb{S}_{\inc}^1$ denote the set of observation and incident directions, respectively, and {$\vv_m$ and $\vt_n$ are given by
\begin{align*}
\vv_m&=[\cos\vartheta_m,\sin\vartheta_m]^T,\quad\vartheta_m=\vartheta_1+\frac{(m-1)}{M-1}(\vartheta_M-\vartheta_1),\\
\intertext{and}
\vt_n&=[\cos\theta_n,\sin\theta_n]^T,\quad\theta_n=\theta_1+\frac{(n-1)}{N-1}(\theta_N-\theta_1),
\end{align*}}
respectively. {Notice that in limited-view problem, one assumes that the sets of incident and observation directions are coincide, i.e., $\mathbb{S}_{\inc}=\mathbb{S}_{\obs}$ because $\vv_n=-\vt_n$ for all $n=1,2,\cdots,N$. However, in the current limited-aperture problem, we assume that the sets and total number of incident and observation directions are different, i.e., $\mathbb{S}_{\inc}\ne\mathbb{S}_{\obs}$ and $N\ne M$. Figure \ref{Configuration} illustrates the comparison between the limited-aperture and limited-view inverse scattering problems.} Throughout this paper, we assume that $\mathbb{S}_{\obs}^1$ and $\mathbb{S}_{\inc}^1$ are connected and proper subsets of $\mathbb{S}^1$.

\begin{figure}
\centering
\end{figure}

As previously, let us generate an MSR matrix $\mathbb{K}\in\mathbb{C}^{M\times N}$ from $\Psi$ as
\[\mathbb{K}=\left[\begin{array}{cccc}
u_\infty(\vv_1,\vt_1) & u_\infty(\vv_1,\vt_2) & \cdots & u_\infty(\vv_1,\vt_N)\\
u_\infty(\vv_2,\vt_1) & u_\infty(\vv_2,\vt_2) & \cdots & u_\infty(\vv_2,\vt_N)\\
\vdots&\vdots&\ddots&\vdots\\
u_\infty(\vv_M,\vt_1) & u_\infty(\vv_M,\vt_2) & \cdots & u_\infty(\vv_M,\vt_N)\\
\end{array}\right].\]
With this, we can design the MUSIC algorithm, analyze the mathematical structure of the imaging function, and explore various MUSIC properties for {two} different cases: dielectric permittivity contrast only ($\eps(\mr)\ne\epsb$ and $\mu(\mr)=\mub$) and magnetic permeability contrast only ($\eps(\mr)=\epsb$ and $\mu(\mr)\ne\mub$). {Both contrast case ($\eps(\mr)\ne\epsb$ and $\mu(\mr)\ne\mub$) can be derived by combining the permittivity and permeability contrast cases.}

\subsection{Dielectric permittivity contrast case: $\eps(\mr)\ne\epsb$ and $\mu(\mr)=\mub$}
First, we consider the dielectric permittivity contrast case. For the sake of simplicity, {we assume that $M,N>S$} and $\eps_s>\epsb$ for all $s$. Then, based on the expression \eqref{Approximation1}, $\mathbb{K}$ can be decomposed as 
\begin{equation}\label{DecompositionK1}
 \mathbb{K}=\mathbb{EDF}^T+\mathbb{O}(\alpha^2)\approx\mathbb{EDF}^T,
\end{equation}
where $\mathbb{O}(\alpha^2)$ is a $M\times N$ matrix, the elements of which are $o(\alpha^2)$, and $\mathbb{D}$ is the $S\times S$ matrix
\[\mathbb{D}=\alpha^2\pi\frac{k^2(1+i)}{4\sqrt{k\pi}}\left[\begin{array}{cccc}
\displaystyle\frac{\eps_1-\epsb}{\sqrt{\epsb\mub}}&0&0&0\\
0&\displaystyle\frac{\eps_2-\epsb}{\sqrt{\epsb\mub}}&0&0\\
\vdots&\vdots&\ddots&\vdots\\
0&0&0&\displaystyle\frac{\eps_S-\epsb}{\sqrt{\epsb\mub}}\\
\end{array}\right]\in\mathbb{C}^{S\times S},\]
and the matrices $\mathbb{E}\in\mathbb{C}^{M\times S}$ and $\mathbb{F}\in\mathbb{C}^{N\times S}$ respectively are
\[\mathbb{E}=\left[\begin{array}{cccc}
e^{-ik\vv_1\cdot\mr_1}&e^{-ik\vv_1\cdot\mr_2}&\cdots&e^{-ik\vv_1\cdot\mr_S}\\
e^{-ik\vv_2\cdot\mr_1}&e^{-ik\vv_2\cdot\mr_2}&\cdots&e^{-ik\vv_2\cdot\mr_S}\\
\vdots&\vdots&\ddots&\vdots\\
e^{-ik\vv_M\cdot\mr_1}&e^{-ik\vv_M\cdot\mr_2}&\cdots&e^{-ik\vv_M\cdot\mr_S}
\end{array}\right]\quad\text{and}\quad
\mathbb{F}=\left[\begin{array}{cccc}
e^{ik\vt_1\cdot\mr_1}&e^{ik\vt_1\cdot\mr_2}&\cdots&e^{ik\vt_1\cdot\mr_S}\\
e^{ik\vt_2\cdot\mr_1}&e^{ik\vt_2\cdot\mr_2}&\cdots&e^{ik\vt_2\cdot\mr_S}\\
\vdots&\vdots&\ddots&\vdots\\
e^{ik\vt_N\cdot\mr_1}&e^{ik\vt_N\cdot\mr_2}&\cdots&e^{ik\vt_N\cdot\mr_S}
\end{array}\right],\]
Let us perform SVD of $\mathbb{K}$ as
\begin{equation}\label{SVD1}
\mathbb{K}=\mathbb{USV}^*=\sum_{n=1}^{N}\sigma_n\mU_n\mV_n^*\approx\sum_{n=1}^{S}\sigma_n\mU_n\mV_n^*.\end{equation}
Then, the first $S-$columns of singular vectors $\set{\mU_1,\mU_2,\cdots,\mU_S}$ and $\set{\mV_1,\mV_2,\cdots,\mV_S}$ span the signal spaces of $\mathbb{KK^*}$ and $\mathbb{K^*K}$, respectively. Correspondingly, we can define orthonormal projection operators onto the noise subspaces as
\[\mathbb{P}_{\noise}^{(\eps)}=\mathbb{I}(S)-\sum_{s=1}^{S}\mU_s\mU_s^*\quad\text{and}\quad\mathbb{Q}_{\noise}^{(\eps)}=\mathbb{I}(S)-\sum_{s=1}^{S}\mV_s\mV_s^*,\]
where $\mathbb{I}(S)$ denotes the $S\times S$ identity matrix.

Based on \eqref{Approximation1}, we define test vectors for $\mr\in\Omega$ as
\[\mf_\eps(\mr)=\frac{1}{\sqrt{M}}\bigg[e^{-ik\vv_1\cdot\mr},e^{-ik\vv_2\cdot\mr},\cdots,e^{-ik\vv_M\cdot\mr}\bigg]^T\quad\text{and}\quad
\mg_\eps(\mr)=\frac{1}{\sqrt{N}}\bigg[e^{ik\vt_1\cdot\mr},e^{ik\vt_2\cdot\mr},\cdots,e^{ik\vt_N\cdot\mr}\bigg]^T.\]
Then, based on \cite{AK2}, it can be concluded that
\begin{align*}
\mf_\eps(\mr)\in\range(\mathbb{KK^*})&\quad\text{if and only if}\quad\mr\in\set{\mr_1,\mr_2,\cdots,\mr_S}\\
\mg_\eps(\mr)\in\range(\mathbb{K^*K})&\quad\text{if and only if}\quad\mr\in\set{\mr_1,\mr_2,\cdots,\mr_S}
\end{align*}
This means that $|\mathbb{P}_{\noise}^{(\eps)}(\mf_\eps(\mr))|=0$ and $|\mathbb{Q}_{\noise}^{(\eps)}(\mg_\eps(\mr))|=0$ when $\mr\in\Sigma$. Thus, the location of $\mr_s\in\Sigma_s$ can be identified by plotting
\begin{equation}\label{Imagingfunction1}
f_{\music}(\mr)=\frac{1}{2}\left(\frac{1}{|\mathbb{P}_{\noise}^{(\eps)}(\mf_\eps(\mr))|}+\frac{1}{|\mathbb{Q}_{\noise}^{(\eps)}(\mg_\eps(\mr))|}\right).
\end{equation}
The resulting plot of $f_{\music}(\mr)$ is expected to exhibit large-magnitude peaks (theoretically, $+\infty$) at $\mr\in\Sigma_s$ for $s=1,2,\cdots,S$.

\begin{rem}[Comparison with traditional MUSIC]
In the traditional MUSIC algorithm, $\mf_\eps(\mr)=\mg_\eps(\mr)$ and $\mathbb{P}_{\noise}^{(\eps)}=\mathbb{Q}_{\noise}^{(\eps)}$ because $\mathbb{K}$ is a complex symmetric matrix. In this case, \eqref{Imagingfunction1} becomes \eqref{ImagingfunctionEps}.
\end{rem}

To explain the feasibility of $f_{\music}(\mr)$, we establish a relationship with an infinite series of integer-order Bessel functions. The derivation is given in Section \ref{sec:A}.

\begin{thm}\label{Theorem1}
Let $\vv_m:=[\cos\vartheta_m,\sin\vartheta_m]^T$, $\vt_n:=[\cos\theta_n,\sin\theta_n]^T$, $\mr:=|\mr|[\cos\phi,\sin\phi]^T$, and $\mr-\mr_s:=|\mr-\mr_s|[\cos\phi_s,\sin\phi_s]^T$. Then, for sufficiently large $k$, $f_{\music}(\mr)$ can be represented as 
\begin{align}
\begin{aligned}\label{Structure1}
f_{\music}(\mr)=&\frac{1}{2}\left(1-\sum_{s=1}^{S}\left|J_0(k|\mr-\mr_s|)+\frac{\Lambda_{\eps}^{(1)}(\mr-\mr_s)}{\vartheta_M-\vartheta_1}\right|^2+O\left(\frac{1}{(\vartheta_M-\vartheta_1)^2}\right)\right)^{-1/2}\\
&+\frac{1}{2}\left(1-\sum_{s=1}^{S}\left|J_0(k|\mr-\mr_s|)+\frac{\Lambda_{\eps}^{(2)}(\mr-\mr_s)}{\theta_N-\theta_1}\right|^2+O\left(\frac{1}{(\theta_N-\theta_1)^2}\right)\right)^{-1/2},
\end{aligned}
\end{align}
 where $\Lambda_{\eps}^{(1)}(\mr-\mr_s)$ and $\Lambda_{\eps}^{(2)}(\mr-\mr_s)$ are given by \eqref{Lambda_Eps1} and \eqref{Lambda_Eps2}, respectively.
\end{thm}

Based on the identified mathematical structure of \eqref{Structure1}, the intrinsic properties of MUSIC for the dielectric permittivity contrast case can be identified.

\begin{discuss}[Feasibility and limitation of MUSIC]\label{Discussion-E1}
The imaging function $f_{\music}(\mr)$ consists of $J_0(k|\mr-\mr_s|)$, $\Lambda_{\eps}^{(1)}(\mr-\mr_s)/(\vartheta_M-\vartheta_1)$, and $\Lambda_{\eps}^{(2)}(\mr-\mr_s)/(\theta_N-\theta_1)$. The first term is independent of the range of incident and observation directions and contributes to the detection, whereas the remaining terms are dependent on the range and disturb the detection. Since $J_0(0)=1$ and $J_p(0)=0$ for every positive integer $p$, the value of $f_{\music}(\mr_s)$ will be sufficiently large to identify $\mr_s$ via the map of $f_{\music}(\mr)$. However, if $\mr\ne\mr_s$ and the range of incident or observations directions is narrow, $f_{\music}(\mr)$ will be dominated by $\Lambda_{\eps}^{(1)}(\mr-\mr_s)/(\vartheta_M-\vartheta_1)$ or $\Lambda_{\eps}^{(2)}(\mr-\mr_s)/(\theta_N-\theta_1)$. Therefore, identifying the location of $\Sigma_s$ should be difficult.
\end{discuss}

\begin{discuss}[Least range of directions]\label{Discussion-E2}
On the basis of \eqref{Structure1}, it will be possible to obtain a good result by eliminating
\[\frac{\Lambda_{\eps}^{(1)}(\mr-\mr_s)}{\vartheta_M-\vartheta_1}=0\quad\text{and}\quad\frac{\Lambda_{\eps}^{(2)}(\mr-\mr_s)}{\theta_N-\theta_1}=0.\]
A possible choice is to select $\vartheta_1$, $\vartheta_M$, $\theta_1$, and $\theta_N$, thereby satisfying
\[\sin\frac{p(\vartheta_M-\vartheta_1)}{2}\cos\frac{p(\vartheta_M+\vartheta_1-2\phi_s+2\pi)}{2}=0\quad\text{and}\quad\sin\frac{q(\theta_N-\theta_1)}{2}\cos\frac{q(\theta_N+\theta_1-2\phi_s)}{2}=0\]
{for every positive integer $p$ and $q$}. Another possible choice is the selection $\vartheta_1=\phi_s$, $\vartheta_M=\pi+\phi_s$, $\theta_1=\phi_s$, and $\theta_N=\pi+\phi_s$. Unfortunately, we do not have a priori information of $\Sigma_s$, this is an ideal condition, but motivated by this observation, one can obtain a good result when $\vartheta_M-\vartheta_1=\pi$ and $\theta_N-\theta_1=\pi$, i.e., if the range of incident and observation directions is wider than $\pi$, refer to Figures \ref{Result1-2} and \ref{Result2-2}.
\end{discuss}

\subsection{Magnetic permeability contrast case: $\eps(\mr)=\epsb$ and $\mu(\mr)\ne\mub$}
Next, we consider the magnetic permeability contrast case. Here, {we assume that $M,N>2S$} and $\mu_s<\mub$ for all $s$. Then, based on the expression \eqref{Approximation2}, $\mathbb{K}$ can be decomposed as 
\begin{equation}\label{DecompositionK2}
 \mathbb{K}=\mathbb{HBG}^T+\mathbb{O}(\alpha^2)\approx\mathbb{HBG}^T,
\end{equation}
where $\mathbb{O}(\alpha^2)$ is a $2M\times2N$ matrix, the elements of which are $o(\alpha^2)$, and $\mathbb{D}$ is the $2S\times2S$ matrix
\[\mathbb{D}=\alpha^2\pi\frac{k^2(1+i)}{4\sqrt{k\pi}}\left[\begin{array}{cccc}
\mathbb{D}_1&0&0&0\\
0&\mathbb{D}_2&0&0\\
\vdots&\vdots&\ddots&\vdots\\
0&0&0&\mathbb{D}_S\\
\end{array}\right]\in\mathbb{R}^{2S\times2S},\quad\mathbb{D}_s=\left[\begin{array}{cc}
\displaystyle\frac{\mub}{\mu_s+\mub}&0\\
0&\displaystyle\frac{\mub}{\mu_s+\mub}
\end{array}\right]\]
and the matrices $\mathbb{H}\in\mathbb{R}^{2M\times2S}$ and $\mathbb{M}\in\mathbb{R}^{2N\times2S}$ respectively are
\begin{align*}
\mathbb{H}&=\left[\begin{array}{ccccccc}
\mH_1^{(1)}&\mH_1^{(2)}&\mH_2^{(1)}&\mH_2^{(2)}&\cdots&\mH_{S}^{(1)}&\mH_{S}^{(2)}
\end{array}\right]\\
&=\left[\begin{array}{ccccc}
(-\vv_1\cdot\me_1)e^{-ik\vv_1\cdot\mr_1}&(-\vv_1\cdot\me_2)e^{-ik\vv_1\cdot\mr_1}&\cdots&(-\vv_1\cdot\me_1)e^{-ik\vv_1\cdot\mr_S}&(-\vv_1\cdot\me_2)e^{-ik\vv_1\cdot\mr_S}\\
(-\vv_2\cdot\me_1)e^{-ik\vv_2\cdot\mr_1}&(-\vv_2\cdot\me_2)e^{-ik\vv_2\cdot\mr_1}&\cdots&(-\vv_2\cdot\me_1)e^{-ik\vv_2\cdot\mr_S}&(-\vv_2\cdot\me_2)e^{-ik\vv_2\cdot\mr_S}\\
\vdots&\vdots&\ddots&\vdots&\vdots\\
(-\vv_M\cdot\me_1)e^{-ik\vv_M\cdot\mr_1}&(-\vv_M\cdot\me_2)e^{-ik\vv_M\cdot\mr_1}&\cdots&(-\vv_M\cdot\me_1)e^{-ik\vv_M\cdot\mr_S}&(-\vv_M\cdot\me_2)e^{-ik\vv_M\cdot\mr_S}\\
\end{array}\right]
\end{align*}
and
\begin{align*}
\mathbb{G}&=\left[\begin{array}{ccccccc}
\mG_1^{(1)}&\mG_1^{(2)}&\mG_2^{(1)}&\mG_2^{(2)}&\cdots&\mG_{S}^{(1)}&\mG_{S}^{(2)}
\end{array}\right]\\
&=\left[\begin{array}{ccccc}
(\vt_1\cdot\me_1)e^{-ik\vv_1\cdot\mr_1}&(\vt_1\cdot\me_2)e^{-ik\vv_1\cdot\mr_1}&\cdots&(\vt_1\cdot\me_1)e^{-ik\vv_1\cdot\mr_S}&(\vt_1\cdot\me_2)e^{-ik\vv_1\cdot\mr_S}\\
(\vt_2\cdot\me_1)e^{-ik\vv_2\cdot\mr_1}&(\vt_2\cdot\me_2)e^{-ik\vv_2\cdot\mr_1}&\cdots&(\vt_2\cdot\me_1)e^{-ik\vv_2\cdot\mr_S}&(\vt_2\cdot\me_2)e^{-ik\vv_2\cdot\mr_S}\\
\vdots&\vdots&\ddots&\vdots&\vdots\\
(\vt_N\cdot\me_1)e^{-ik\vv_N\cdot\mr_1}&(\vt_N\cdot\me_2)e^{-ik\vv_N\cdot\mr_1}&\cdots&(\vt_N\cdot\me_1)e^{-ik\vv_N\cdot\mr_S}&(\vt_N\cdot\me_2)e^{-ik\vv_N\cdot\mr_S}\\
\end{array}\right],
\end{align*}
where $\me_1=[1,0]^T$ and $\me_2=[0,1]^T$.

As in the permittivity contrast case, let us perform SVD of $\mathbb{K}$ as
\begin{equation}\label{SVD2}
\mathbb{K}=\mathbb{USV}^*=\sum_{n=1}^{N}\sigma_n\mU_n\mV_n^*\approx\sum_{n=1}^{2S}\sigma_n\mU_n\mV_n^*.
\end{equation}
Then, the first $2S-$columns of singular vectors $\set{\mU_1,\mU_2,\cdots,\mU_{2S}}$ and $\set{\mV_1,\mV_2,\cdots,\mV_{2S}}$ span the signal spaces of $\mathbb{KK^*}$ and $\mathbb{K^*K}$, respectively. Hence, we can define orthonormal projection operators onto the noise subspaces as
\[\mathbb{P}_{\noise}^{(\mu)}=\mathbb{I}(S)-\sum_{s=1}^{2S}\mU_s\mU_s^*\quad\text{and}\quad\mathbb{Q}_{\noise}^{(\mu)}=\mathbb{I}(S)-\sum_{s=1}^{2S}\mV_s\mV_s^*,\]

Now, on the basis of the expression \eqref{Approximation2}, let us introduce the following test vectors for $\mr\in\Omega$, $\vx_1,\vx_2\in\mathbb{R}^2\backslash\set{\mathbf{0}}$,
\begin{align*}
\mf_\mu(\mr)&=\frac{1}{\sqrt{C_1}}\bigg[(-\vv_1\cdot\vx_1)e^{-ik\vv_1\cdot\mr},(-\vv_2\cdot\vx_1)e^{-ik\vv_2\cdot\mr},\cdots,(-\vv_M\cdot\vx_1)e^{-ik\vv_M\cdot\mr}\bigg]^T
\intertext{and}
\mg_\mu(\mr)&=\frac{1}{\sqrt{C_2}}\bigg[(\vt_1\cdot\vx_2)e^{ik\vt_1\cdot\mr},(\vt_2\cdot\vx_2)e^{ik\vt_2\cdot\mr},\cdots,(\vt_N\cdot\vx_2)e^{ik\vt_N\cdot\mr}\bigg]^T,
\end{align*}
where
\[C_1=\frac{\vartheta_M-\vartheta_1}{2}+\frac{1}{2}\cos(\vartheta_M+\vartheta_1)\sin(\vartheta_M-\vartheta_1)\quad\text{and}\quad C_2=\frac{\theta_N-\theta_1}{2}+\frac{1}{2}\cos(\theta_N+\theta_1)\sin(\theta_N-\theta_1),\]
respectively. Based on \cite{AK2}, if $\vx_1$ and $\vx_2$ are the linear combination of $\me_1$ and $\me_2$, then we obtain
\begin{align*}
\mf_\mu(\mr)\in\range(\mathbb{KK^*})&\quad\text{if and only if}\quad\mr\in\set{\mr_1,\mr_2,\cdots,\mr_S}\\
\mg_\mu(\mr)\in\range(\mathbb{K^*K})&\quad\text{if and only if}\quad\mr\in\set{\mr_1,\mr_2,\cdots,\mr_S}
\end{align*}
This means that $|\mathbb{P}_{\noise}^{(\mu)}(\mf_\mu(\mr))|=0$ and $|\mathbb{Q}_{\noise}^{(\mu)}(\mg_\mu(\mr))|=0$ when $\mr\in\Sigma$. Thus, the location of $\mr_s\in\Sigma_s$ can be identified by plotting
\begin{equation}\label{Imagingfunction2}
f_{\music}(\mr)=\frac{1}{2}\left(\frac{1}{|\mathbb{P}_{\noise}^{(\mu)}(\mf_\mu(\mr))|}+\frac{1}{|\mathbb{Q}_{\noise}^{(\mu)}(\mg_\mu(\mr))|}\right).
\end{equation}
The resulting plot of $f_{\music}(\mr)$ is expected to exhibit large-magnitude peaks (theoretically, $+\infty$) at $\mr\in\Sigma_s$ for $s=1,2,\cdots,S$.

\begin{rem}[Selection of test vector]
Based on the structure of $\mathbb{H}$ and $\mathbb{M}$, $\vx_1$ and $\vx_2$ must be a linear combination of $\me_1$ and $\me_2$. Roughly speaking, if one has a priori information, i.e., the size of inhomogeneities is already known, the selection $\vx_1=\me_1$ and $\vx_2=\me_2$ will guarantee a good result. Unfortunately, since we do not have target information, estimating optimal $\vx_1$ and $\vx_2$ requires large computational costs; refer to \cite{PL1,HSZ1}. Hence, based on \cite{P-MUSIC1,HSZ1,P-MUSIC5,P-SUB3}, we apply the test vectors $\set{\mf_\eps(\mr),\mg_\eps(\mr)}$ instead of $\set{\mf_\mu(\mr),\mg_\mu(\mr)}$.
\end{rem}

Now, we explore the mathematical structure of $f_{\music}(\mr)$ by establishing a relationship with an infinite series of integer-order Bessel functions. The derivation is given in Section \ref{sec:B}.

\begin{thm}\label{Theorem2}
Let $\vv_m:=[\cos\vartheta_m,\sin\vartheta_m]^T$, $\vt_n:=[\cos\theta_n,\sin\theta_n]^T$, $\mr=|\mr|[\cos\phi,\sin\phi]^T$, and $\mr-\mr_s:=|\mr-\mr_s|[\cos\phi_s,\sin\phi_s]^T$. Then, for sufficiently large $k$ and test vectors $\set{\mf_\eps(\mr),\mg_\eps(\mr)}$, $f_{\music}(\mr)$ can be represented as 
\begin{align}
\begin{aligned}\label{Structure2}
f_{\music}(\mr)=&\frac{1}{2}\left(1-\sum_{s=1}^{S}\sum_{h=1}^{2}\left|iJ_1(k|\mr-\mr_s|)\bigg(\frac{\mr-\mr_s}{|\mr-\mr_s|}\cdot\me_h\bigg)+\frac{\Lambda_{\mu}^{(1,h)}(\mr-\mr_s)}{\vartheta_M-\vartheta_1}\right|^2+O\left(\frac{1}{(\vartheta_M-\vartheta_1)^2}\right)\right)^{-1/2}\\
&+\frac{1}{2}\left(1-\sum_{s=1}^{S}\sum_{h=1}^{2}\left|iJ_1(k|\mr-\mr_s|)\bigg(\frac{\mr-\mr_s}{|\mr-\mr_s|}\cdot\me_h\bigg)+\frac{\Lambda_{\mu}^{(2,h)}(\mr-\mr_s)}{\theta_N-\theta_1}\right|^2+O\left(\frac{1}{(\theta_N-\theta_1)^2}\right)\right)^{-1/2},
\end{aligned}
\end{align}
where $\Lambda_{\mu}^{(1,h)}(\mr-\mr_s)$ and $\Lambda_{\mu}^{(2,h)}(\mr-\mr_s)$ for $h=1,2$ are given by \eqref{Lambda_Mu1-1}, \eqref{Lambda_Mu1-2}, \eqref{Lambda_Mu2-1}, and \eqref{Lambda_Mu2-2}.
\end{thm}

Based on the identified mathematical structure of \eqref{Structure2}, we can discuss some properties of MUSIC for the magnetic permeability contrast case.

\begin{figure}
\centering
\subfigure[Case 1]{\includegraphics[width=.250\textwidth]{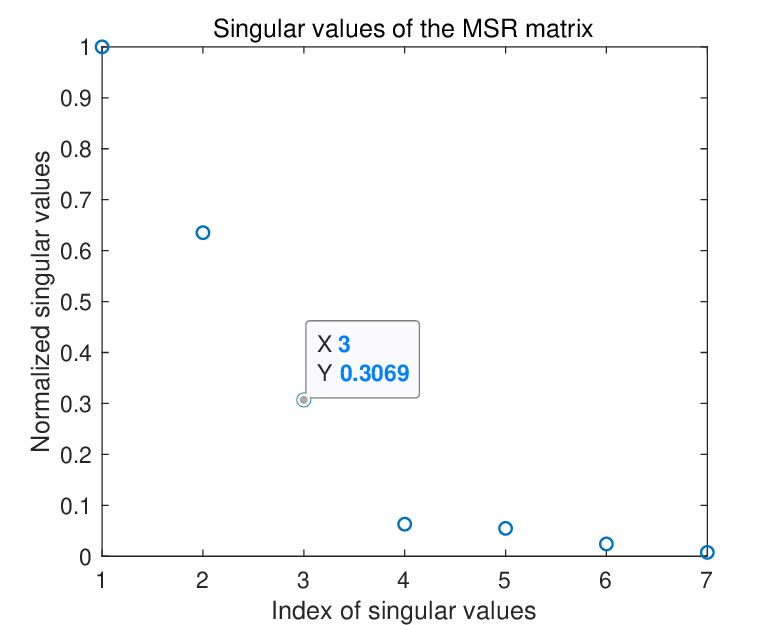}\hfill\includegraphics[width=.250\textwidth]{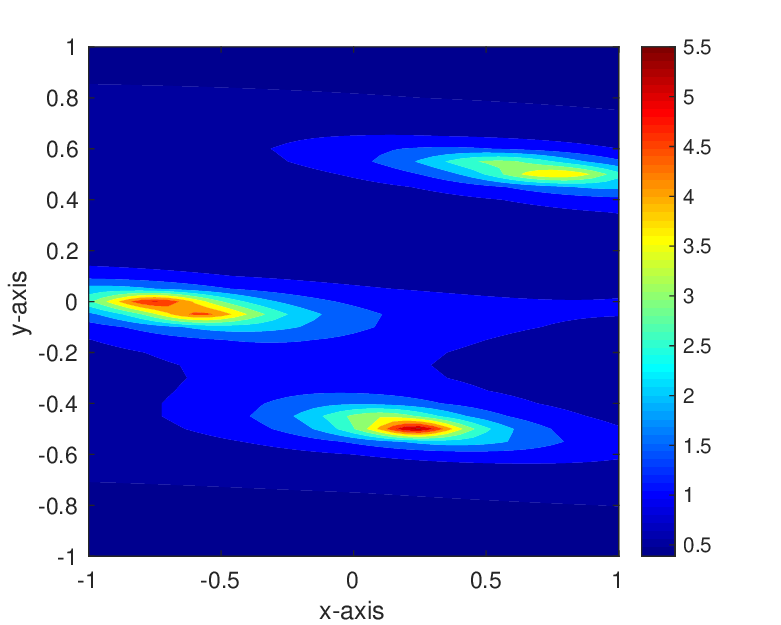}}\hfill
\subfigure[Case 2]{\includegraphics[width=.250\textwidth]{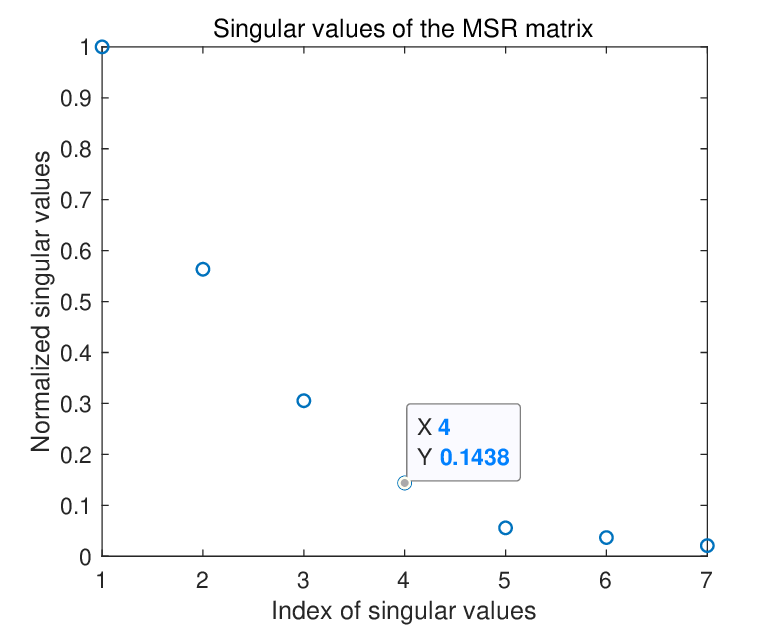}\hfill\includegraphics[width=.250\textwidth]{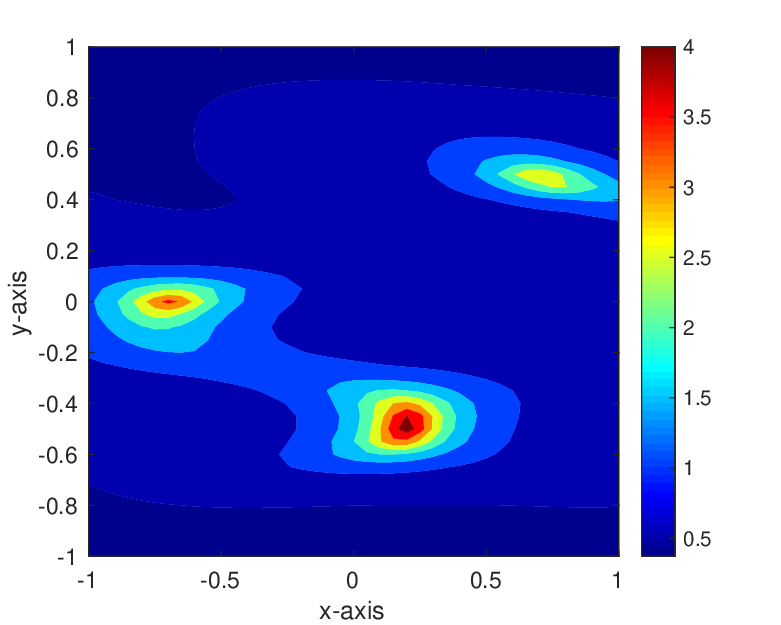}}\\
\subfigure[Case 3]{\includegraphics[width=.250\textwidth]{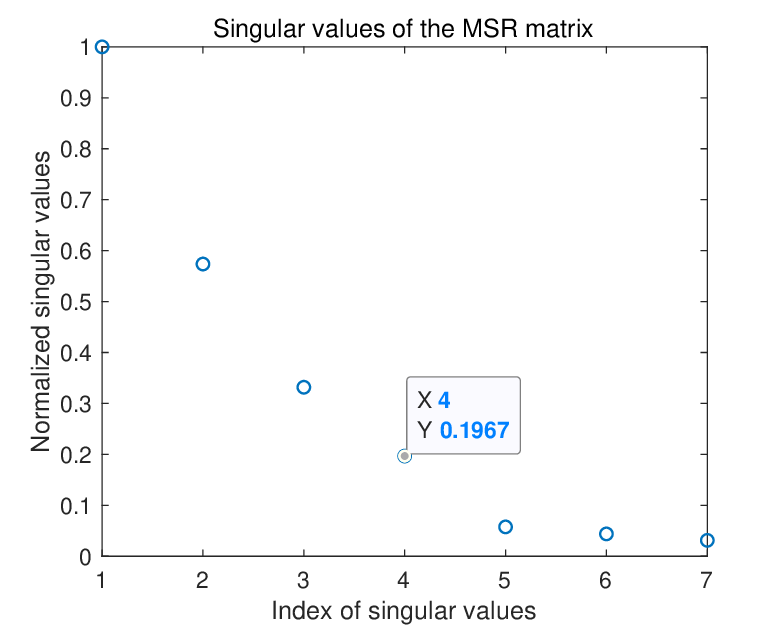}\hfill\includegraphics[width=.250\textwidth]{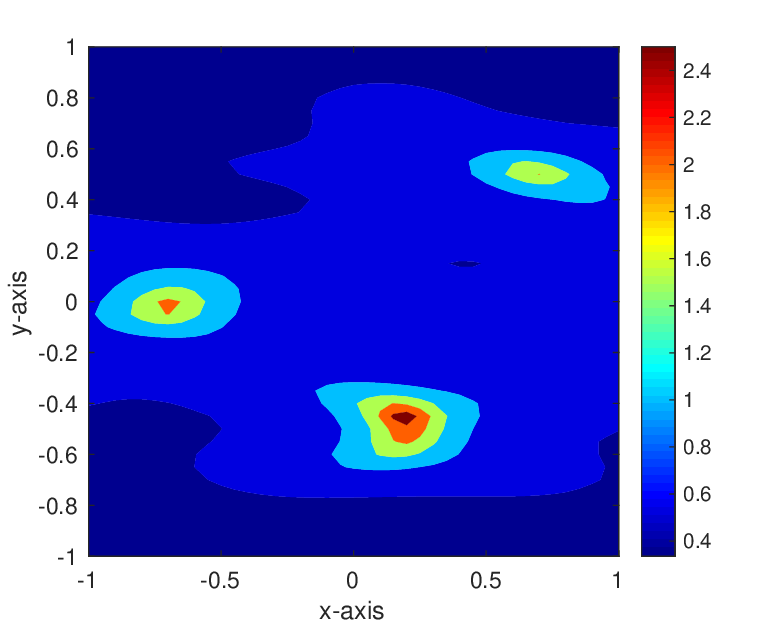}}\hfill
\subfigure[Case 4]{\includegraphics[width=.250\textwidth]{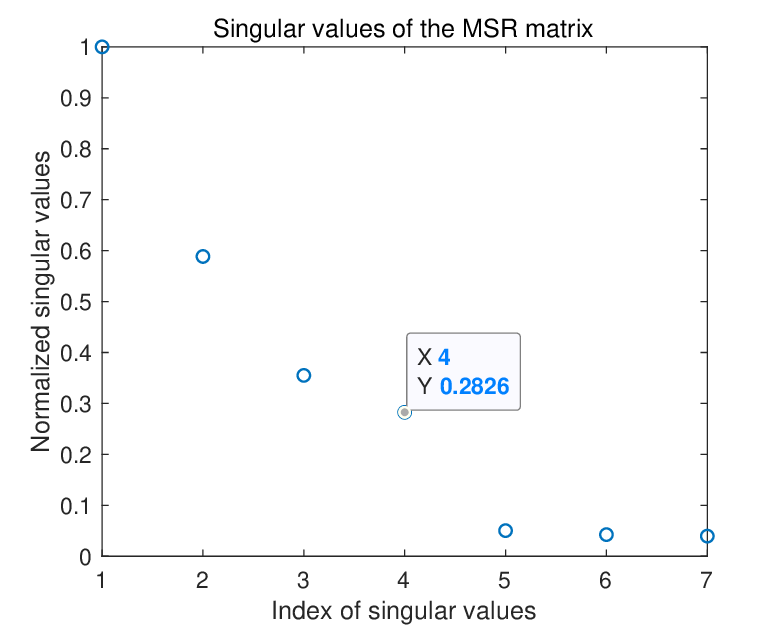}\hfill\includegraphics[width=.250\textwidth]{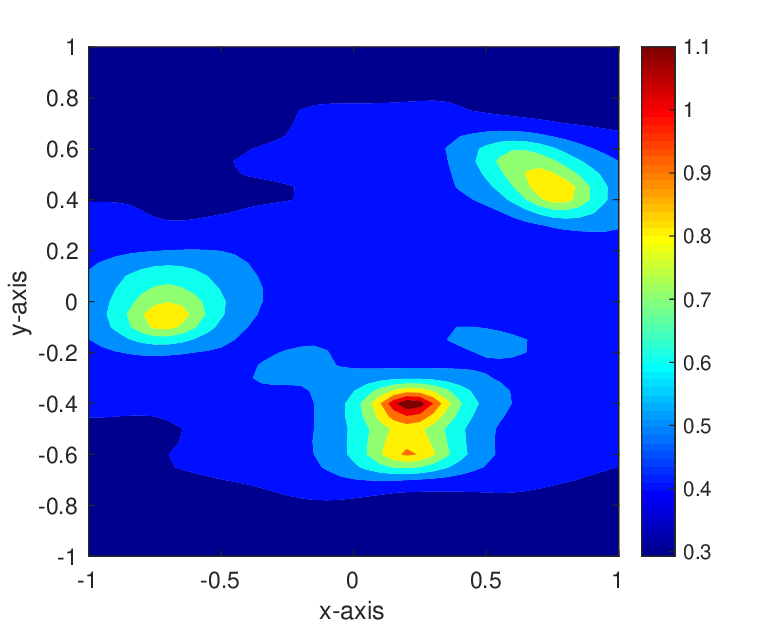}}
\caption{\label{Result4-1}(Example \ref{MU2}) Distribution of singular values (first and third columns) and maps of $f_{\music}(\mr)$ for $k=2\pi/0.4$ (second and fourth columns).}
\end{figure}

\begin{figure}
\centering
\subfigure[Case 5]{\includegraphics[width=.250\textwidth]{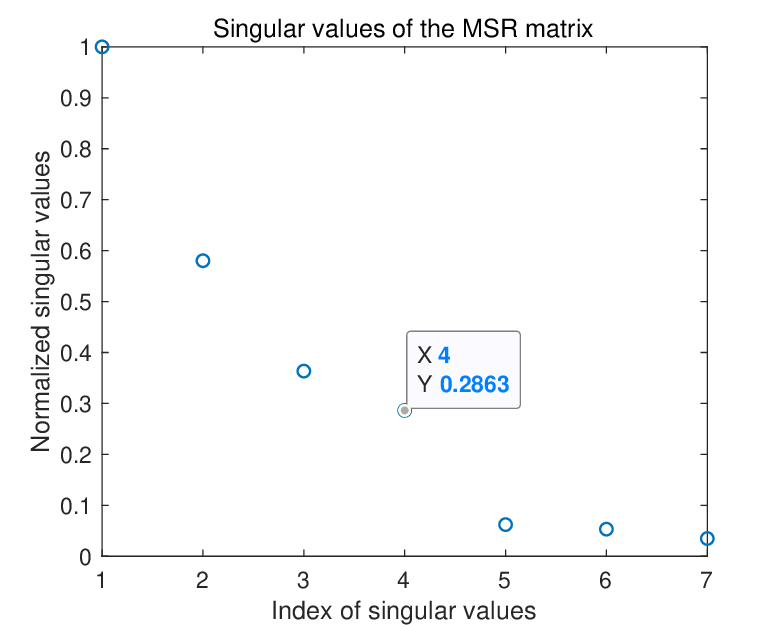}\hfill\includegraphics[width=.250\textwidth]{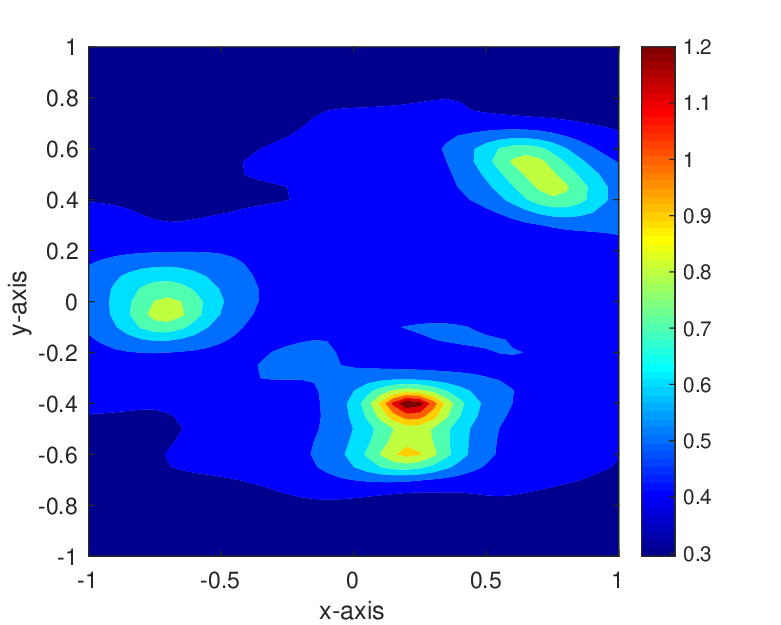}}\hfill
\subfigure[Case 6]{\includegraphics[width=.250\textwidth]{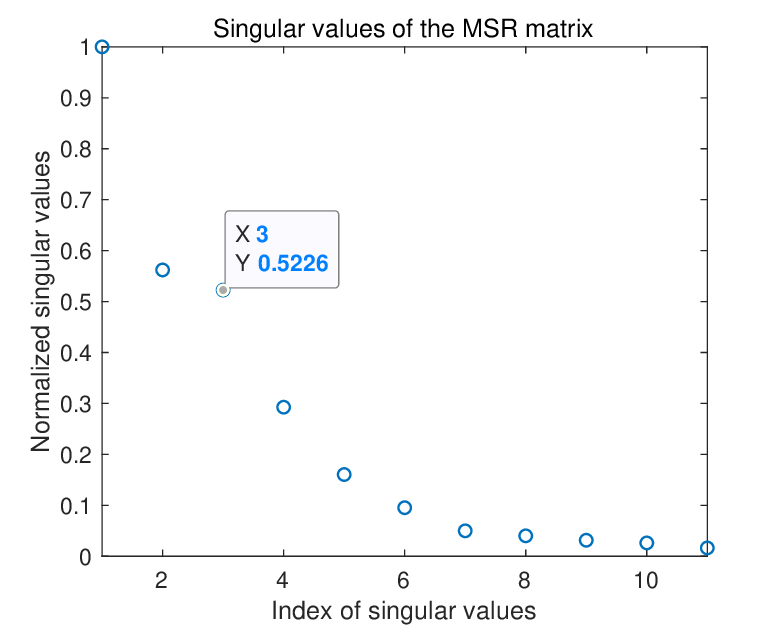}\hfill\includegraphics[width=.250\textwidth]{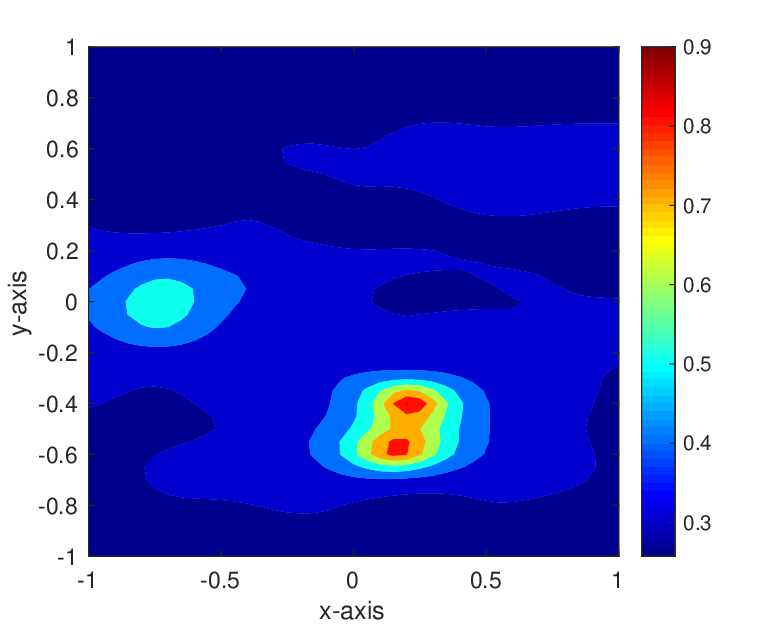}}\\
\subfigure[Case 7]{\includegraphics[width=.250\textwidth]{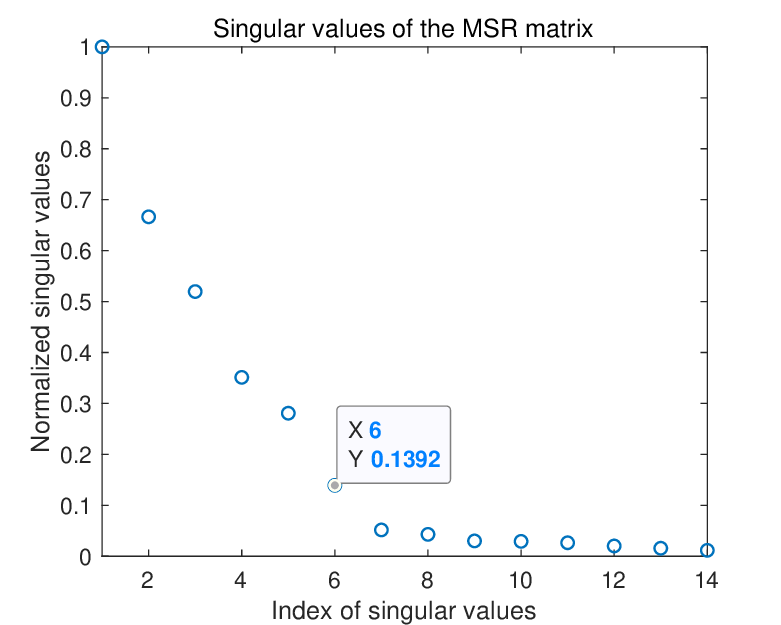}\hfill\includegraphics[width=.250\textwidth]{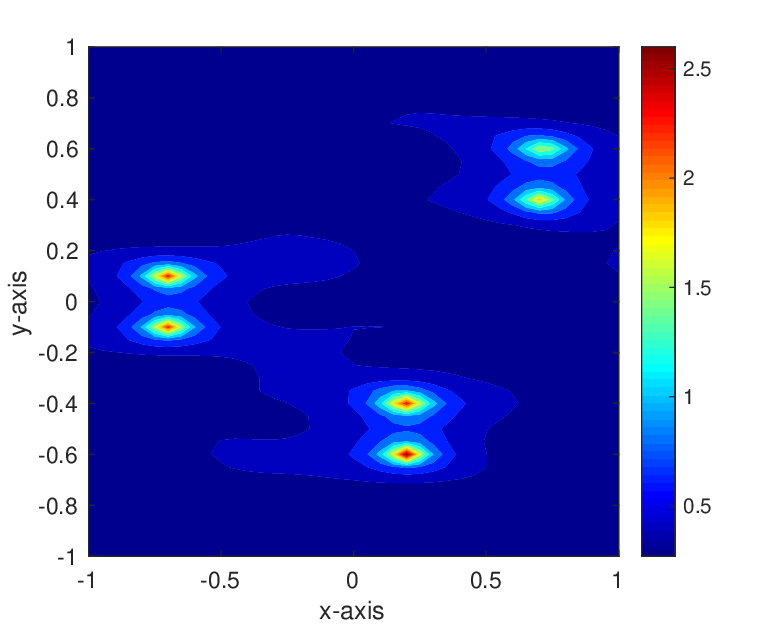}}\hfill
\subfigure[Case 8]{\includegraphics[width=.250\textwidth]{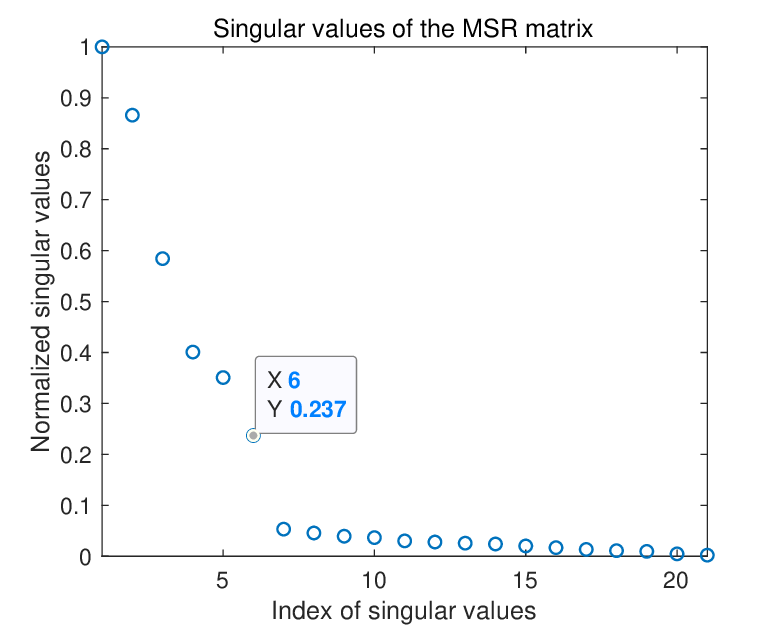}\hfill\includegraphics[width=.250\textwidth]{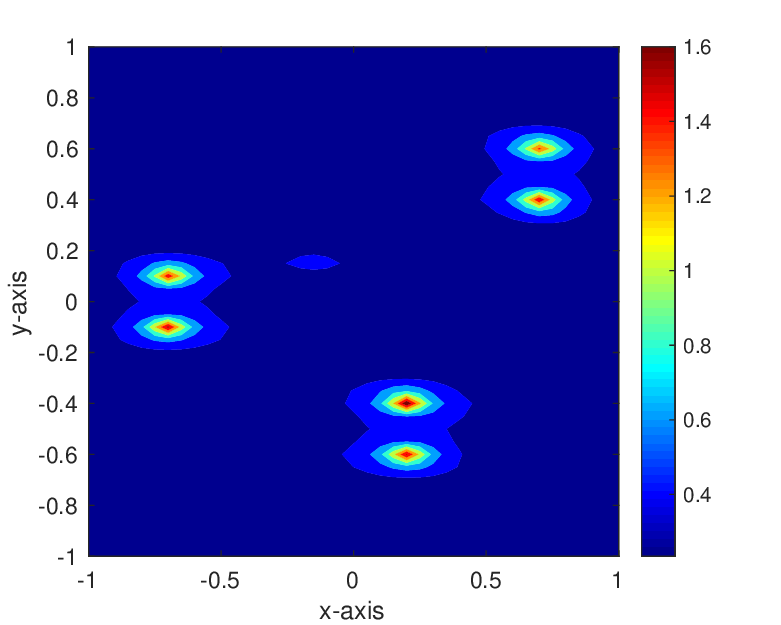}}
\caption{\label{Result4-2}(Example \ref{MU2}) Distribution of singular values (first and third columns) and maps of $f_{\music}(\mr)$ for $k=2\pi/0.4$ (second and fourth columns).}
\end{figure}

\begin{discuss}[Feasibility and limitation of MUSIC]\label{Discussion-M1}
The imaging function $f_{\music}(\mr)$ consists of
\[iJ_1(k|\mr-\mr_s|)\bigg(\frac{\mr-\mr_s}{|\mr-\mr_s|}\cdot\me_h\bigg),\quad\frac{\Lambda_{\mu}^{(1,h)}(\mr-\mr_s)}{\vartheta_M-\vartheta_1},\quad\text{and}\quad\frac{\Lambda_{\mu}^{(2,h)}(\mr-\mr_s)}{\theta_N-\theta_1}\]
for $h=1,2$. The first term is independent of the range of incident and observation directions but does not contribute to the detection because $J_1(0)=0$, whereas the remaining terms are dependent. It is interesting to observe that, since $\Lambda_{\mu}^{(1,h)}(\mr-\mr_s)$ and $\Lambda_{\mu}^{(2,h)}(\mr-\mr_s)$ contain the factors $J_0(k|\mr-\mr_s|)$ and $J_0(0)=1$, it is possible to identify $\mr_s$ via the map of $f_{\music}(\mr)$ if the range of incident or observation directions is narrow. However, if $\mr\ne\mr_s$ and the range of incident or observation directions is narrow, refer to Figures \ref{Result3-1} and \ref{Result4-1}. Otherwise, if the range is sufficiently wide, $f_{\music}(\mr)$ will be dominated by the first term containing $J_1(k|\mr-\mr_s|)$, so that two large-magnitude peaks will appear in the neighborhood of $\mr_s$; refer to Figures \ref{Result3-2} and \ref{Result4-2}. Although this result does not yield, one cannot identify a true location with this result, based on the property of $J_1$.
\end{discuss}

\begin{discuss}[Least range of directions]\label{Discussion-M2}
Now, let us eliminate the terms that are dependent on the range of directions, i.e.,
\[\frac{\Lambda_{\mu}^{(1,h)}(\mr-\mr_s)}{\vartheta_M-\vartheta_1}=0\quad\text{and}\quad\frac{\Lambda_{\mu}^{(2,h)}(\mr-\mr_s)}{\theta_N-\theta_1}=0\]
for all $h=1,2$. One possible choice is to select $\vartheta_1$, $\vartheta_M$, $\theta_1$, and $\theta_N$, thereby satisfying
\[\left\{\begin{array}{l}
\displaystyle\medskip\sin\frac{\vartheta_M-\vartheta_1}{2}\cos\frac{\vartheta_M+\vartheta_1}{2}=0,\quad \sin(\vartheta_M-\vartheta_1)\cos(\vartheta_M+\vartheta_1-\phi_s)=0,\\
\displaystyle\medskip\sin\frac{\vartheta_M-\vartheta_1}{2}\cos\frac{\vartheta_M+\vartheta_1-\phi_s}{2}=0,\quad\sin(\theta_N-\theta_1)\cos(\theta_N+\theta_1-\phi_s)=0,\\
\displaystyle\medskip\sin\frac{(1-p)(\vartheta_M-\vartheta_1)}{2}\cos\frac{(1-p)(\vartheta_M+\vartheta_1)+2p\phi_s-2p\pi}{2}=0,\\
\displaystyle\medskip\sin\frac{(1+n)(\vartheta_M-\vartheta_1)}{2}\cos\frac{(1+p)(\vartheta_M+\vartheta_1)-2p\phi_s+2p\pi}{2}=0,\\
\displaystyle\medskip\sin\frac{(1-p)(\theta_N-\theta_1)}{2}\cos\frac{(1-p)(\theta_N+\theta_1)+2p\phi_s-(2p+1)\pi}{2}=0,\\
\displaystyle\sin\frac{(1+n)(\theta_N-\theta_1)}{2}\cos\frac{(1+p)(\theta_N+\theta_1)-2p\phi_s+(2p-1)\pi}{2}=0,
\end{array}\right.\]
{for every positive integer $p$ and $q$}. Similar to the permittivity contrast case, another possible choice is to select $\vartheta_1=\phi_s$, $\vartheta_M=\pi+\phi_s$, $\theta_1=\phi_s$, and $\theta_N=\pi+\phi_s$. Therefore, if the range of incident and observation directions is wider than $\pi$, two large-magnitude peaks will appear in the neighborhood of $\mr_s$, and correspondingly, it will be possible to identify the location of $\mr_s$, which is the middle point of two peaks; refer to Figures \ref{Result3-2} and \ref{Result4-2}.
\end{discuss}

\section{Simulation results}\label{sec:4}
\begin{figure}
\centering
\subfigure[Case 1]{\includegraphics[width=.250\textwidth]{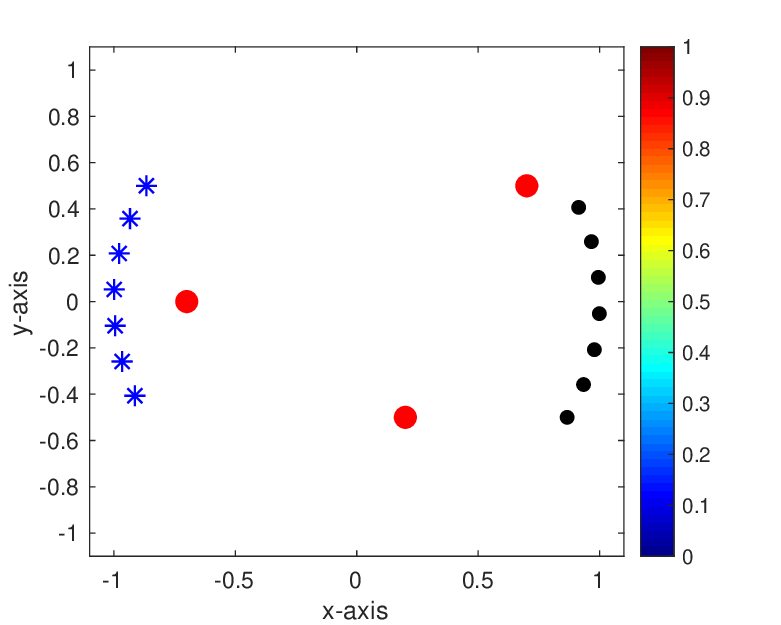}}\hfill
\subfigure[Case 2]{\includegraphics[width=.250\textwidth]{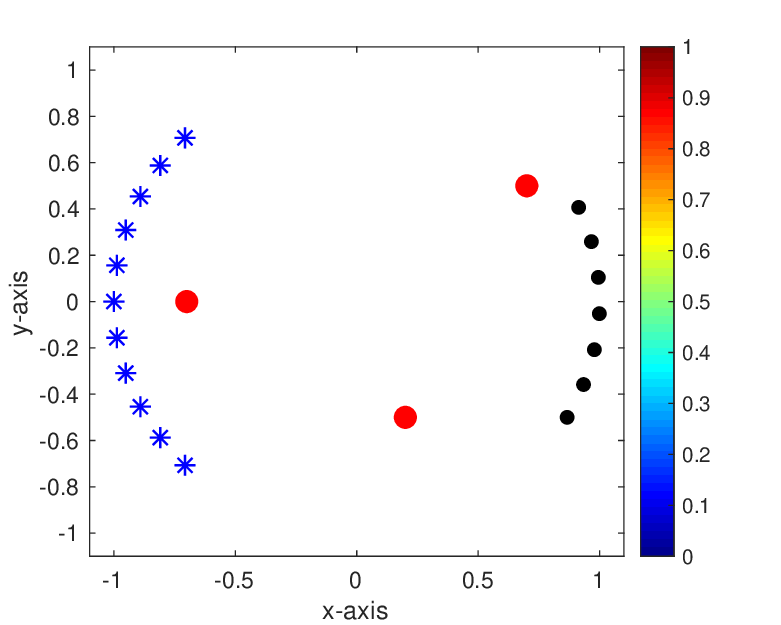}}\hfill
\subfigure[Case 3]{\includegraphics[width=.250\textwidth]{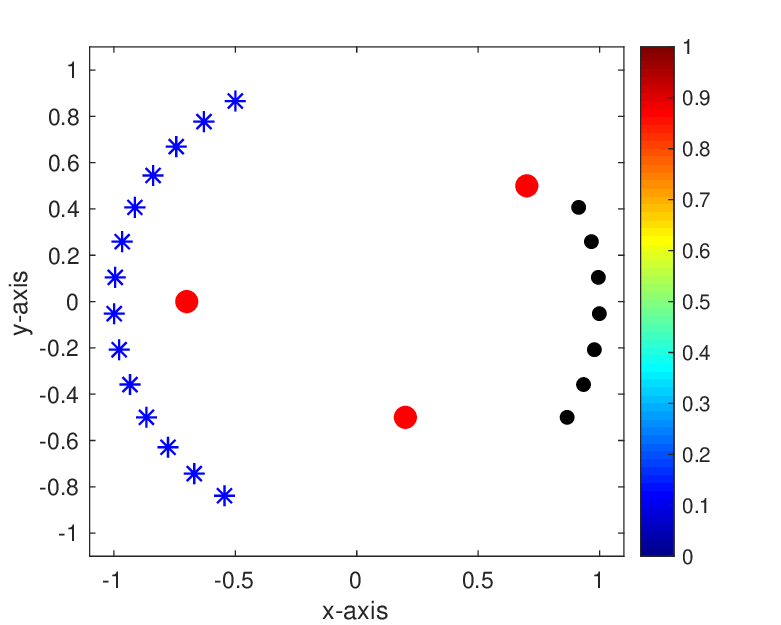}}\hfill
\subfigure[Case 4]{\includegraphics[width=.250\textwidth]{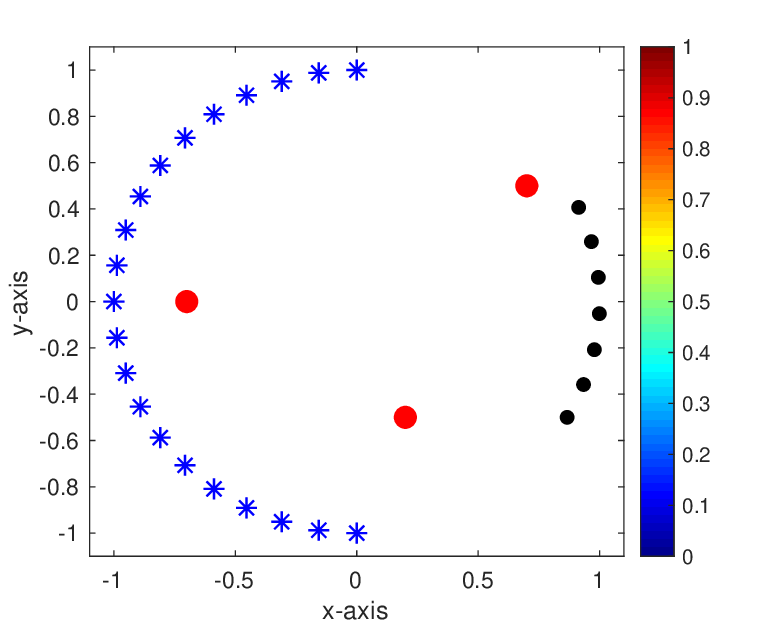}}\\
\subfigure[Case 5]{\includegraphics[width=.250\textwidth]{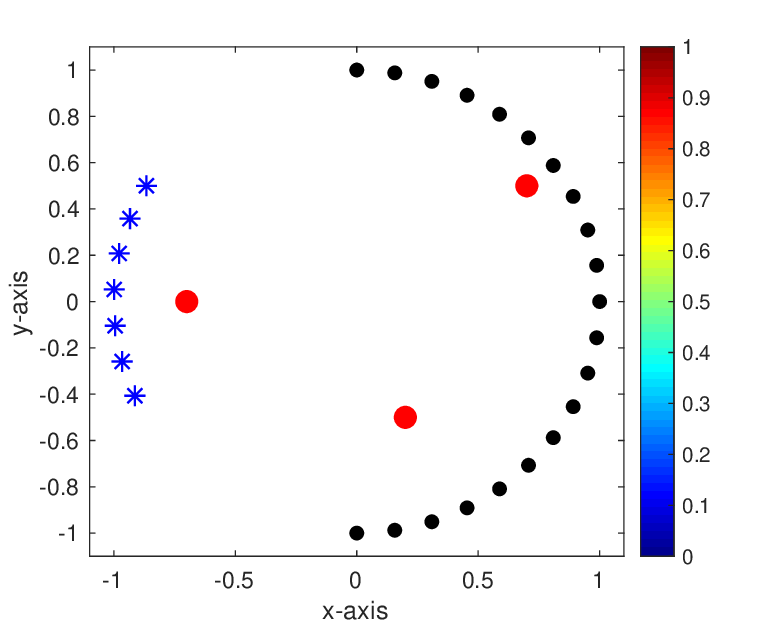}}\hfill
\subfigure[Case 6]{\includegraphics[width=.250\textwidth]{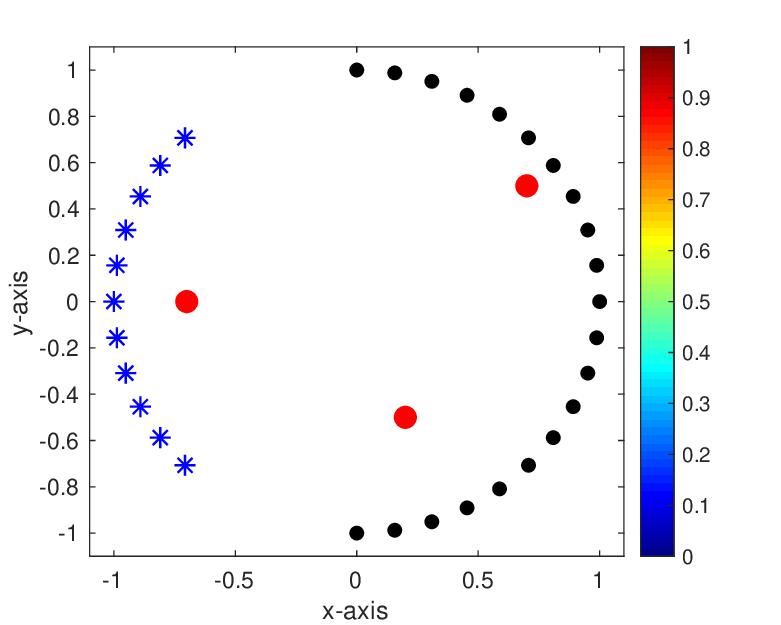}}\hfill
\subfigure[Case 7]{\includegraphics[width=.250\textwidth]{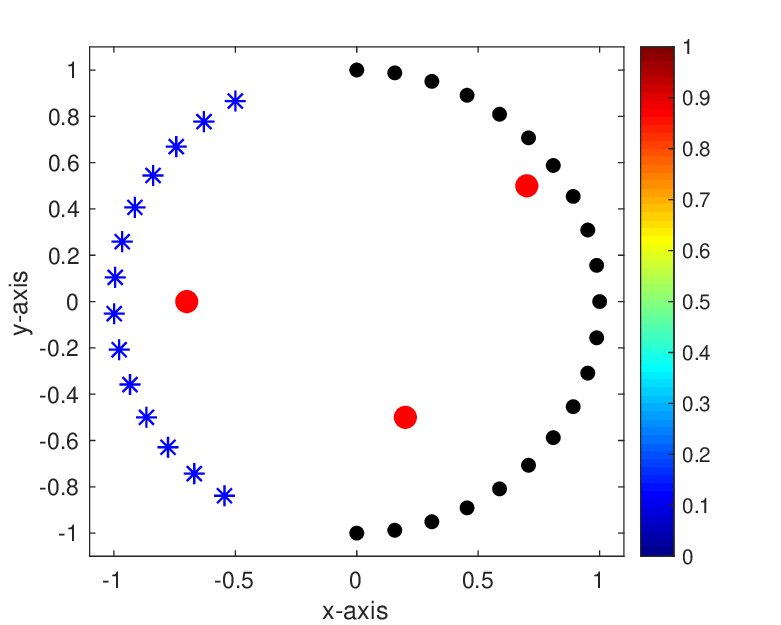}}\hfill
\subfigure[Case 8]{\includegraphics[width=.250\textwidth]{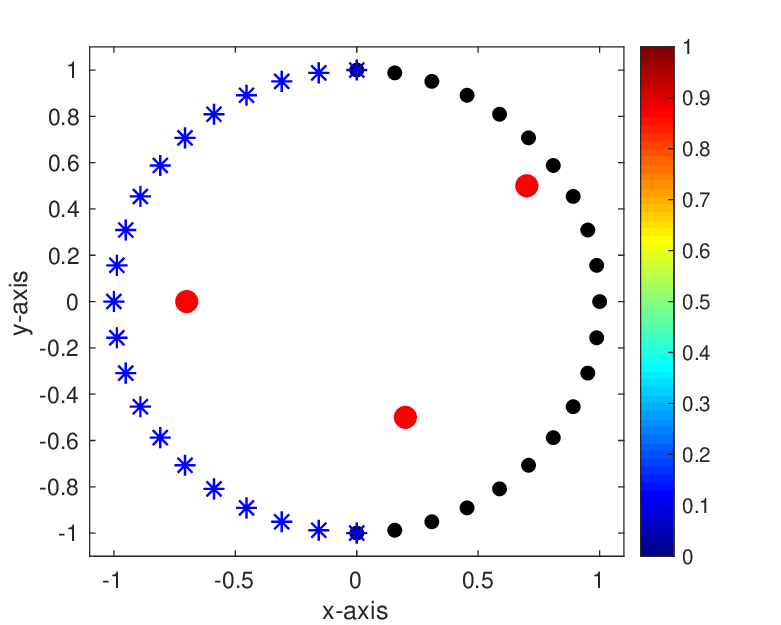}}
\caption{\label{Setting}Illustration of the configuration range. Red-colored circles denote the inhomogeneities $\Sigma_s$, $s=1,2,3$, and black-colored circles and blue-colored stars denote the incident and observation directions, respectively.}
\end{figure}

In this section, we exhibit a set of simulation results to validate the investigated results of Theorem \ref{Theorem1} and Theorem \ref{Theorem2}. We select $S=3$ small circles with the same radii {$\alpha_s\equiv0.1$}, permittivity $\eps_s$, and permeability $\mu_s$. The locations $\mr_s$ are $\mr_1=[0.7,0.5]^T$, $\mr_2=[-0.7,0.0]^T$, and $\mr_3=[0.2,-0.5]^T$, and the background permittivity and permeability are $\epsb=\mub=1$. The far-field pattern data $u_{\infty}(\vv_j,\vt_l)$ of $\mathbb{K}$ at wavelength $k=2\pi/0.4$ are generated by solving the Foldy-Lax formulation introduced in \cite{HSZ4} to avoid an \textit{inverse crime}. After the generation of the far-field pattern, $\SI{20}{\deci\bel}$ white Gaussian random noise is added to the unperturbed data through the MATLAB command \texttt{awgn} included in the signal processing package. For incident and observation direction setting, eight different configurations are chosen, as shown in Figure \ref{Setting}.

\begin{figure}
\centering
\subfigure[Case 1]{\includegraphics[width=.250\textwidth]{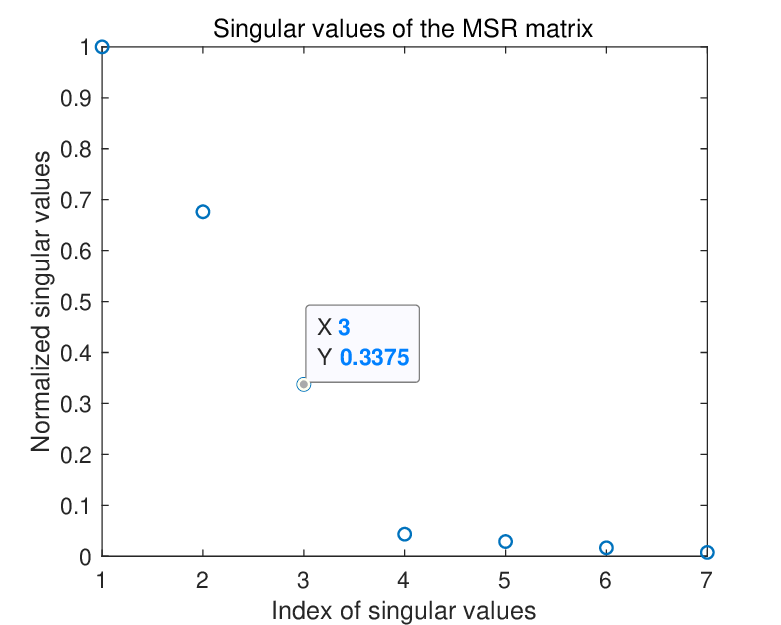}\hfill\includegraphics[width=.250\textwidth]{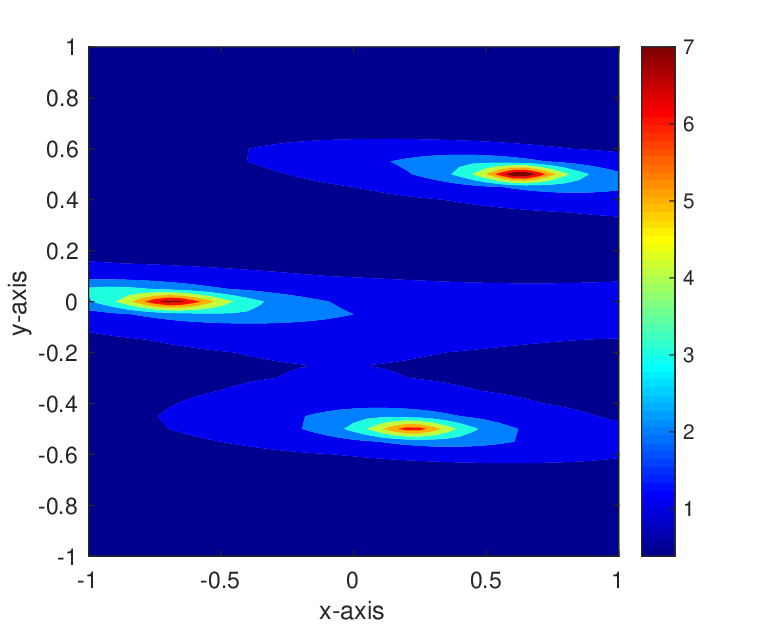}}\hfill
\subfigure[Case 2]{\includegraphics[width=.250\textwidth]{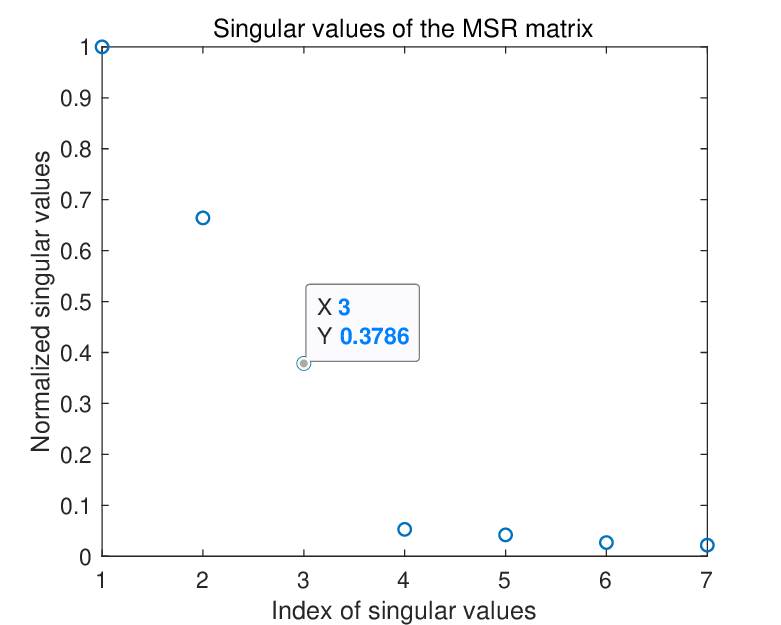}\hfill\includegraphics[width=.250\textwidth]{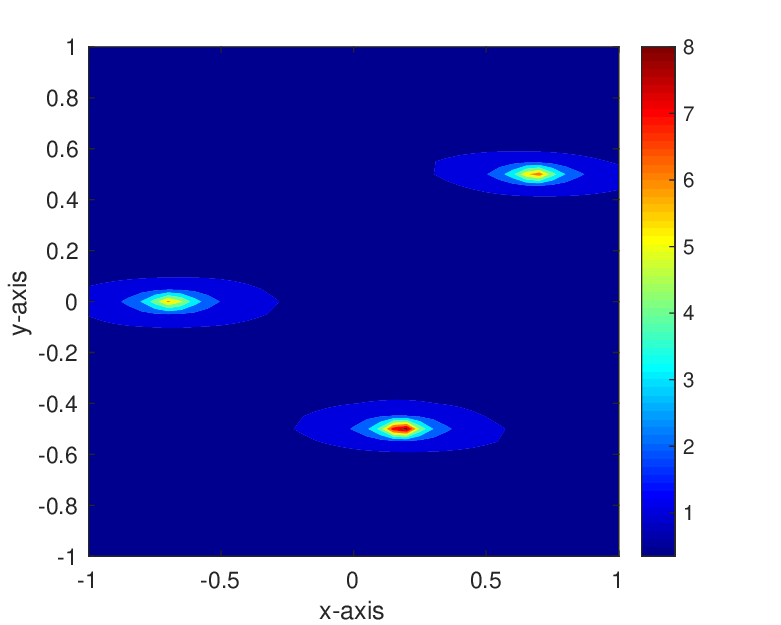}}\\
\subfigure[Case 3]{\includegraphics[width=.250\textwidth]{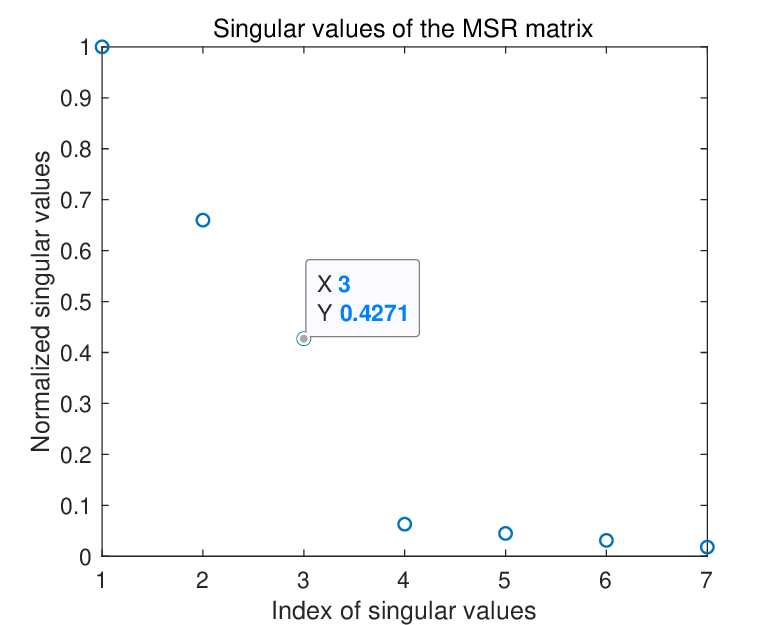}\hfill\includegraphics[width=.250\textwidth]{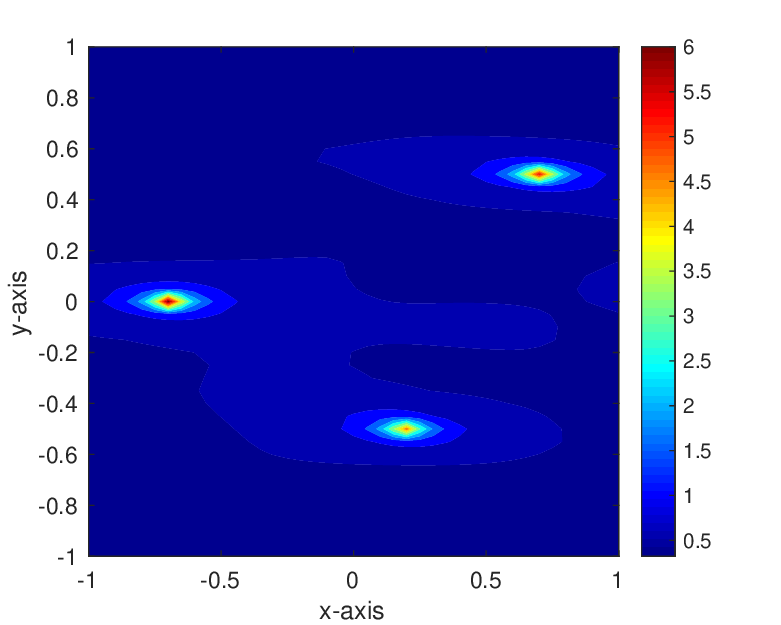}}\hfill
\subfigure[Case 4]{\includegraphics[width=.250\textwidth]{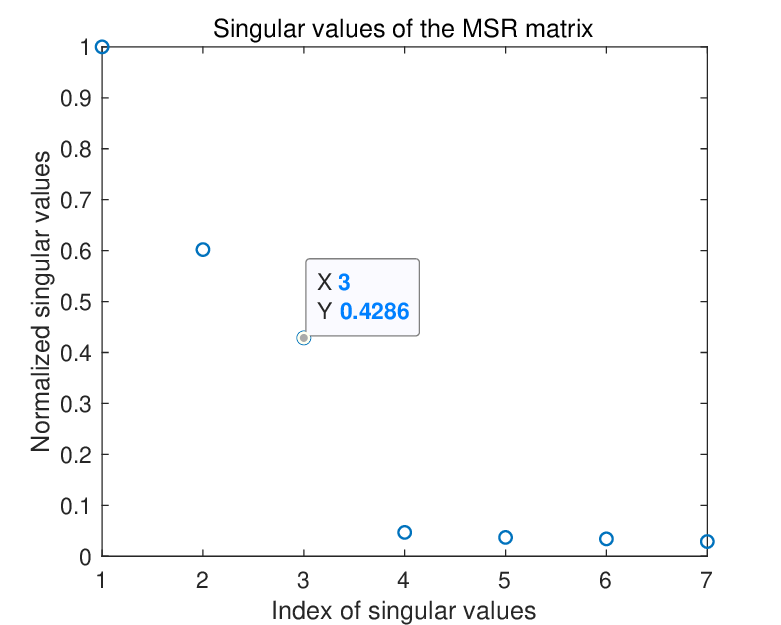}\hfill\includegraphics[width=.250\textwidth]{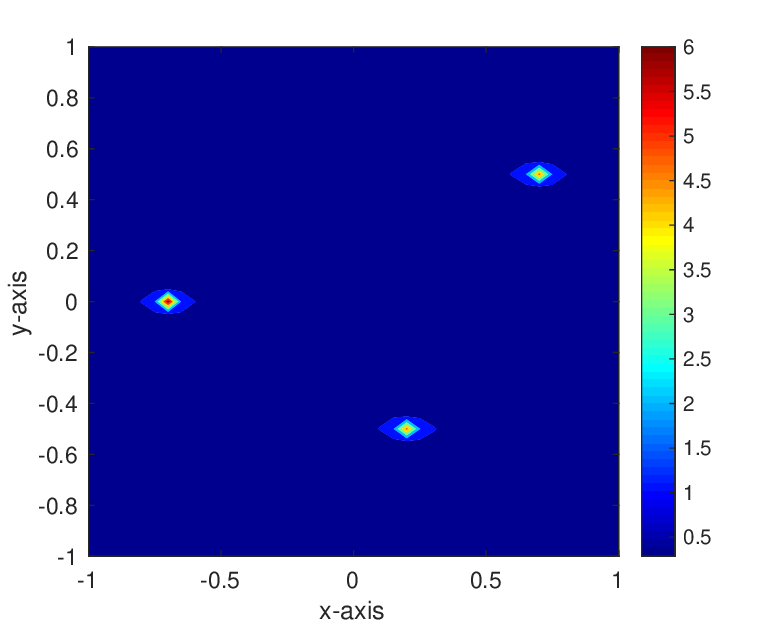}}
\caption{\label{Result1-1}(Example \ref{EPS1}) Distribution of singular values (first and third columns) and maps of $f_{\music}(\mr)$ for $k=2\pi/0.4$ (second and fourth columns).}
\end{figure}

\begin{figure}
\centering
\subfigure[Case 5]{\includegraphics[width=.250\textwidth]{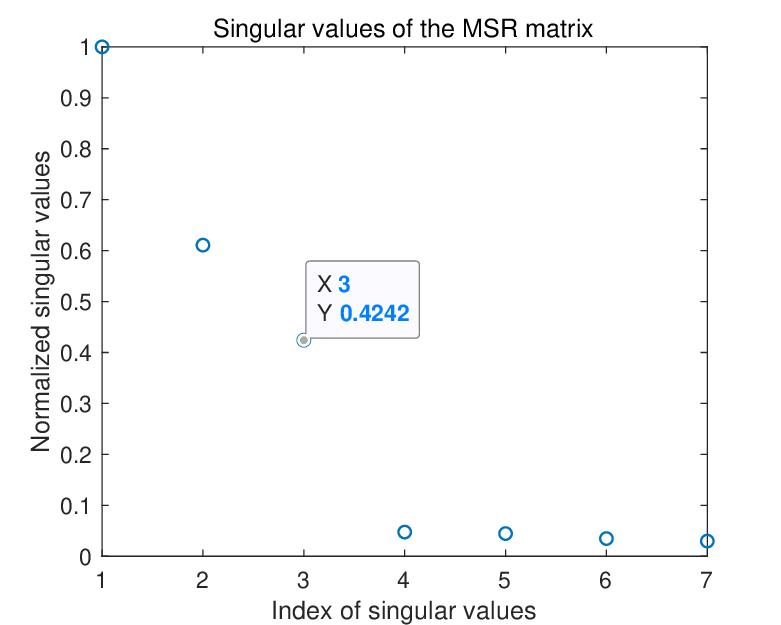}\hfill\includegraphics[width=.250\textwidth]{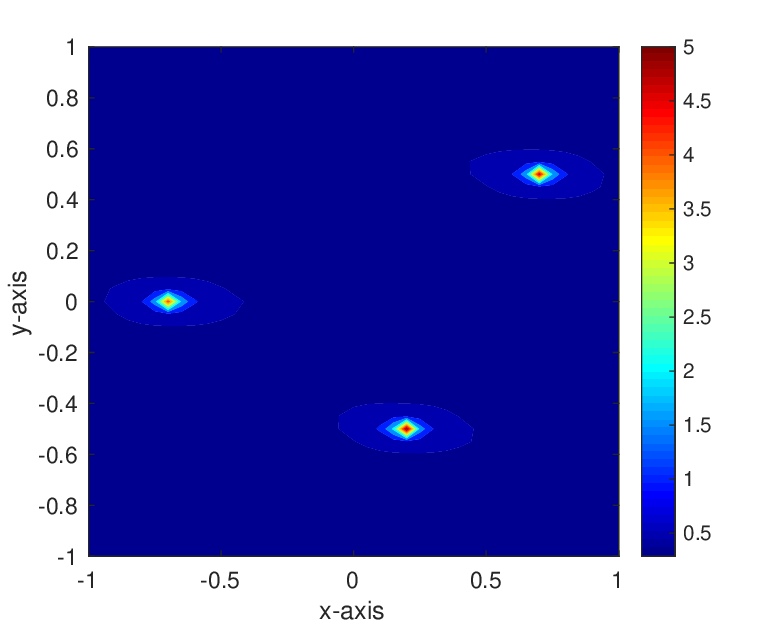}}\hfill
\subfigure[Case 6]{\includegraphics[width=.250\textwidth]{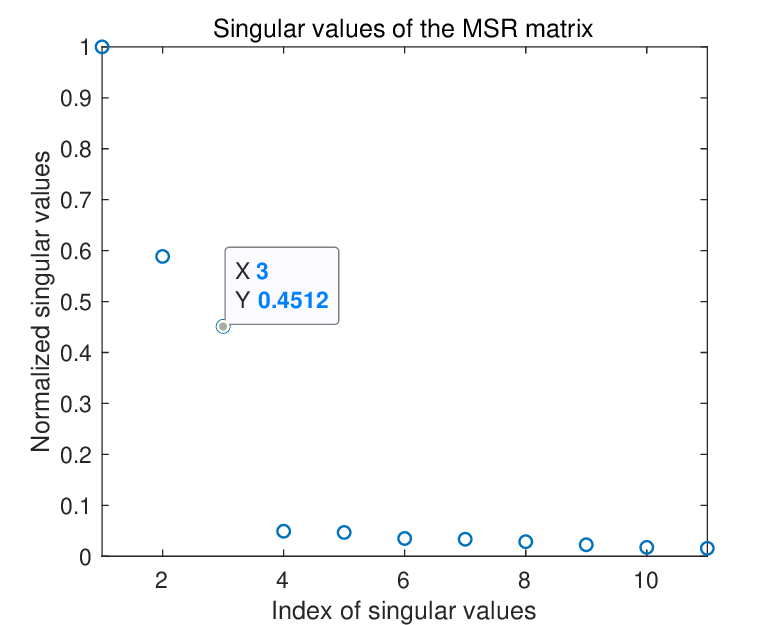}\hfill\includegraphics[width=.250\textwidth]{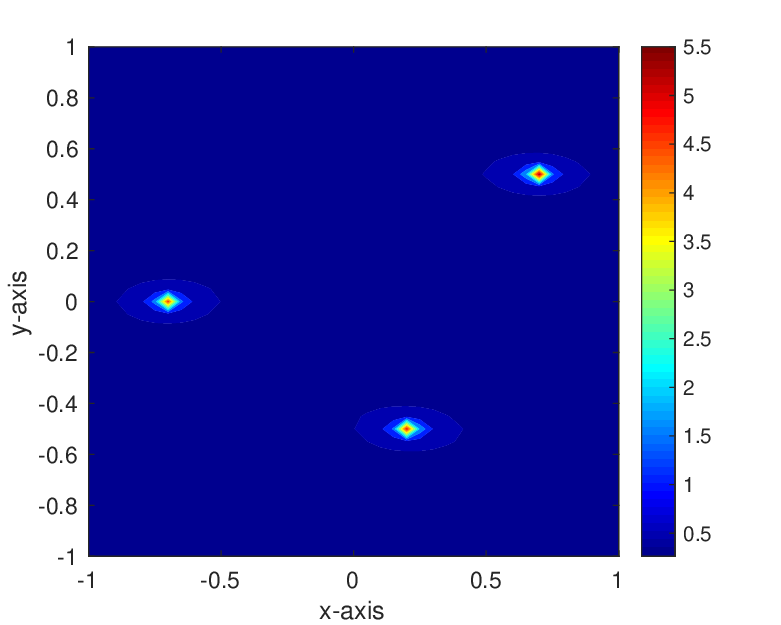}}\\
\subfigure[Case 7]{\includegraphics[width=.250\textwidth]{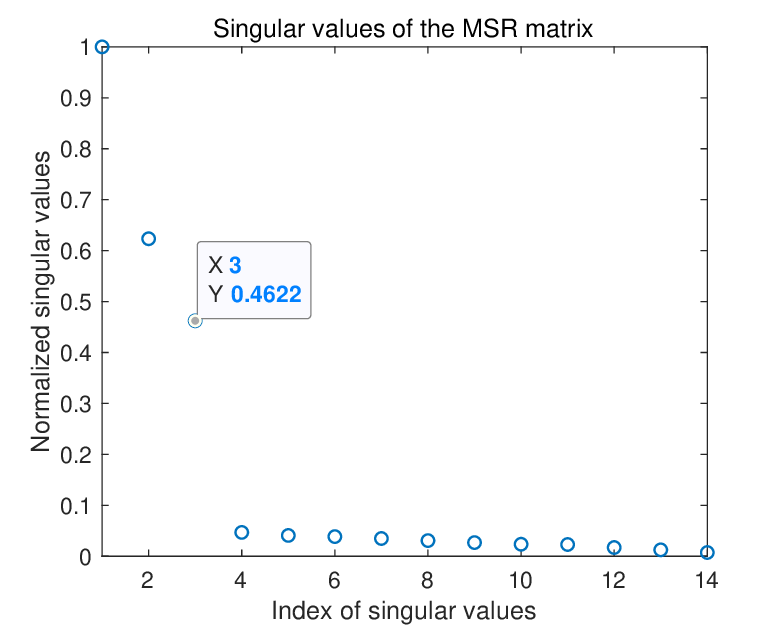}\hfill\includegraphics[width=.250\textwidth]{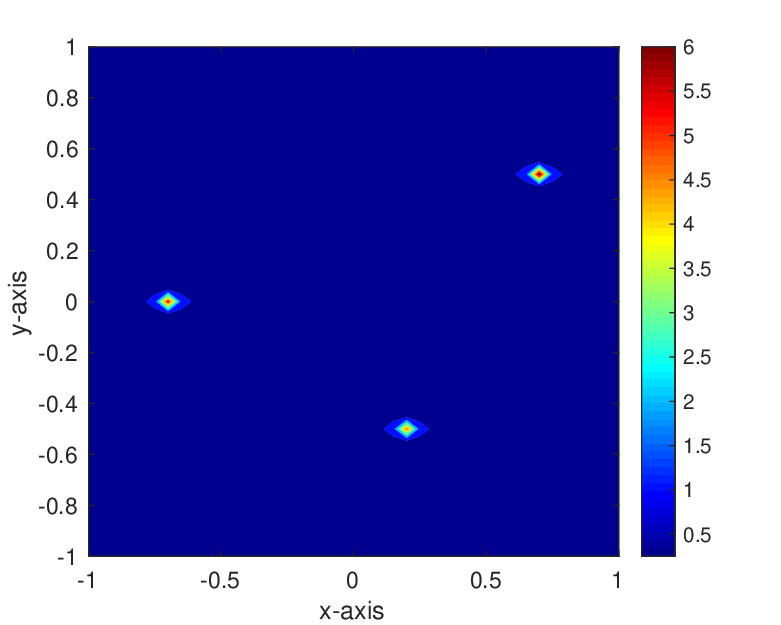}}\hfill
\subfigure[Case 8]{\includegraphics[width=.250\textwidth]{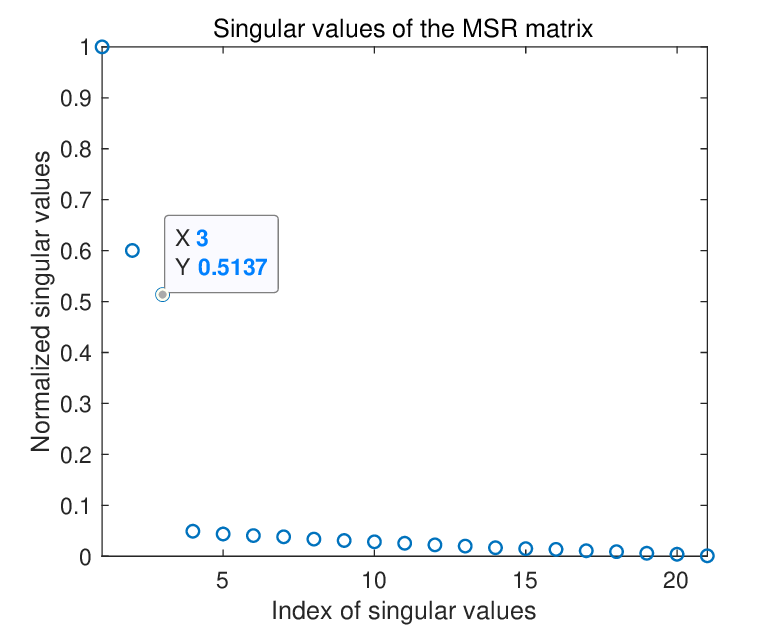}\hfill\includegraphics[width=.250\textwidth]{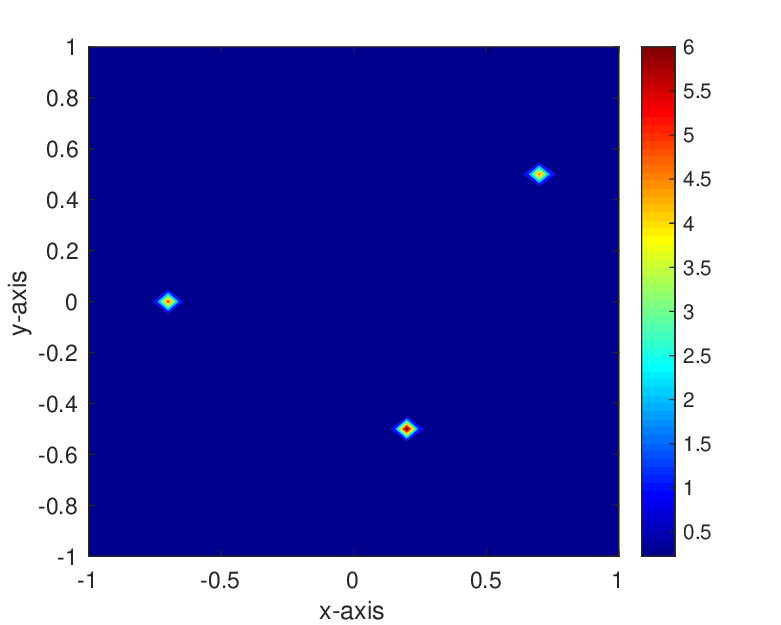}}
\caption{\label{Result1-2}(Example \ref{EPS1}) Distribution of singular values (first and third columns) and maps of $f_{\music}(\mr)$ for $k=2\pi/0.4$ (second and fourth columns).}
\end{figure}

\begin{ex}[Permittivity contrast case: same permittivity]\label{EPS1}
Figure \ref{Result1-1} shows the distribution of singular values of $\mathbb{K}$ and the maps of $f_{\music}(\mr)$ for Cases 1-4 when $\eps_s\equiv5$ and $\mu_s=\mub$. Based on the results, we can select three singular values for generating the projection operator, and every location $\mr_s$ can be identified clearly. Notice that, for Case 1, due to the appearance of the blurring effect in the neighborhood of $\mr_s$, identifying the exact location of $\Sigma_s$ is difficult, but it is possible to recognize its existence.

Figure \ref{Result1-2} shows the distribution of singular values of $\mathbb{K}$ and the maps of $f_{\music}(\mr)$ for Cases 5-8. In contrast to the previous results, since the range of incident direction is $\pi$, the location of $\Sigma_s$ can be identified clearly. This result partially supports Property \ref{Discussion-E2}.
\end{ex}

\begin{figure}
\centering
\subfigure[Case 1]{\includegraphics[width=.250\textwidth]{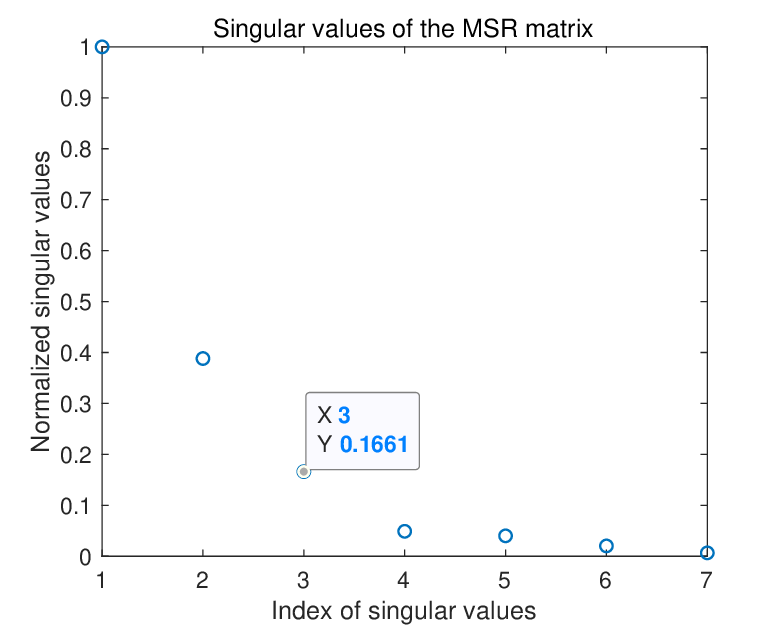}\hfill\includegraphics[width=.250\textwidth]{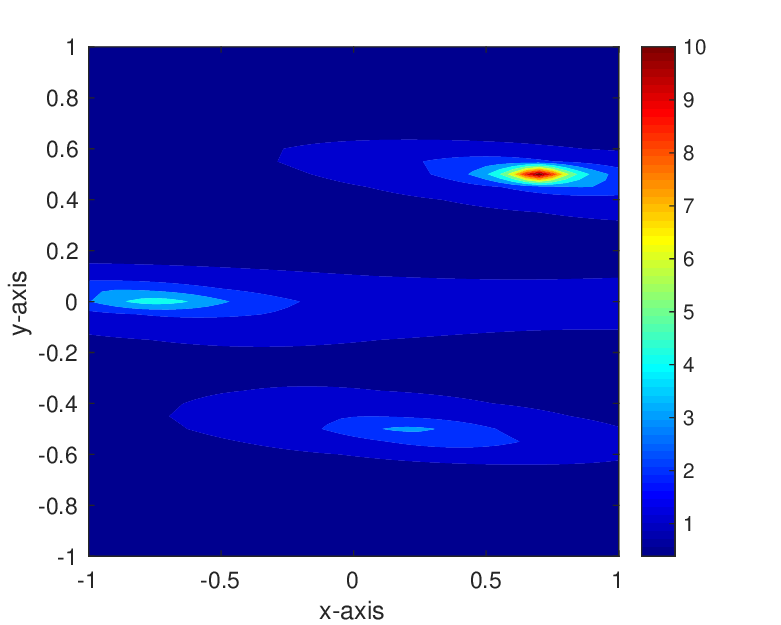}}\hfill
\subfigure[Case 2]{\includegraphics[width=.250\textwidth]{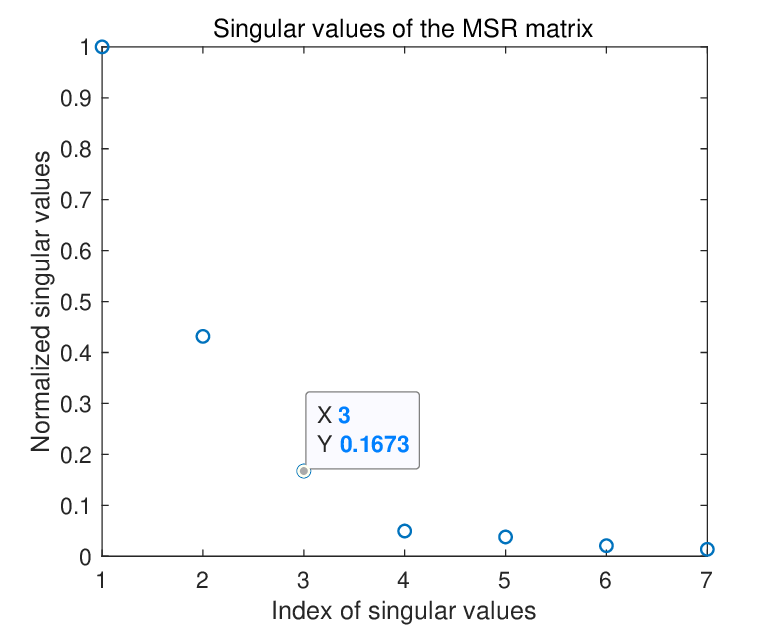}\hfill\includegraphics[width=.250\textwidth]{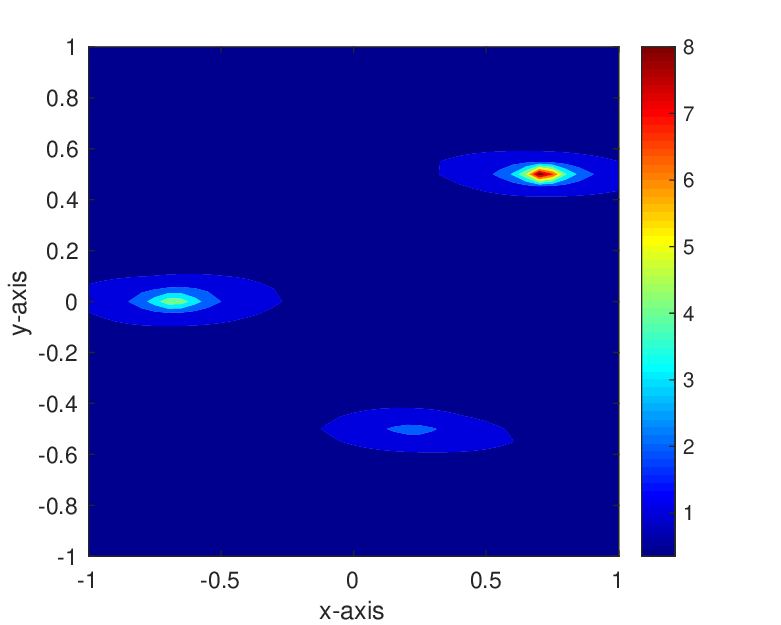}}\\
\subfigure[Case 3]{\includegraphics[width=.250\textwidth]{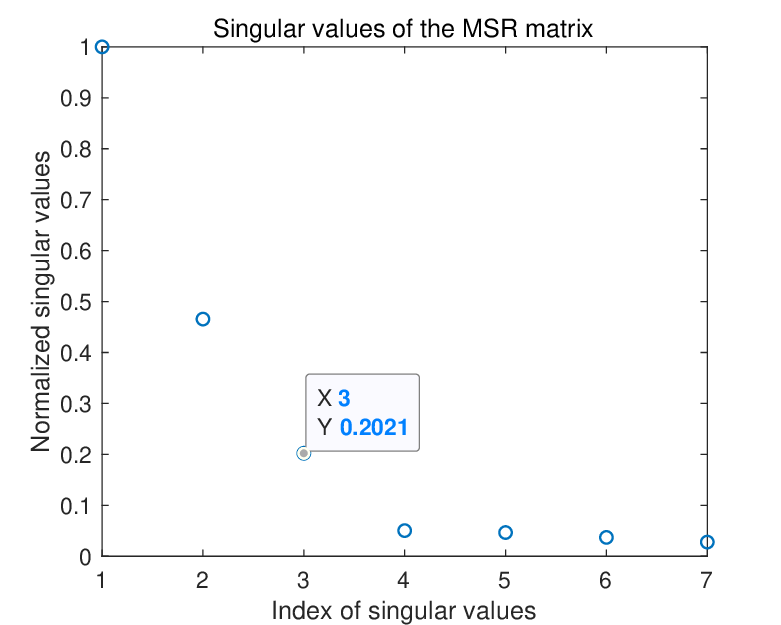}\hfill\includegraphics[width=.250\textwidth]{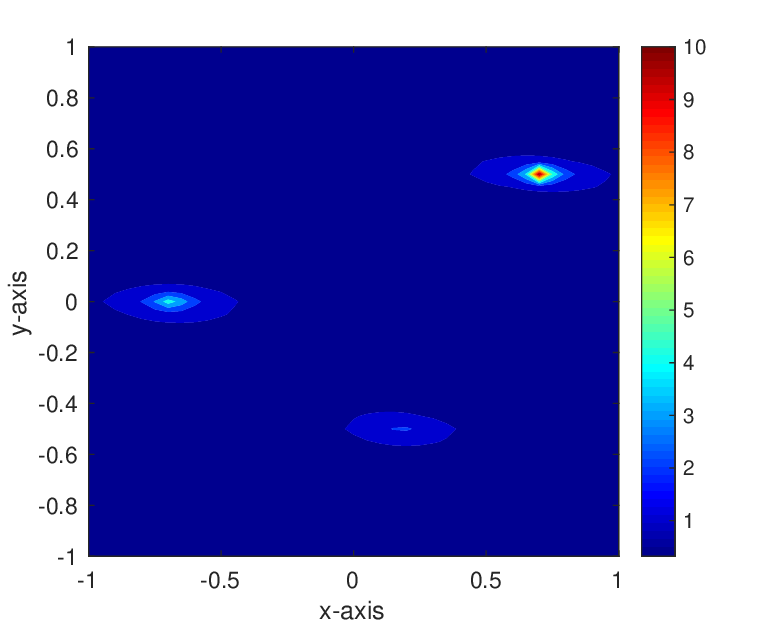}}\hfill
\subfigure[Case 4]{\includegraphics[width=.250\textwidth]{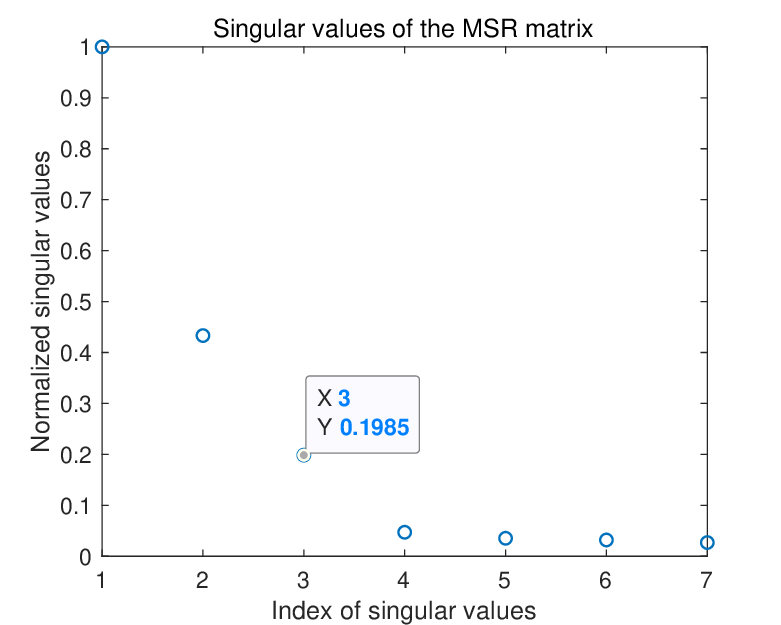}\hfill\includegraphics[width=.250\textwidth]{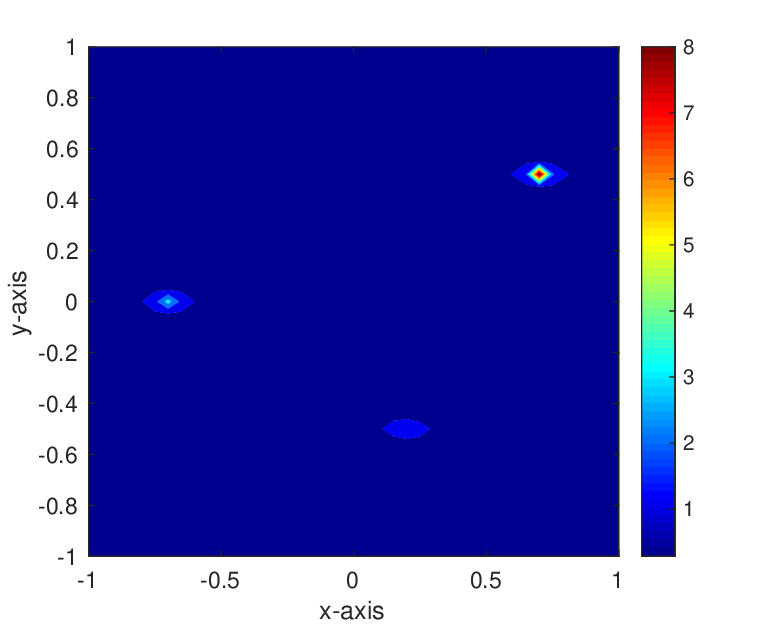}}
\caption{\label{Result2-1}(Example \ref{EPS2}) Distribution of singular values (first and third columns) and maps of $f_{\music}(\mr)$ for $k=2\pi/0.4$ (second and fourth columns).}
\end{figure}

\begin{figure}
\centering
\subfigure[Case 5]{\includegraphics[width=.250\textwidth]{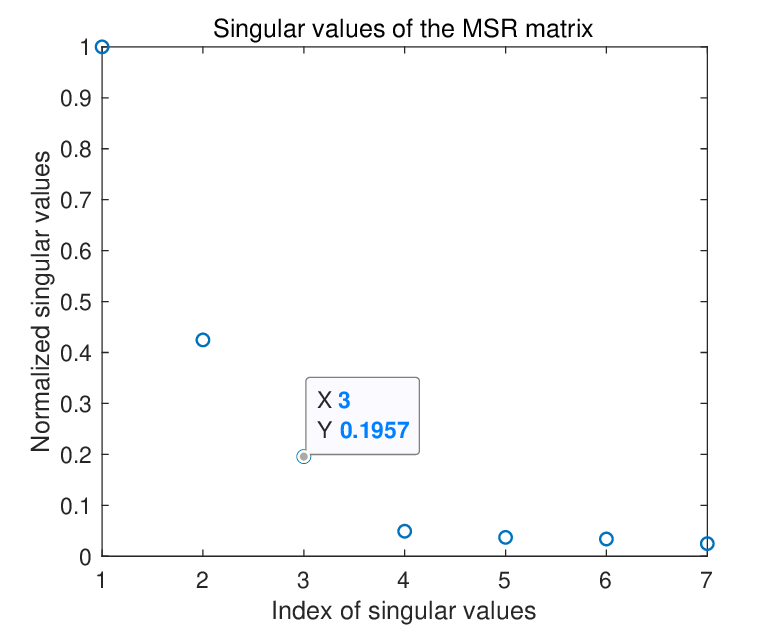}\hfill\includegraphics[width=.250\textwidth]{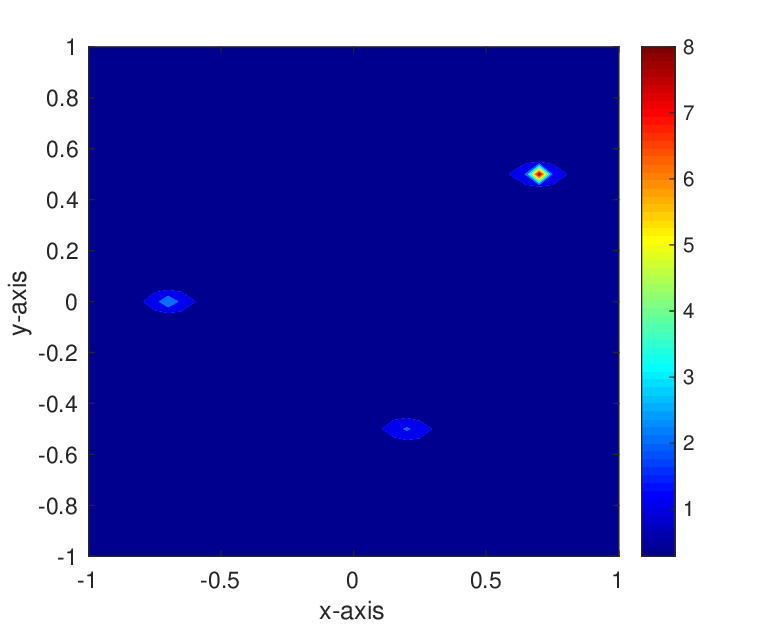}}\hfill
\subfigure[Case 6]{\includegraphics[width=.250\textwidth]{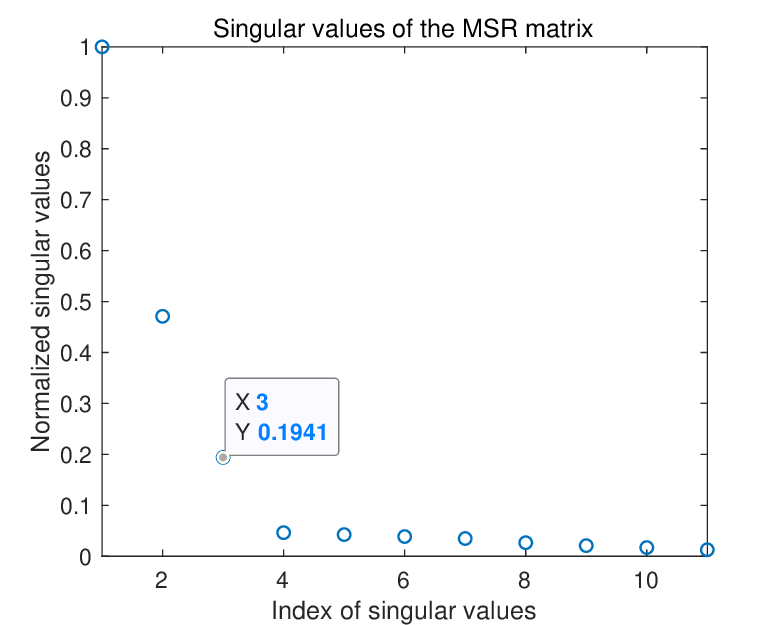}\hfill\includegraphics[width=.250\textwidth]{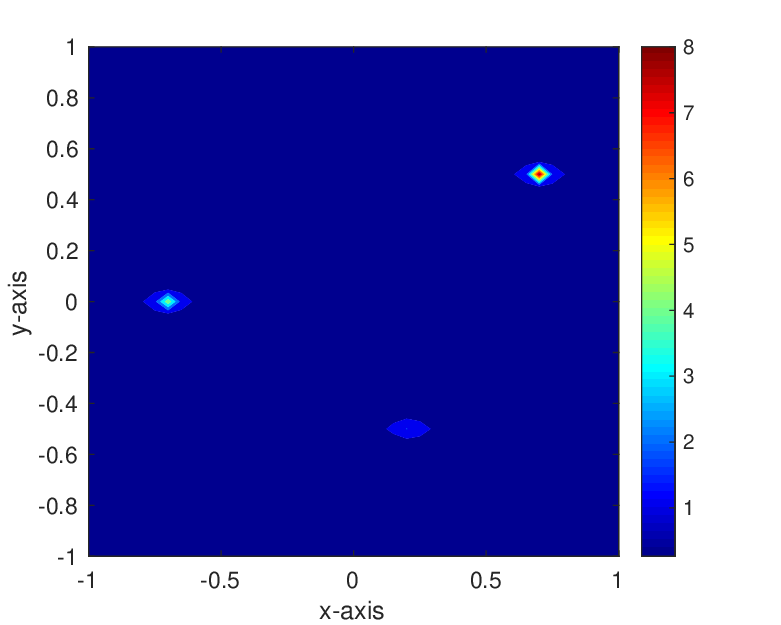}}\\
\subfigure[Case 7]{\includegraphics[width=.250\textwidth]{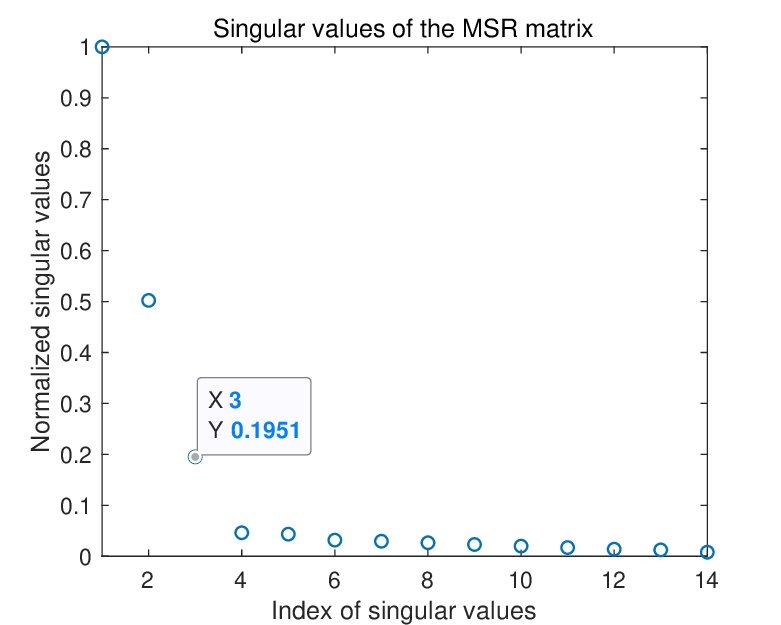}\hfill\includegraphics[width=.250\textwidth]{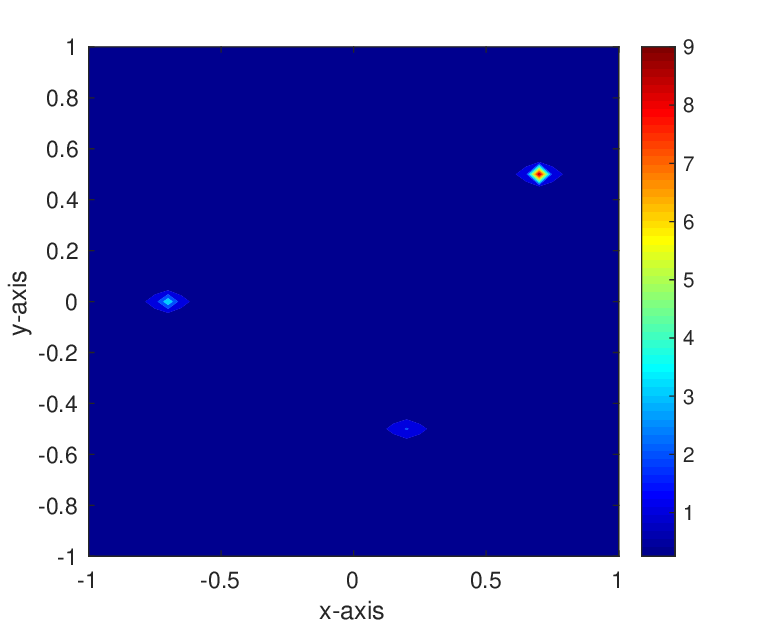}}\hfill
\subfigure[Case 8]{\includegraphics[width=.250\textwidth]{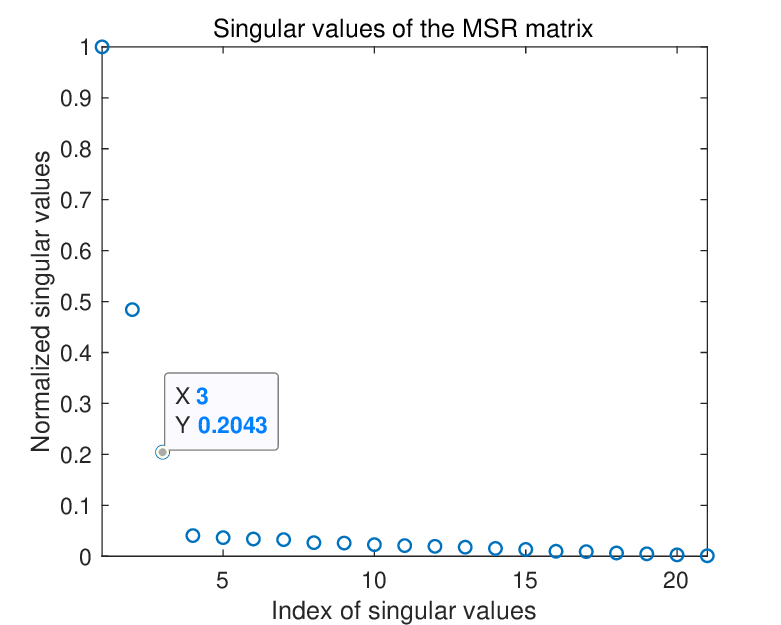}\hfill\includegraphics[width=.250\textwidth]{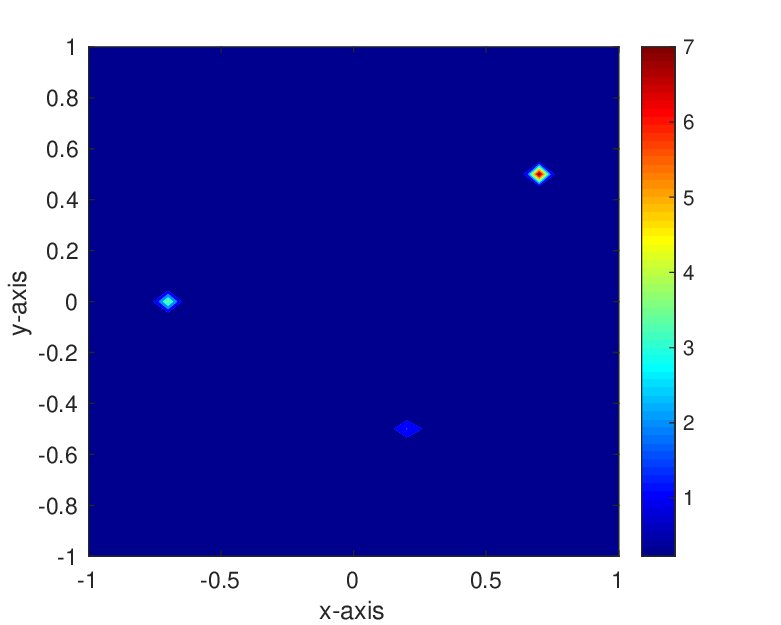}}
\caption{\label{Result2-2}(Example \ref{EPS2}) Distribution of singular values (first and third columns) and maps of $f_{\music}(\mr)$ for $k=2\pi/0.4$ (second and fourth columns).}
\end{figure}

\begin{ex}[Permittivity contrast case: different permittivities]\label{EPS2}
Here, we consider the imaging results for different permittivities. For this, we set $\eps_1=5$, $\eps_2=3$, $\eps_3=2$, and $\mu_s=\mub$. Figures \ref{Result2-1} and \ref{Result2-2} show the distribution of singular values of $\mathbb{K}$ and the maps of $f_{\music}(\mr)$ for Cases 1-4. Compared with the results in Example \ref{EPS1}, the selection of three nonzero singular values is more difficult. Moreover, since $\eps_3<\eps_1,\eps_2$, the value of $f_{\music}(\mr_3)$ is significantly smaller than the values of $f_{\music}(\mr_1)$ and $f_{\music}(\mr_2)$, which means it is possible to recognize the existence of $\Sigma_3$, but its exact location cannot be identified.
\end{ex}

\begin{figure}
\centering
\subfigure[Case 1]{\includegraphics[width=.250\textwidth]{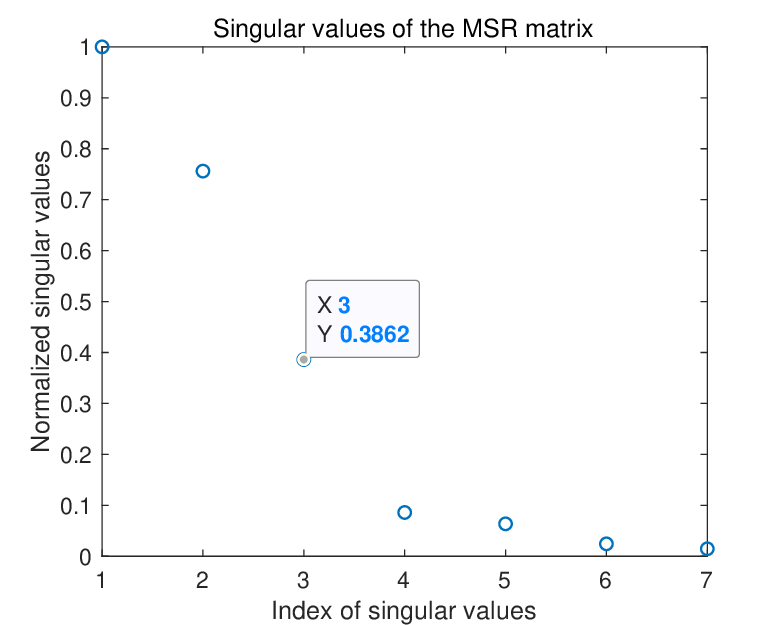}\hfill\includegraphics[width=.250\textwidth]{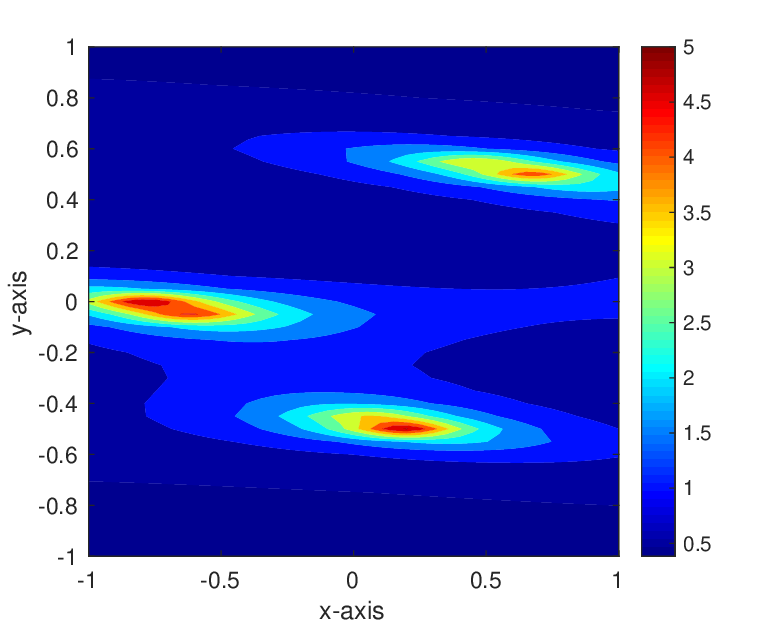}}\hfill
\subfigure[Case 2]{\includegraphics[width=.250\textwidth]{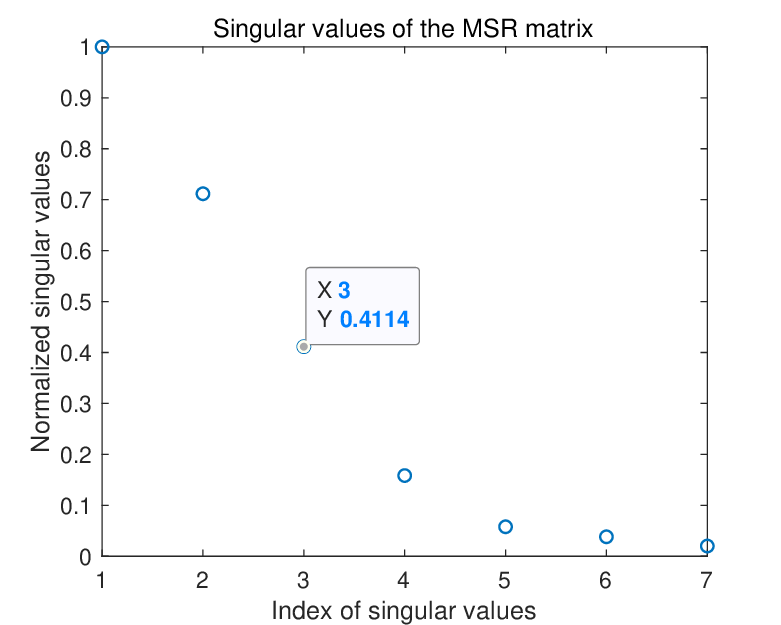}\hfill\includegraphics[width=.250\textwidth]{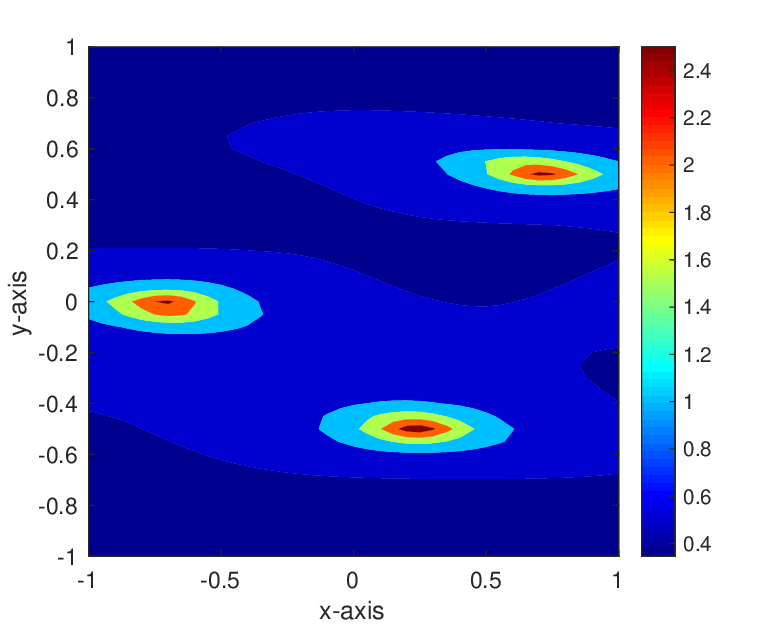}}\\
\subfigure[Case 3]{\includegraphics[width=.250\textwidth]{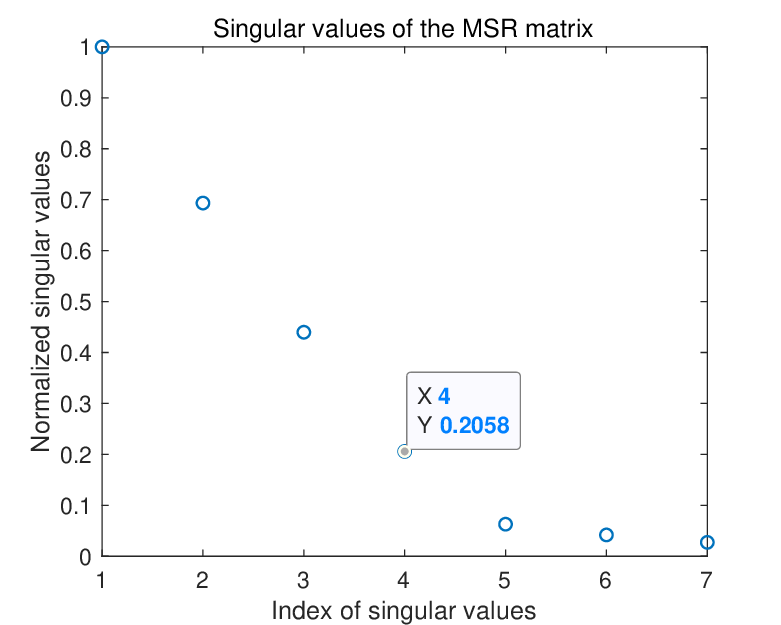}\hfill\includegraphics[width=.250\textwidth]{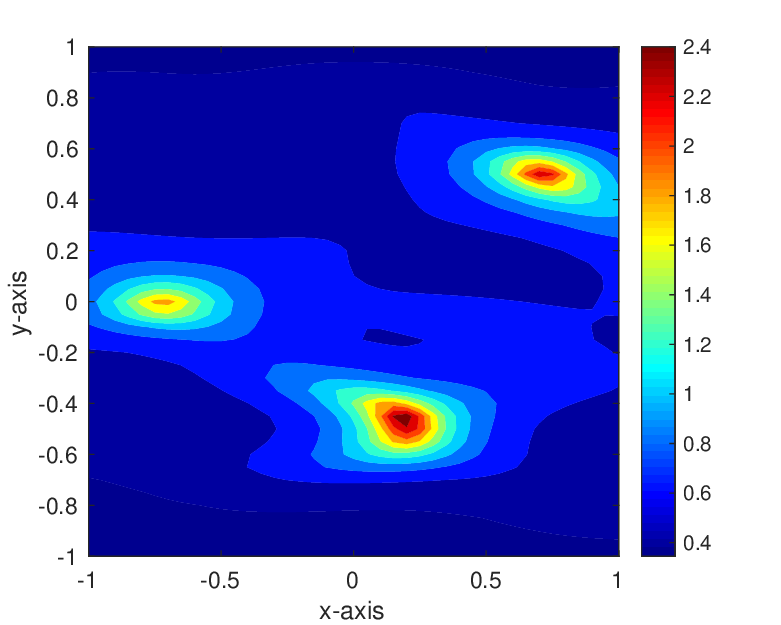}}\hfill
\subfigure[Case 4]{\includegraphics[width=.250\textwidth]{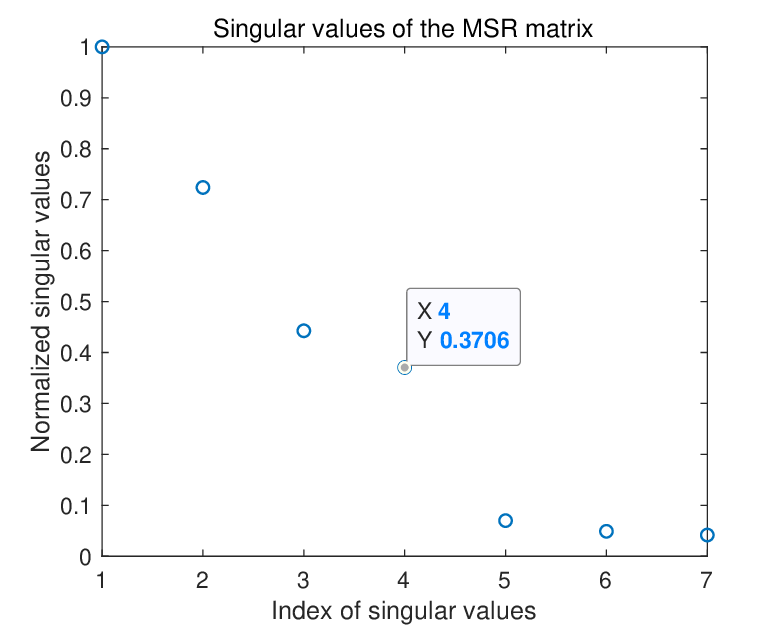}\hfill\includegraphics[width=.250\textwidth]{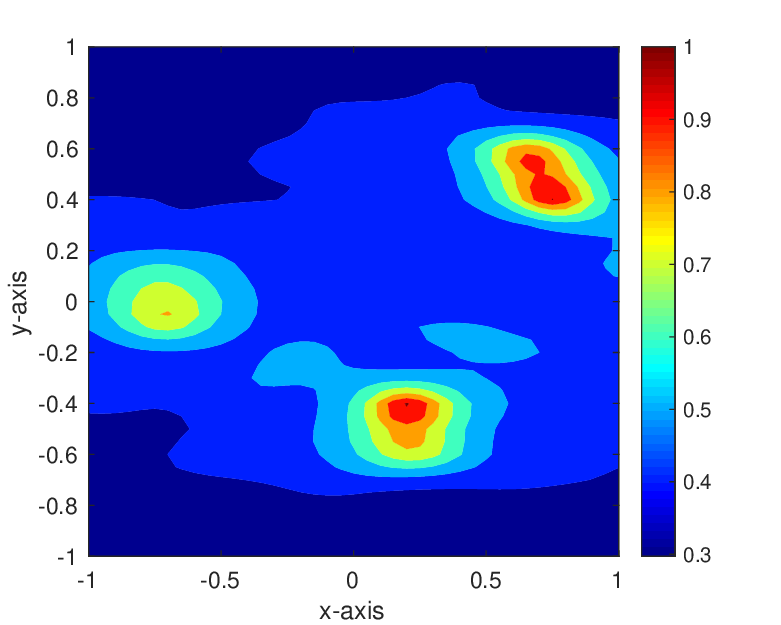}}
\caption{\label{Result3-1}(Example \ref{MU1}) Distribution of singular values (first and third columns) and maps of $f_{\music}(\mr)$ for $k=2\pi/0.4$ (second and fourth columns).}
\end{figure}

\begin{figure}
\centering
\subfigure[Case 5]{\includegraphics[width=.250\textwidth]{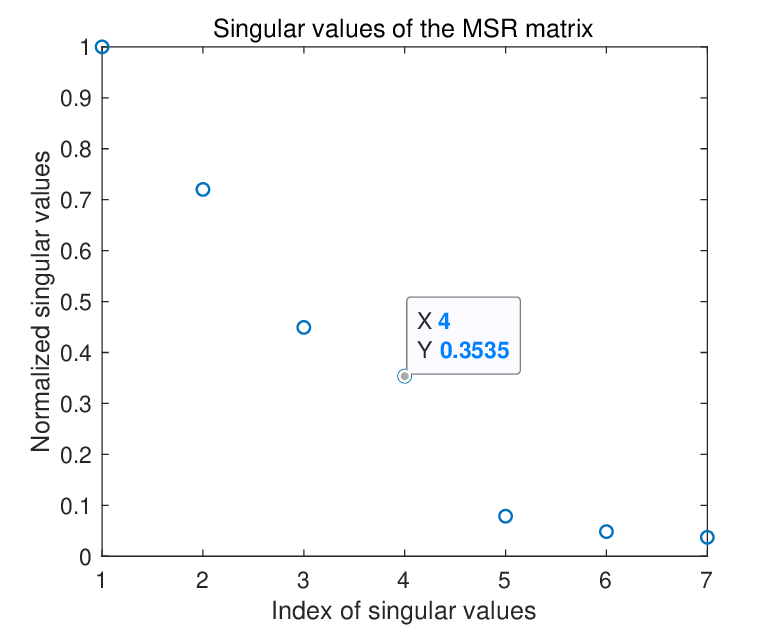}\hfill\includegraphics[width=.250\textwidth]{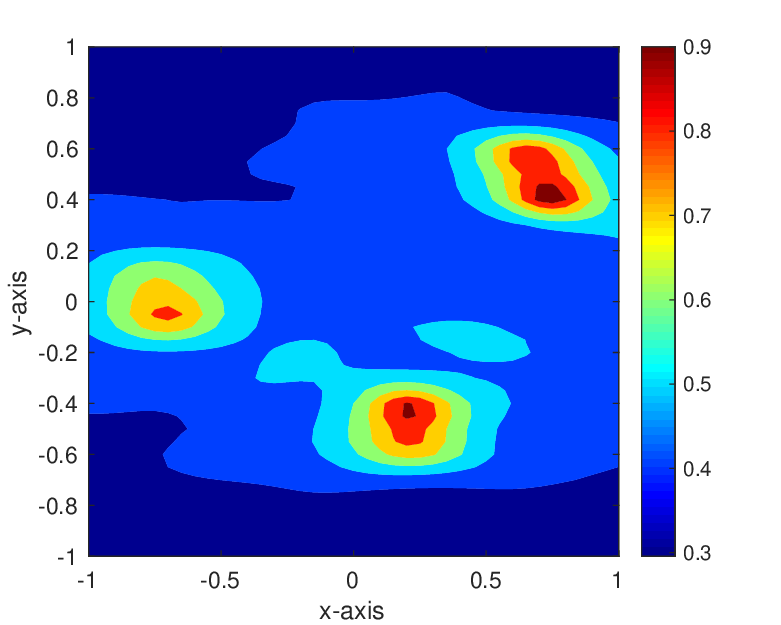}}\hfill
\subfigure[Case 6]{\includegraphics[width=.250\textwidth]{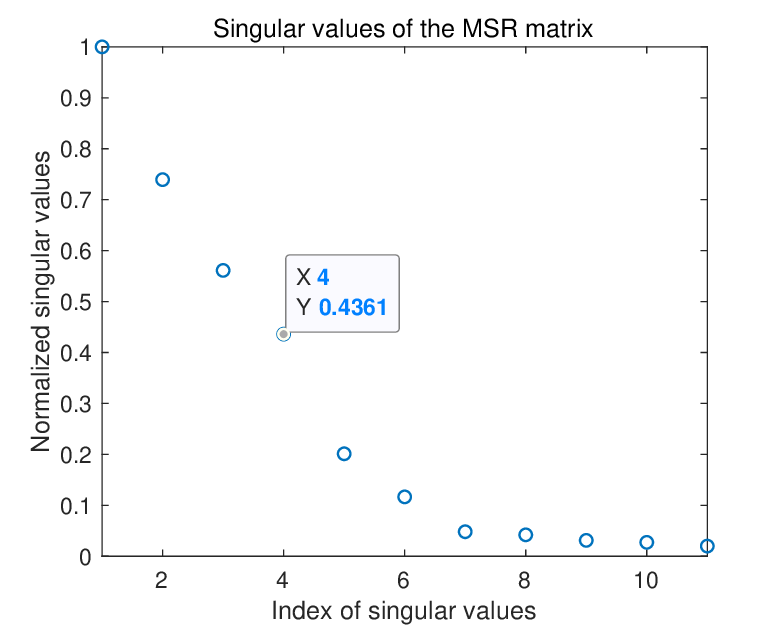}\hfill\includegraphics[width=.250\textwidth]{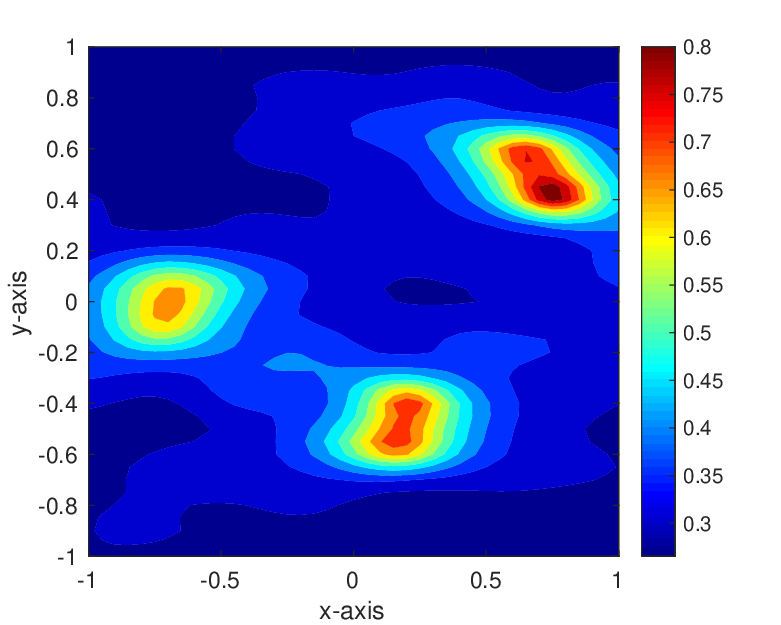}}\\
\subfigure[Case 7]{\includegraphics[width=.250\textwidth]{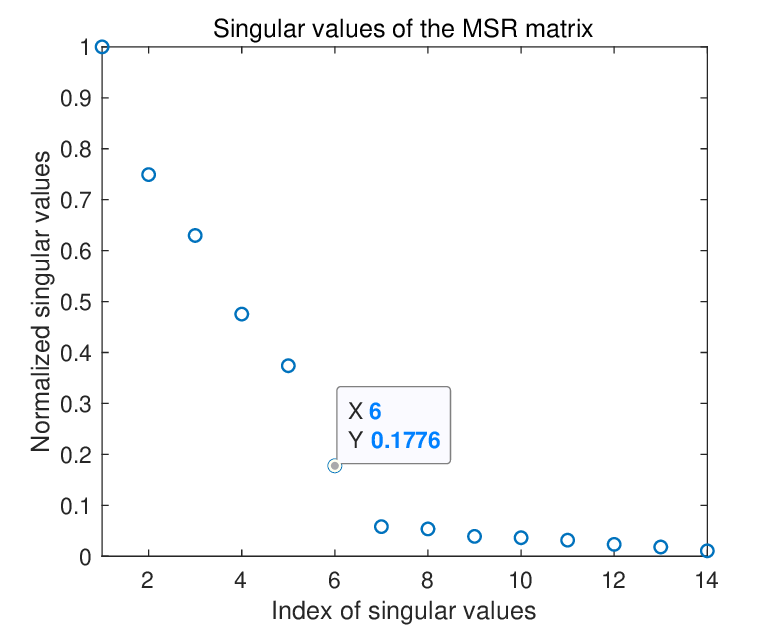}\hfill\includegraphics[width=.250\textwidth]{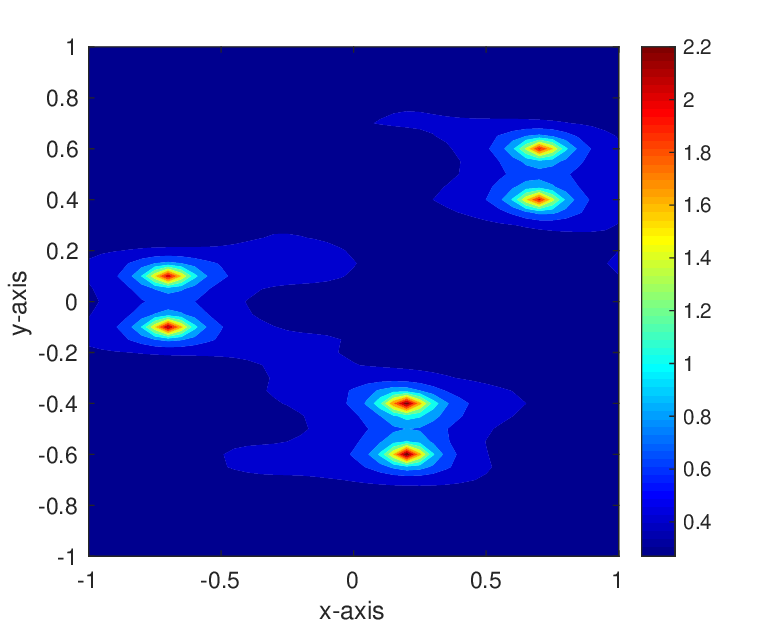}}\hfill
\subfigure[Case 8]{\includegraphics[width=.250\textwidth]{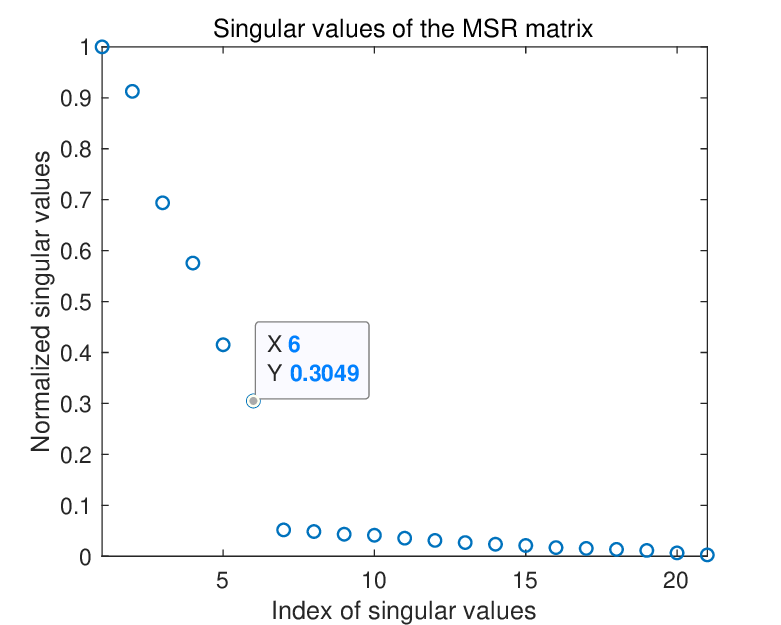}\hfill\includegraphics[width=.250\textwidth]{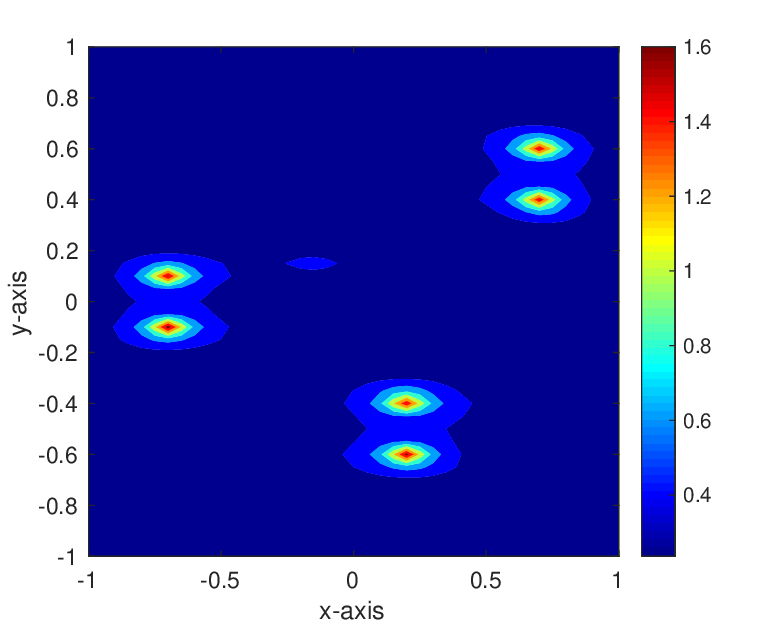}}
\caption{\label{Result3-2}(Example \ref{MU1}) Distribution of singular values (first and third columns) and maps of $f_{\music}(\mr)$ for $k=2\pi/0.4$ (second and fourth columns).}
\end{figure}

\begin{ex}[Permeability contrast case: same permeability]\label{MU1}
Figure \ref{Result3-1} shows the distribution of singular values of $\mathbb{K}$ and the maps of $f_{\music}(\mr)$ for Cases 1-4 when $\eps_s=\epsb$ and $\mu_s=5$. In contrast to the permittivity contrast case, only three or four singular values can be selected to generate the projection operator onto the noise subspace (theoretically, six singular values must be chosen for a proper generation). For Cases 1-3, since the range of incident and observation directions is narrower than $\pi$, the existence of $\Sigma_s$ can be recognized due to the factor $J_0(|k|\mr-\mr_s|)$ from the Discussion \ref{Discussion-M1}, but its exact location cannot be identified. For Case 4, since the range of observation direction is $\pi$, two large-magnitude peaks can be observed in the neighborhood of $\mr_s$ due to the factor $J_1(|k|\mr-\mr_s|)$; refer to Discussion \ref{Discussion-M2}.

Figure \ref{Result3-2} shows the distribution of singular values of $\mathbb{K}$ and the maps of $f_{\music}(\mr)$ for Cases 5-8. Since the range of incident direction is $\pi$, two large-magnitude peaks can be observed in the neighborhood of $\mr_s$. Notice that, when the range of observation direction is close to $\pi$, the two large-magnitude peaks are clearer because the value of $f_{\music}(\mr)$ is dominated by $J_1(|k|\mr-\mr_s|)$; refer to Discussion \ref{Discussion-M1}.
\end{ex}

\begin{ex}[Permeability contrast case: different permeabilities]\label{MU2}
For the final simulation, we consider the imaging results with different permeabilities. For this, we set $\mu_1=5$, $\mu_2=3$, $\mu_3=2$, and $\eps_s=\epsb$. In contrast to the permittivity contrast case in Example \ref{EPS2}, since the far-field pattern data are influenced by the inverse proportion of $\mu_s$, the location of $\mr_3$ can be easily identified, but $\mr_1$ cannot. It is interesting to observe that, for Case 6, three singular values can be selected to generate the projection operator onto the noise subspace, which means that some inhomogeneities (here, $\Sigma_1$) cannot be recognized via MUSIC. Hence, the development of an effective technique for selecting nonzero singular values must be considered at the beginning of the imaging procedure.
\end{ex}

\section{Concluding remarks}\label{sec:5}
In this paper, we designed a MUSIC-type imaging technique to identify the location of small inhomogeneities, the permittivity or permeability of which differs from the background. To explain the applicability of MUSIC in the limited-aperture problem, the mathematical structure of the imaging function was derived by establishing a relationship with an infinite series of first-kind integer-order Bessel functions and the range of incident and observation directions. Based on the explored structure, we examined the various properties of MUSIC and least condition of the range of incident and observation directions to guarantee good results. Various simulation results with noisy data were conducted to validate the theoretical results.

Unfortunately, on the basis of the theoretical and simulation results, it is hard to say that the imaging results via MUSIC do not guarantee complete location identification of the inhomogeneities (especially for the magnetic permeability contrast case). Therefore, improvement of the imaging/detection performance is an interesting research subject that should be examined. 

It has been confirmed that MUSIC is an effective non-iterative imaging technique in real-world microwave imaging; refer to \cite{P-MUSIC6}. Application to real-world microwave imaging with a limited-aperture measurement configuration will be a forthcoming work. {Finally, the design of the imaging algorithm, the analysis of mathematical structure, and the numerical simulations conducted in this paper could be extended to the three-dimensional limited-aperture inverse scattering problem.}

\section*{Acknowledgment}
The author would like to acknowledge two anonymous referees for their comments that help to increase the quality of the paper. This research was supported by the National Research Foundation of Korea (NRF) grant funded by the Korean government (MSIT) (NRF-2020R1A2C1A01005221).

\appendix
\section{Derivation of Theorem \ref{Theorem1}}\label{sec:A}
On the basis of \eqref{DecompositionK1} and \eqref{SVD1}, the following relationship holds
\[\mU_s\approx e^{i\gamma_s^{(1)}}\mf_\eps(\mr_s),\quad\mV_s\approx e^{-i\gamma_s^{(2)}}\overline{\mg}_\eps(\mr_s),\quad\text{and}\quad\gamma_s^{(1)}+\gamma_s^{(2)}=\arg\left(\frac{\alpha^2k^2(1+i)(\eps_s-\epsb)\pi}{4\sqrt{MNk\epsb\mub\pi}}\right),\]
With this, we can evaluate
 \begin{align*}
 \mathbb{P}_{\noise}^{(\eps)}(\mf_\eps(\mr))&=\left(\mathbb{I}(S)-\sum_{s=1}^{S}\mU_s\mU_s^*\right)\mf_\eps(\mr)=\left(\mathbb{I}(S)-\frac{1}{M}\sum_{s=1}^{S}\mf_\eps(\mr_s)\mf_\eps(\mr_s)^*\right)\mf_\eps(\mr)\\
 &=\frac{1}{\sqrt{M}}\left[\begin{array}{c}
  e^{-ik\vv_1\cdot\mr} \\
  e^{-ik\vv_2\cdot\mr} \\
  \vdots \\
  e^{-ik\vv_M\cdot\mr}
 \end{array}
\right]-\frac{1}{M\sqrt{M}}\sum_{s=1}^{S}\left[\begin{array}{c}
\medskip\displaystyle e^{-ik\vv_1\cdot\mr}+\sum_{m\in\mathcal{M}_1}e^{-ik\vv_1\cdot\mr_s}e^{-ik\vv_m\cdot(\mr-\mr_s)}\\
\medskip\displaystyle e^{-ik\vv_2\cdot\mr}+\sum_{m\in\mathcal{M}_2}e^{-ik\vv_2\cdot\mr_s}e^{-ik\vv_m\cdot(\mr-\mr_s)}\\
\medskip\vdots \\
\medskip\displaystyle e^{-ik\vv_M\cdot\mr}+\sum_{m\in\mathcal{M}_M}e^{-ik\vv_M\cdot\mr_s}e^{-ik\vv_m\cdot(\mr-\mr_s)}\\
\end{array}\right],
\end{align*}
where $\mathcal{M}_p=\set{1,2,\cdots,M}\backslash\set{p}$. Since $e^{-ik\vv_p\cdot\mr}=e^{-ik\vv_p\cdot\mr_s}e^{-ik\vv_p\cdot(\mr-\mr_s)}$ and the following relation
\begin{align}
\begin{aligned}\label{Expansion1}
 \frac{1}{M}\sum_{m=1}^{M}e^{-ik\vv_m\cdot\mr}&=\frac{1}{\vartheta_M-\vartheta_1}\int_{\mathbb{S}_{\obs}^1}e^{-ik\vv\cdot\mr}d\vv\\
 &=J_0(k|\mr|)+\frac{4}{\vartheta_M-\vartheta_1}\sum_{p=1}^{\infty}\frac{i^p}{p}J_p(k|\mr|)\sin\frac{p(\vartheta_M-\vartheta_1)}{2}\cos\frac{p(\vartheta_M+\vartheta_1-2\phi+2\pi)}{2},
\end{aligned}
\end{align}
holds uniformly (see \cite{P-SUB3}, for instance), we can derive
\begin{align*}
 &\frac{1}{M}\sum_{s=1}^{S}\left(e^{-ik\vv_p\cdot\mr}+\sum_{p\in\mathcal{M}_p}e^{-ik\vv_p\cdot\mr_s}e^{-ik\vv_p\cdot(\mr-\mr_s)}\right)=\sum_{s=1}^{S}e^{-ik\vv_p\cdot\mr_s}\left(\frac{1}{M}\sum_{m=1}^{M}e^{-ik\vv_m\cdot(\mr-\mr_s)}\right)\\
 &=\sum_{s=1}^{S}e^{-ik\vv_p\cdot\mr_s}\left(J_0(k|\mr-\mr_s|)+\frac{4}{\vartheta_M-\vartheta_1}\sum_{p=1}^{\infty}\frac{i^p}{p}J_p(k|\mr-\mr_s|)\sin\frac{p(\vartheta_M-\vartheta_1)}{2}\cos\frac{p(\vartheta_M+\vartheta_1-2\phi_s+2\pi)}{2}\right).
\end{align*}
Now, let us denote
\begin{equation}\label{Lambda_Eps1}
\Lambda_{\eps}^{(1)}(\mr-\mr_s):=4\sum_{p=1}^{\infty}\frac{i^p}{p}J_p(k|\mr-\mr_s|)\sin\frac{p(\vartheta_M-\vartheta_1)}{2}\cos\frac{p(\vartheta_M+\vartheta_1-2\phi_s+2\pi)}{2}.
\end{equation}
As such, we obtain
\[\mathbb{P}_{\noise}^{(\eps)}(\mf_\eps(\mr))=\frac{1}{\sqrt{M}}
 \left[
 \begin{array}{c}
 \medskip\displaystyle e^{-ik\vv_1\cdot\mr}-\sum_{s=1}^{S}e^{-ik\vv_1\cdot\mr_s}\left(J_0(k|\mr-\mr_s|)+\frac{\Lambda_{\eps}^{(1)}(\mr-\mr_s)}{\vartheta_M-\vartheta_1}\right)\\
 \medskip\displaystyle e^{-ik\vv_2\cdot\mr}-\sum_{s=1}^{S}e^{-ik\vv_2\cdot\mr_s}\left(J_0(k|\mr-\mr_s|)+\frac{\Lambda_{\eps}^{(1)}(\mr-\mr_s)}{\vartheta_M-\vartheta_1}\right)\\
 \medskip\vdots \\
 \displaystyle e^{-ik\vv_M\cdot\mr}-\sum_{s=1}^{S}e^{-ik\vv_M\cdot\mr_s}\left(J_0(k|\mr-\mr_s|)+\frac{\Lambda_{\eps}^{(1)}(\mr-\mr_s)}{\vartheta_M-\vartheta_1}\right)
 \end{array}
 \right].\]
Now, let us denote
\[\Phi_\eps(\mr-\mr_s)=J_0(k|\mr-\mr_s|)+\frac{\Lambda_{\eps}^{(1)}(\mr-\mr_s)}{\vartheta_M-\vartheta_1}\]
by which we obtain
\begin{multline*}
\mathbb{P}_{\noise}^{(\eps)}(\mf_\eps(\mr))\cdot\overline{\mathbb{P}_{\noise}^{(\eps)}(\mf_\eps(\mr))}=\frac{1}{M}\sum_{m=1}^{M}\left(1-\sum_{s=1}^{S}\Big(\Phi_\eps(\mr-\mr_s)e^{-ik\vv_m\cdot(\mr_s-\mr)}+\overline{\Phi_\eps(\mr-\mr_s)}e^{ik\vv_m\cdot(\mr-\mr_s)}\Big)\right.\\
\left.+\sum_{s=1}^{S}\sum_{s'=1}^{S}\Phi_\eps(\mr-\mr_s)\overline{\Phi_\eps(\mr_{s'}-\mr)}e^{-ik\vv_m\cdot(\mr_s-\mr_{s'})}\right).
\end{multline*}
Applying \eqref{Expansion1}, we can evaluate
\begin{align}
\begin{aligned}\label{Term1}
&\frac{1}{M}\sum_{m=1}^{M}\sum_{s=1}^{S}\Phi_\eps(\mr-\mr_s)e^{-ik\vv_m\cdot(\mr-\mr_s)}=\sum_{s=1}^{S}\Phi_\eps(\mr-\mr_s)\overline{\Phi_\eps(\mr-\mr_s)}\\
&\frac{1}{M}\sum_{m=1}^{M}\sum_{s=1}^{S}\overline{\Phi_\eps(\mr-\mr_s)}e^{ik\vv_m\cdot(\mr-\mr_s)}=\sum_{s=1}^{S}\overline{\Phi_\eps(\mr-\mr_s)}\Phi_\eps(\mr-\mr_s).
\end{aligned}
\end{align}

Now, let us assume that $s\ne s'$ and $\mr_{s}-\mr_{s'}=|\mr_{s}-\mr_{s'}|[\cos\xi,\sin\xi]^T$. Then we obtain
\begin{multline*}
\frac{1}{M}\sum_{m=1}^{M}e^{-ik\vv_m\cdot(\mr_{s}-\mr_{s'})}=J_0(k|\mr_{s}-\mr_{s'}|)\\+\frac{4}{\vartheta_M-\vartheta_1}\sum_{p=1}^{\infty}\frac{i^p}{p}J_p(k|\mr_{s}-\mr_{s'}|)\sin\frac{p(\vartheta_M-\vartheta_1)}{2}\cos\frac{p(\vartheta_M+\vartheta_1-2\xi+2\pi)}{2}.
\end{multline*}
Since $k$ satisfies \eqref{Separation}, the following asymptotic form of Bessel function holds
\begin{equation}\label{AsymptoticBessel}
J_0(k|\mr_{s}-\mr_{s'}|)=\sqrt{\frac{2}{k\pi|\mr_{s}-\mr_{s'}|}}\cos\left(k|\mr_{s}-\mr_{s'}|)-\frac{\pi}{4}\right).
\end{equation}
Furthermore, since for $p>0$ and $x\ne0$ (see \cite{L}), it can be expressed as
\begin{equation}\label{UpperboundBessel}
|J_p(x)|\leq\sup\left\{\frac{b}{\sqrt[3]{p}},\frac{c}{\sqrt[3]{|x|}}:b=0.674 885\ldots,~c=0.785747\ldots\right\},
\end{equation}
We can also observe that, for $s\ne s'$,
\begin{equation}\label{Term2}
\sum_{s=1}^{S}\sum_{s'=1}^{S}\Phi_\eps(\mr-\mr_s)\overline{\Phi_\eps(\mr_{s'}-\mr)}e^{-ik\vv_m\cdot(\mr_s-\mr_{s'})}=O\left(\frac{1}{(\vartheta_M-\vartheta_1)^2}\right).
\end{equation}
If $s=s'$, then $e^{-ik\vv_m\cdot(\mr_s-\mr_{s'})}=1$, and, on the basis of \eqref{Term2}, we can examine that
\begin{equation}\label{Term3}
\sum_{s=1}^{S}\sum_{s'=1}^{S}\Phi_\eps(\mr-\mr_s)\overline{\Phi_\eps(\mr_{s'}-\mr)}e^{-ik\vv_m\cdot(\mr_s-\mr_{s'})}=\sum_{s=1}^{S}\Phi_\eps(\mr-\mr_s)\overline{\Phi_\eps(\mr-\mr_s)}=O\left(\frac{1}{(\vartheta_M-\vartheta_1)^2}\right).
\end{equation}
Finally, by combining \eqref{Term1} and \eqref{Term3}, we obtain
\begin{align*}
\mathbb{P}_{\noise}^{(\eps)}(\mf_\eps(\mr))\cdot\overline{\mathbb{P}_{\noise}^{(\eps)}(\mf_\eps(\mr))}&=1-\sum_{s=1}^{S}\Phi_\eps(\mr-\mr_s)\overline{\Phi_\eps(\mr-\mr_s)}+O\left(\frac{1}{(\vartheta_M-\vartheta_1)^2}\right)\\
&=1-\sum_{s=1}^{S}\left|J_0(k|\mr-\mr_s|)+\frac{\Lambda_{\eps}^{(1)}(\mr-\mr_s)}{\vartheta_M-\vartheta_1}\right|^2+O\left(\frac{1}{(\vartheta_M-\vartheta_1)^2}\right).
\end{align*}
Similarly, by letting
\begin{equation}\label{Lambda_Eps2}
\Lambda_{\eps}^{(2)}(\mr-\mr_s):=4\sum_{q=1}^{\infty}\frac{i^q}{q}J_q(k|\mr-\mr_s|)\sin\frac{q(\theta_N-\theta_1)}{2}\cos\frac{q(\theta_N+\theta_1-2\phi_s)}{2},
\end{equation}
we can obtain
\[\mathbb{Q}_{\noise}^{(\eps)}(\mg_\eps(\mr))\cdot\overline{\mathbb{Q}_{\noise}^{(\eps)}(\mg_\eps(\mr))}=1-\sum_{s=1}^{S}\left|J_0(k|\mr-\mr_s|)+\frac{\Lambda_{\eps}^{(2)}(\mr-\mr_s)}{\theta_N-\theta_1}\right|^2+O\left(\frac{1}{(\theta_N-\theta_1)^2}\right).\]
Hence, \eqref{Structure1} is derived. This completes the proof.

\section{Derivation of Theorem \ref{Theorem2}}\label{sec:B}
On the basis of \eqref{DecompositionK2} and \eqref{SVD2}, the following relationship holds
\begin{align*}
\mU_{2s-1}\approx e^{i\gamma_{s}^{(1,1)}}\mH_{s}^{(1)}(\mr_s),\quad\mU_{2s}\approx e^{i\gamma_{s}^{(1,2)}}\mH_{s}^{(2)}(\mr_s),\quad\mV_{2s-1}\approx e^{-i\gamma_{s}^{(2,1)}}\mG_{s}^{(1)}(\mr_s),\quad\mV_{2s}\approx e^{-i\gamma_{s}^{(2,2)}}\mG_{s}^{(2)}(\mr_s)
\end{align*}
and
\[\gamma_{s}^{(1,1)}+\gamma_{s}^{(2,1)}=\gamma_{s}^{(1,2)}+\gamma_{s}^{(2,2)}=\arg\left(\frac{\alpha^2k^2(1+i)\mub\pi}{2(\mu_s+\mub)\sqrt{MNk\pi}}\right).\]
With this, we can evaluate
 \begin{align*}
 \mathbb{P}_{\noise}^{(\mu)}(\mf_\eps(\mr))&=\left(\mathbb{I}(S)-\sum_{s=1}^{2S}\mU_s\mU_s^*\right)\mf_\eps(\mr)=\left(\mathbb{I}(S)-\frac{1}{C_1}\sum_{s=1}^{S}\Big(\mH_{s}^{(1)}(\mr_s)\mH_{s}^{(1)}(\mr_s)^*+\mH_{s}^{(2)}(\mr_s)\mH_{s}^{(2)}(\mr_s)^*\Big)\right)\mf_\eps(\mr)\\
 &=\frac{1}{\sqrt{M}}\left[\begin{array}{c}
  e^{-ik\vv_1\cdot\mr} \\
  e^{-ik\vv_2\cdot\mr} \\
  \vdots \\
  e^{-ik\vv_M\cdot\mr}
 \end{array}
\right]-\frac{1}{C_1\sqrt{M}}\sum_{s=1}^{S}\left[\begin{array}{c}
K_1\\K_2\\\vdots\\K_M
 \end{array}\right],
\end{align*}
where
\begin{align*}
K_m=&(-\vv_m\cdot\me_1)^2e^{-ik\vv_m\cdot\mr}+(-\vv_m\cdot\me_1)e^{-ik\vv_m\cdot\mr_s}\sum_{p\in\mathcal{M}_m}(-\vv_p\cdot\me_1)e^{-ik\vv_p\cdot(\mr-\mr_s)}\\
&+(-\vv_m\cdot\me_2)^2e^{-ik\vv_m\cdot\mr}+(-\vv_m\cdot\me_2)e^{-ik\vv_m\cdot\mr_s}\sum_{p\in\mathcal{M}_m}(-\vv_p\cdot\me_2)e^{-ik\vv_p\cdot(\mr-\mr_s)}
\end{align*}
for $m=1,2,\cdots,M$. Since $e^{-ik\vv_m\cdot\mr}=e^{-ik\vv_m\cdot\mr_s}e^{-ik\vv_m\cdot(\mr-\mr_s)}$, the following relations hold uniformly (see \cite{P-SUB3} for instance)
\begin{align}
\begin{aligned}\label{Expansion2-1}
 \int_{\mathbb{S}_{\obs}^1}(-\vv\cdot\me_1)&e^{-ik\vv\cdot\mr}d\vv=2J_0(k|\mr|)\sin\frac{\vartheta_M-\vartheta_1}{2}\cos\frac{\vartheta_M+\vartheta_1}{2}\\
 &+iJ_1(k|\mr|)\bigg[(\vartheta_M-\vartheta_1)\bigg(\frac{\mr}{|\mr|}\cdot\me_1\bigg)+\sin(\vartheta_M-\vartheta_1)\cos(\vartheta_M+\vartheta_1-\phi)\bigg]\\
 &+2\sum_{p=2}^{\infty}i^pJ_p(k|\mr|)\bigg[\frac{1}{1-p}\sin\frac{(1-p)(\vartheta_M-\vartheta_1)}{2}\cos\frac{(1-p)(\vartheta_M+\vartheta_1)+2p\phi-2p\pi}{2}\\
 &+\frac{1}{1+p}\sin\frac{(1+n)(\vartheta_M-\vartheta_1)}{2}\cos\frac{(1+p)(\vartheta_M+\vartheta_1)-2p\phi+2p\pi}{2}\bigg]
 \end{aligned}
 \end{align}
 and
 \begin{align}
\begin{aligned}\label{Expansion2-2}
 \int_{\mathbb{S}_{\obs}^1}(-\vv\cdot\me_2)&e^{-ik\vv\cdot\mr}d\vv=2J_0(k|\mr|)\sin\frac{\vartheta_M-\vartheta_1}{2}\cos\frac{\vartheta_M+\vartheta_1-\pi}{2}\\
 &+iJ_1(k|\mr|)\bigg[(\vartheta_M-\vartheta_1)\bigg(\frac{\mr}{|\mr|}\cdot\me_2\bigg)+\sin(\vartheta_M-\vartheta_1)\sin(\vartheta_M+\vartheta_1-\phi)\bigg]\\
 &+2\sum_{p=2}^{\infty}i^pJ_p(k|\mr|)\bigg[\frac{1}{1-p}\sin\frac{(1-p)(\vartheta_M-\vartheta_1)}{2}\cos\frac{(1-p)(\vartheta_M+\vartheta_1)+2p\phi-(2p+1)\pi}{2}\\
 &+\frac{1}{1+p}\sin\frac{(1+n)(\vartheta_M-\vartheta_1)}{2}\cos\frac{(1+p)(\vartheta_M+\vartheta_1)-2p\phi+(2p-1)\pi}{2}\bigg],
 \end{aligned}
 \end{align}
we can derive
\begin{align*}
\frac{1}{C_1}K_m=&(-\vv_m\cdot\me_1)e^{-ik\vv_m\cdot\mr_s}\frac{1}{C_1}\sum_{m'=1}^{M}(\vv_{m'}\cdot\me_1)e^{-ik\vv_{m'}\cdot(\mr-\mr_s)}\\
&+(-\vv_m\cdot\me_2)e^{-ik\vv_m\cdot\mr_s}\frac{1}{C_1}\sum_{m'=1}^{M}(\vv_{m'}\cdot\me_2)e^{-ik\vv_{m'}\cdot(\mr-\mr_s)}\\
=&(-\vv_m\cdot\me_1)e^{-ik\vv_m\cdot\mr_s}\left\{iJ_1(k|\mr-\mr_s|)\bigg(\frac{\mr-\mr_s}{|\mr-\mr_s|}\cdot\me_1\bigg)+\frac{\Lambda_{\mu}^{(1,1)}(\mr-\mr_s)}{\vartheta_M-\vartheta_1}\right\}\\
&+(-\vv_m\cdot\me_2)e^{-ik\vv_m\cdot\mr_s}\left\{iJ_1(k|\mr-\mr_s|)\bigg(\frac{\mr-\mr_s}{|\mr-\mr_s|}\cdot\me_2\bigg)+\frac{\Lambda_{\mu}^{(1,2)}(\mr-\mr_s)}{\vartheta_M-\vartheta_1}\right\},
\end{align*}
where
\begin{align}
\begin{aligned}\label{Lambda_Mu1-1}
\Lambda_{\mu}^{(1,1)}(\mr-\mr_s):=&2J_0(k|\mr-\mr_s|)\sin\frac{\vartheta_M-\vartheta_1}{2}\cos\frac{\vartheta_M+\vartheta_1}{2}+iJ_1(k|\mr-\mr_s|)\sin(\vartheta_M-\vartheta_1)\cos(\vartheta_M+\vartheta_1-\phi_s)\\
&+2\sum_{p=2}^{\infty}i^pJ_p(k|\mr-\mr_s|)\bigg[\frac{1}{1-p}\sin\frac{(1-p)(\vartheta_M-\vartheta_1)}{2}\cos\frac{(1-p)(\vartheta_M+\vartheta_1)+2p\phi_s-2p\pi}{2}\\
 &+\frac{1}{1+p}\sin\frac{(1+n)(\vartheta_M-\vartheta_1)}{2}\cos\frac{(1+p)(\vartheta_M+\vartheta_1)-2p\phi_s+2p\pi}{2}\bigg]
 \end{aligned}
\end{align}
and
\begin{align}
\begin{aligned}\label{Lambda_Mu1-2}
\Lambda_{\mu}^{(1,2)}(\mr-\mr_s):=&2J_0(k|\mr-\mr_s|)\sin\frac{\vartheta_M-\vartheta_1}{2}\cos\frac{\vartheta_M+\vartheta_1-\phi_s}{2}
 +iJ_1(k|\mr-\mr_s|)\sin(\vartheta_M-\vartheta_1)\sin(\vartheta_M+\vartheta_1-\phi_s)\\
 &+2\sum_{p=2}^{\infty}i^pJ_p(k|\mr-\mr_s|)\bigg[\frac{1}{1-p}\sin\frac{(1-p)(\vartheta_M-\vartheta_1)}{2}\cos\frac{(1-p)(\vartheta_M+\vartheta_1)+2p\phi_s-(2p+1)\pi}{2}\\
 &+\frac{1}{1+p}\sin\frac{(1+n)(\vartheta_M-\vartheta_1)}{2}\cos\frac{(1+p)(\vartheta_M+\vartheta_1)-2p\phi_s+(2p-1)\pi}{2}\bigg]
\end{aligned}
\end{align}

With this, we can express $\mathbb{P}_{\noise}^{(\mu)}(\mf_\eps(\mr))$ as
 \begin{align*}
 \mathbb{P}_{\noise}^{(\mu)}(\mf_\eps(\mr))&=\frac{1}{\sqrt{M}}\left[\begin{array}{c}
  \displaystyle e^{-ik\vv_1\cdot\mr}-\sum_{s=1}^{S}\sum_{h=1}^{2}(-\vv_m\cdot\me_h)e^{-ik\vv_m\cdot\mr_s}\left\{iJ_1(k|\mr-\mr_s|)\bigg(\frac{\mr-\mr_s}{|\mr-\mr_s|}\cdot\me_h\bigg)+\frac{\Theta_h}{\vartheta_M-\vartheta_1}\right\}\\
  \displaystyle e^{-ik\vv_2\cdot\mr}-\sum_{s=1}^{S}\sum_{h=1}^{2}(-\vv_m\cdot\me_h)e^{-ik\vv_m\cdot\mr_s}\left\{iJ_1(k|\mr-\mr_s|)\bigg(\frac{\mr-\mr_s}{|\mr-\mr_s|}\cdot\me_h\bigg)+\frac{\Theta_h}{\vartheta_M-\vartheta_1}\right\}\\
  \vdots \\
  \displaystyle e^{-ik\vv_M\cdot\mr}-\sum_{s=1}^{S}\sum_{h=1}^{2}(-\vv_m\cdot\me_h)e^{-ik\vv_m\cdot\mr_s}\left\{iJ_1(k|\mr-\mr_s|)\bigg(\frac{\mr-\mr_s}{|\mr-\mr_s|}\cdot\me_h\bigg)+\frac{\Theta_h}{\vartheta_M-\vartheta_1}\right\}\\
 \end{array}
\right].
\end{align*}
Now, let us denote
\[\Phi_\mu(\mr-\mr_s,h)=iJ_1(k|\mr-\mr_s|)\bigg(\frac{\mr-\mr_s}{|\mr-\mr_s|}\cdot\me_h\bigg)+\frac{\Theta_h}{\vartheta_M-\vartheta_1}\]
then
\begin{align*}
&\mathbb{P}_{\noise}^{(\mu)}(\mf_\eps(\mr))\cdot\overline{\mathbb{P}_{\noise}^{(\mu)}(\mf_\eps(\mr))}\\
&=\frac{1}{M}\sum_{m=1}^{M}\bigg\{1-\sum_{s=1}^{S}\sum_{h=1}^{2}(-\vv_m\cdot\me_h)e^{-ik\vv_m\cdot(\mr-\mr_s)}\overline{\Phi_\mu(\mr-\mr_s,h)}-\sum_{s=1}^{S}\sum_{h=1}^{2}(-\vv_m\cdot\me_h)e^{ik\vv_m\cdot(\mr-\mr_s)}\Phi_\mu(\mr-\mr_s,h)\\
&+\bigg(\sum_{s=1}^{S}\sum_{h=1}^{2}(-\vv_m\cdot\me_h)e^{-ik\vv_m\cdot\mr_s}\Phi_\mu(\mr-\mr_s,h)\bigg)\bigg(\sum_{s'=1}^{S}\sum_{h''=1}^{2}(-\vv_m\cdot\me_{h''})e^{ik\vv_m\cdot\mr_{s'}}\overline{\Phi_\mu(\mr-\mr_s,h'')}\bigg)\bigg\}.
\end{align*}
Applying \eqref{Expansion2-1} and \eqref{Expansion2-2}, we obtain
\begin{multline}\label{Term4-1}
\frac{1}{M}\sum_{m=1}^{M}\sum_{s=1}^{S}\sum_{h=1}^{2}(-\vv_m\cdot\me_h)e^{-ik\vv_m\cdot(\mr-\mr_s)}\overline{\Phi_\mu(\mr-\mr_s,h)}\\
\approx\frac{1}{C_1}\sum_{s=1}^{S}\sum_{h=1}^{2}\int_{\mathbb{S}_{\obs}^1}(-\vv\cdot\me_h)e^{-ik\vv\cdot(\mr-\mr_s)}\overline{\Phi_\mu(\mr-\mr_s,h)}d\vv=\sum_{s=1}^{S}\sum_{h=1}^{2}\Phi_\mu(\mr-\mr_s,h)\overline{\Phi_\mu(\mr-\mr_s,h)}
\end{multline}
and
\begin{multline}\label{Term4-2}
\frac{1}{M}\sum_{m=1}^{M}\sum_{s=1}^{S}\sum_{h=1}^{2}(-\vv_m\cdot\me_h)e^{ik\vv_m\cdot(\mr-\mr_s)}\Phi_\mu(\mr-\mr_s,h)\\
\approx\frac{1}{C_1}\sum_{s=1}^{S}\sum_{h=1}^{2}\int_{\mathbb{S}_{\obs}^1}(-\vv\cdot\me_h)e^{ik\vv\cdot(\mr-\mr_s)}\Phi_\mu(\mr-\mr_s,h)d\vv=\sum_{s=1}^{S}\sum_{h=1}^{2}\overline{\Phi_\mu(\mr-\mr_s,h)}\Phi_\mu(\mr-\mr_s,h).
\end{multline}

Let $s\ne s'$ and $\mr_s-\mr_{s'}=|\mr_s-\mr_{s'}|[\cos\xi,\sin\xi]^T$. Since $k$ satisfies \eqref{Separation}, based on \eqref{AsymptoticBessel} and \eqref{UpperboundBessel}, we observe that
\begin{align*}
\frac{1}{M}\sum_{m=1}^{M}(-\vv_m\cdot\me_h)e^{-ik\vv_m\cdot(\mr_s-\mr_{s'})}&\approx\frac{1}{C_1}\int_{\mathbb{S}_{\obs}^1}(-\vv\cdot\me_h)e^{-ik\vv\cdot(\mr_s-\mr_{s'})}d\vv\\
&=iJ_1(k|\mr_s-\mr_{s'}|)\bigg(\frac{\mr_s-\mr_{s'}}{|\mr_s-\mr_{s'}|}\cdot\me_h\bigg)+\frac{\Theta_h}{\vartheta_M-\vartheta_1}=O\left(\frac{1}{(\vartheta_M-\vartheta_1)^2}\right).
\end{align*}
Let $s=s'$. Then, based on the orthonormality of singular vectors, we obtain
\begin{align}
\begin{aligned}\label{Term5}
\frac{1}{M}&\sum_{m=1}^{M}\bigg(\sum_{s=1}^{S}\sum_{h=1}^{2}(-\vv_m\cdot\me_h)e^{-ik\vv_m\cdot\mr_s}\Phi_\mu(\mr-\mr_s,h)\bigg)\bigg(\sum_{s'=1}^{S}\sum_{h''=1}^{2}(-\vv_m\cdot\me_{h''})e^{ik\vv_m\cdot\mr_{s'}}\overline{\Phi_\mu(\mr-\mr_s,h'')}\bigg)\\
=&\frac{1}{M}\sum_{s=1}^{S}\sum_{m=1}^{M}\bigg(\sum_{h=1}^{2}(-\vv_m\cdot\me_h)e^{-ik\vv_m\cdot\mr_s}\Phi_\mu(\mr-\mr_s,h)\bigg)\bigg(\sum_{h''=1}^{2}(-\vv_m\cdot\me_{h''})e^{ik\vv_m\cdot\mr_s}\overline{\Phi_\mu(\mr-\mr_s,h'')}\bigg)\\
\approx&\sum_{s=1}^{S}\sum_{h=1}^{2}\Phi_\mu(\mr-\mr_s,h)\overline{\Phi_\mu(\mr-\mr_s,h'')}\bigg(\frac{1}{C_1}\int_{\mathbb{S}_{\obs}^1}(-\vv\cdot\me_h)^2d\vv\bigg)
=\sum_{s=1}^{S}\sum_{h=1}^{2}\Phi_\mu(\mr-\mr_s,h)\overline{\Phi_\mu(\mr-\mr_s,h'')}.
\end{aligned}
\end{align}
Now, by combining \eqref{Term4-1}, \eqref{Term4-2}, and \eqref{Term5}, we obtain
\begin{align*}
\mathbb{P}_{\noise}^{(\mu)}(\mf_\eps(\mr))\cdot\overline{\mathbb{P}_{\noise}^{(\mu)}(\mf_\eps(\mr))}&=1-\sum_{s=1}^{S}\sum_{h=1}^{2}\Phi_\mu(\mr-\mr_s,h)\overline{\Phi_\mu(\mr-\mr_s,h)}+O\left(\frac{1}{(\vartheta_M-\vartheta_1)^2}\right)\\
&=1-\sum_{s=1}^{S}\sum_{h=1}^{2}\left|iJ_1(k|\mr-\mr_s|)\bigg(\frac{\mr-\mr_s}{|\mr-\mr_s|}\cdot\me_h\bigg)+\frac{\Theta_h}{\vartheta_M-\vartheta_1}\right|^2+O\left(\frac{1}{(\vartheta_M-\vartheta_1)^2}\right).
\end{align*}

Now, let
\begin{align}
\begin{aligned}\label{Lambda_Mu2-1}
 \Lambda_{\mu}^{(2,1)}(\mr-\mr_s):=&2J_0(k|\mr-\mr_s|)\sin\frac{\theta_N-\theta_1}{2}\cos\frac{\theta_N+\theta_1}{2}+iJ_1(k|\mr-\mr_s|)\sin(\theta_N-\theta_1)\cos(\theta_N+\theta_1-\phi_s)\\
&+2\sum_{p=2}^{\infty}i^pJ_p(k|\mr-\mr_s|)\bigg[\frac{1}{1-p}\sin\frac{(1-p)(\theta_N-\theta_1)}{2}\cos\frac{(1-p)(\theta_N+\theta_1)+2p\phi_s-2p\pi}{2}\\
 &+\frac{1}{1+p}\sin\frac{(1+n)(\theta_N-\theta_1)}{2}\cos\frac{(1+p)(\theta_N+\theta_1)-2p\phi_s+2p\pi}{2}\bigg], \end{aligned}
 \end{align}
 and
 \begin{align}
\begin{aligned}\label{Lambda_Mu2-2}
\Lambda_{\mu}^{(2,2)}(\mr-\mr_s):=&2J_0(k|\mr-\mr_s|)\sin\frac{\theta_N-\theta_1}{2}\cos\frac{\theta_N+\theta_1-\phi_s}{2}
 +iJ_1(k|\mr-\mr_s|)\sin(\theta_N-\theta_1)\sin(\theta_N+\theta_1-\phi_s)\\
 &+2\sum_{p=2}^{\infty}i^pJ_p(k|\mr-\mr_s|)\bigg[\frac{1}{1-p}\sin\frac{(1-p)(\theta_N-\theta_1)}{2}\cos\frac{(1-p)(\theta_N+\theta_1)+2p\phi_s-(2p+1)\pi}{2}\\
 &+\frac{1}{1+p}\sin\frac{(1+n)(\theta_N-\theta_1)}{2}\cos\frac{(1+p)(\theta_N+\theta_1)-2p\phi_s+(2p-1)\pi}{2}\bigg].
 \end{aligned}
 \end{align}
Then, with a similar process, we can derive
\[\mathbb{Q}_{\noise}^{(\mu)}(\mf_\eps(\mr))\cdot\overline{\mathbb{Q}_{\noise}^{(\mu)}(\mf_\eps(\mr))}=1-\sum_{s=1}^{S}\sum_{h=1}^{2}\left|iJ_1(k|\mr-\mr_s|)\bigg(\frac{\mr-\mr_s}{|\mr-\mr_s|}\cdot\me_h\bigg)+\frac{\Lambda_{\mu}^{(2,h)}(\mr-\mr_s)}{\theta_N-\theta_1}\right|^2+O\left(\frac{1}{(\theta_N-\theta_1)^2}\right).\] 
Hence, \eqref{Structure2} is derived.

\end{document}